\documentclass[a4paper,reqno]{amsart} 
\usepackage{amssymb,epic,epsfig,changebar}
\usepackage[all]{xy}
\SelectTips{cm}{10}
\renewcommand{\labelenumi}{(\arabic{enumi})}

\renewcommand{\P}{{\mathbb P}}
\newcommand{\Z}{{\mathbb Z}}
\newcommand{\C}{{\mathbb C}} 
\newcommand{\T}{{\mathbb T}} 
\newcommand{\kj}{\mathcal{J}}
\newcommand{\kc}{\mathcal{C}} 
\newcommand{\ko}{\mathcal{O}} 

\newcommand{\kd}{\mathcal{D}}

\newcommand{\kp}{\mathcal{P}}
\newcommand{\fm}{{\mathfrak m}}

\newcommand{\fp}{{\mathfrak p}}
\newcommand{\lra}{\longrightarrow} 
\newcommand{\ep}{\textit{ep}}
\newcommand{\wes}{\text{wes}}
\newcommand{\ses}{\text{ses}}
\newcommand{\es}{\textit{es}}

\newcommand{\fix}{\textit{fix}}
\newcommand{\mt}{\text{mt}}

\newcommand{\ef}{\text{ef}}
\newcommand{\Keps}{K[\varepsilon]}
\newcommand{\Rbar}{\overline{R}}

\newcommand{\Hbar}{\overline{H}}
\newcommand{\fmbar}{\overline{\fm}}
\newcommand{\osigma}{\overline{\sigma}}
\newcommand{\otau}{\overline{\tau}}
\newcommand{\oXi}{\overline{\Xi}}

\DeclareMathAlphabet{\mathsc}{U}{rsfs}{m}{n}
\newcommand{\sA}{\mathsc{A}}

\newcommand{\sC}{\mathsc{C}}

\newcommand{\sG}{\mathsc{G}}

\newcommand{\sR}{\mathsc{R}}

\newcommand{\kT}{\mathcal{T}}

\DeclareMathOperator{\con}{con}
\DeclareMathOperator{\Hom}{Hom}
\DeclareMathOperator{\Spec}{Spec}
\DeclareMathOperator{\Quot}{Quot}
\DeclareMathOperator{\Proj}{Proj}
\DeclareMathOperator{\Bl}{B\ell}
\DeclareMathOperator{\ES}{ES}
\DeclareMathOperator{\Ess}{Ess}
\DeclareMathOperator{\gr}{gr}
\DeclareMathOperator{\Gr}{Gr}
\DeclareMathOperator{\Ker}{ker}
\DeclareMathOperator{\Der}{Der}
\DeclareMathOperator{\mult}{mult}

\DeclareMathOperator{\Char}{char}
\DeclareMathOperator{\ord}{ord}
\DeclareMathOperator{\id}{id}

\DeclareMathOperator{\Def}{\kd\!\!\:\it{ef}}

\newcommand{\uDef}{\underline{\Def}\!\,}
\newcommand{\Defem}{\Def^{\text{\it em}}}
\newcommand{\Defes}{\Def^{\text{\it es}}}
\newcommand{\Defessec}{\Def^{\text{\it es,sec}}}
\newcommand{\uDefessec}{\uDef^{\text{\it es,sec}}}
\newcommand{\Defep}{\Def^{\text{\it ep}}}

\newcommand{\Defwes}{\Def^{\text{\it wes}}}
\newcommand{\uDefem}{\uDef^{\text{\it em}}}
\newcommand{\uDefes}{\uDef^{\text{\it es}}}
\newcommand{\uDefep}{\uDef^{\text{\it ep}}}

\newcommand{\Defsec}{\Def^{\text{\it sec}}}
\newcommand{\uDefsec}{\uDef^{\text{\it sec}}}
\renewcommand{\sec}{\textit{sec}}

\newcommand{\Ties}{T^{1,\text{\it es}}}
\newcommand{\Tiesec}{T^{1,\text{\it es, sec}}}

\newcommand{\Tiwes}{T^{1,\text{\it wes}}}

\newcommand{\Msec}{M^{\text{\it sec}}}
\newcommand{\Tisec}{T^{1,\text{\it sec}}}
\newcommand{\Tnsec}{T^{0,\text{\it sec}}}
\newcommand{\Bes}{B^{\text{\it es}}}
\newcommand{\Ies}{I^{\text{\it es}}}

\newcommand{\Iep}{I^{\text{\it ep}}}
\newcommand{\Tep}{T^{\text{\it ep}}}
\newcommand{\bTep}{\bT^{\text{\it ep}}}

\newcommand{\Iesf}{I^{\text{\it es}}_{\fix}}

\newcommand{\ba}{\boldsymbol{a}}

\newcommand{\bb}{\boldsymbol{b}}
\newcommand{\bbeta}{\boldsymbol{\beta}}

\newcommand{\bd}{\boldsymbol{d}}

\newcommand{\bm}{\boldsymbol{m}}

\newcommand{\bnull}{\boldsymbol{0}}
\newcommand{\bone}{\boldsymbol{1}}

\newcommand{\bord}{\boldsymbol{\ord}}

\newcommand{\bR}{\boldsymbol{R}}

\newcommand{\bs}{\boldsymbol{s}}
\newcommand{\bt}{\boldsymbol{t}}
\newcommand{\bT}{\boldsymbol{T}}

\newcommand{\bw}{\boldsymbol{w}}
\newcommand{\bx}{\boldsymbol{x}}
\newcommand{\bW}{\boldsymbol{W}}

\newcommand{\by}{\boldsymbol{y}}

\newcommand{\bz}{\boldsymbol{z}}
\numberwithin{equation}{section}

\newenvironment{Enumerate}
{\begin{list}{\labelenumi}{\usecounter{enumi}%
\setlength{\itemsep}{4pt}%
\setlength{\labelsep}{0pt}\setlength{\leftmargin}{0pt}%
\setlength{\labelwidth}{0pt}\setlength{\listparindent}{2pt}%
}}%
{\end{list}}

\begin{document}

\newtheorem{lemma}{Lemma}[section]
\newtheorem{proposition}[lemma]{Proposition}
\newtheorem{theorem}[lemma]{Theorem}
\newtheorem{corollary}[lemma]{Corollary}
\newtheorem{satz}[lemma]{Satz}

\theoremstyle{definition}
\newtheorem{definition}[lemma]{Definition}
\newtheorem{algorithm}[lemma]{Algorithm}
\newtheorem{example}[lemma]{Example}
\newtheorem{examples}[lemma]{Examples}
\newtheorem{remarks}[lemma]{Remark}
\newtheorem{remark}[lemma]{Remark}

\title[Equisingular Deformations in Arbitrary Characteristic]{Equisingular
  Deformations of Plane Curves in Arbitrary Characteristic}
\author{Antonio Campillo}
\address{Departamento de Algebra, Geometria y Topologia\\
Universidad de Valladolid\\Facultad de Ciencias\\
E -- 47005   Valladolid}
\author{Gert-Martin Greuel}
\address{
Fachbereich Mathematik\\
TU Kaiserslautern\\
Erwin-Schr\"odinger-Stra{\ss}e\\
D -- 67663 Kaiserslautern}
\author{Christoph Lossen}
\address{
Fachbereich Mathematik\\
TU Kaiserslautern\\
Erwin-Schr\"odinger-Stra{\ss}e\\
D -- 67663 Kaiserslautern}

\vskip2cm

%\date{}
\maketitle
\vskip1cm

\centerline{Dedicated to Joseph Steenbrink on the occasion of his sixtieth birthday}

\begin{abstract}
In this paper we develop the theory of equisingular 
deformations of plane curve singularities in arbitrary characteristic. 
We study equisingular deformations of the parametrization and 
of the equation and show that the base space of
its semiuniveral deformation is smooth in both cases. 
Our approach through deformations of the parametrization is elementary and 
we show that equisingular deformations of the parametrization
form a linear subfunctor of all deformations of the parametrization. 
This gives additional information, even in characteristic zero, the case 
which was treated by J. Wahl. The methods and proofs extend easily
to good characteristic, that is, when the characteristic does
not divide the multiplicity of any branch of the singularity.

In bad characteristic, however, new phenomena occur and we are naturally
led to consider weakly trivial respectively weakly equisingular deformations, 
that is, those which become trivial respectively equisingular
after a finite and dominant base change. The semiuniversal base space for
weakly equisingular deformations is, in general, not smooth but becomes 
smooth after a finite and purely inseparable base extension. 
For the proof of this fact we introduce some
constructions which may have further applications in
the theory of singularities in positive characteristic.
\end{abstract}

\thispagestyle{empty}
\vskip0.5cm

\tableofcontents

\thispagestyle{empty}
\vfill

\hfill \today

\newpage

\section*{Introduction}

We develop in this paper the theory of equisingular deformations of plane
algebroid curve singularities over an algebraically closed field $K$ of
arbitrary characteristic. If the curve singularity is given by the equation
$f=0$ where $f\in P=K[[x,y]]$ is a reduced formal power series we study deformations
of the local analytic ring $R=P/\langle f\rangle$, as well as deformations of the map $\varphi:P\to \overline{R}$
where $\overline{R}$ is the integral closure of $R$ in its total ring of
fractions. The first are called {\em deformations of the equation} and the
latter {\em deformations of the parametrization}. Since $P$
(resp. $\overline{R}$) are regular local (resp. semilocal) rings, deformations of the parametrization are very
simple objects and the semiuniversal object of the corresponding functor
$\uDef_{\Rbar\leftarrow P}$ of isomorphism classes of deformations of
\mbox{$\varphi:P\to\overline{P}$} can be
explicitely described in terms of a $K$--basis of its tangent space
$T^1_{\Rbar\leftarrow P}$. The same holds for deformations with sections $\uDefsec_{\Rbar\leftarrow P}$.

Equisingular deformations of the
parametrization $\varphi$ (with section) are defined by requiring equimultiplicity (along some
sections) for each infinitely near point of $P$ on $R$, in a compatible
manner. Of course, we have to consider only the finite set of essential
infinitely near points occurring in an embedded good resolution of $P/\langle
f\rangle$. We prove that the functor $\uDefes_{\Rbar\leftarrow P}$  of
{\em equisingular deformations of the parametrization} is a {\em linear} subfunctor of
$\uDefsec_{\Rbar\leftarrow P}$ and, therefore, has also an explicit description in
terms of a $K$--basis of its tangent space $\Ties_{\Rbar\leftarrow P}$
(Theorem \ref{theo3.8}). In particular, the base space of the semiuniversal
deformation of the parametrization is smooth. Furthermore, the
linearity allows an easy proof of the openness of versality for equisingular
deformations of the parametrization (Corollary \ref{coro3.13}).

The relation between deformations of the parametrization and deformations of the equation
is based on the fact that the deformation of $R$ can be  uniquely recovered
from  the deformation of  $\varphi: P\to \Rbar$. That is, the deformation functor
$\uDef_{\Rbar\leftarrow P}$ is natural isomorphic to the functor
$\uDef_{\Rbar\leftarrow R}$ of {\em deformations of the normalization}, that
is, of the normalization map $R\to \Rbar$. In the same way we get an
isomorphism $\uDefsec_{\Rbar\leftarrow P}\cong \uDefsec_{\Rbar\leftarrow R}$
for the corresponding deformations with section (Proposition \ref{prop:1.3})
and below we do not distinguish between these two functors. By forgetting $\Rbar$ we have a natural transformation
$\uDefsec_{\Rbar\leftarrow P}\to \uDefsec_R$ and we denote the image of
\mbox{$\uDefessec_{R\leftarrow P}$} in \mbox{$\uDefsec_R$}  by $\uDefessec_R$. We show that
equisingular sections of deformations of $R$ are unique (Proposition
\ref{prop:sigma unique}) and, hence, by forgetting the section, $\uDefessec_R$ is isomorphic to 
$\uDefes_R$. The latter is the functor of (isomorphism classes of) {\em equisingular
deformations of the equation} (or {\em of $R$}), which is a subfunctor of $\uDef_R$, the
(usual) deformations of $R$. The transformation $\uDefes_{\Rbar\leftarrow P}\to
\uDefes_R$ from equisingular deformations of the parametrization to
equisinglar deformations of the equation is, in general, not an isomorphism. However, we show that it is
smooth (Theorem \ref{theo:def of eqn}). This implies the first main result, that the base space of the
semiuniversal equisingular deformation of $R$ is smooth of dimension equal to
the vector space dimension of its tangent space $\Ties_R$, in any characteristic.

We have defined equisingular deformations of $R$ as those which lift to
deformations of the parametrization $P\to \Rbar$ such that this lifting is
equisingular along some sections. While the equisingular sections of
deformations of $R$ are unique, the lifting to equisingular sections of
deformations of $\Rbar$ are in general not unique in positive
characteristic. Indeed, the behaviour of equisingular deformations of the
equation depend, in contrast to equisingular deformations of the
parametrization, very much on the characteristic $p$ of the field $K$. We say
that the {\em characteristic is good} (with respect to $R$) if $p=0$ or if
$p>0$ and $p$ does not divide the multiplicity of any branch of $R$. We prove
that, if $p$ is good, then $\uDefes_{\Rbar\leftarrow R}\cong \uDefes_R$, hence
the lifting of equisingular sections of deformations of $R$ to those of
deformations of $\Rbar$ is unique up to isomorphism (Theorem \ref{theo:def of
  eqn}). Moreover, in this case we can show that the base space of the semiuniversal
equisingular deformation of $R$ can be represented by a (smooth) {\em algebraic} subscheme of the
(algebraic) base space of the semiuniversal deformation of $R$.

The theory of equisingular deformations of plane curve singularities in
characteristic zero has been initiated by J. Wahl in his thesis
(cf. \cite{Wa}). Wahl's approach is different from ours as he considers only
deformations of the equation and defines equisingularity by requiring
equimultiplicity of the equation of the reduced total transform along sections
through all essential infinitely near points of $P$ on $R$. Although
equimultiplicity of the parametrization is usually stronger than
equimultiplicity of the equation, one can prove that Wahl's functor
$\overline{ES}$ and our functor $\uDefes_R$ are isomorphic
  (cf. \cite{GrLoShu}). Thus, we get, in characteristic zero, a new proof of
  Wahl's result that the equisingularity stratum (which coincides then with
  the $\mu$--constant stratum, where $\mu$ is the Milnor number) in the base space of the semiuniversal
  deformation of $R$ is smooth. As mentioned above, the same result holds if the
  characteristic is good.

Our approach through deformations of the
parametrization appears to be quite simple and provides, even in characteristic
zero, additional information. This can be seen clearly in section
\ref{sec:exact seq} of this paper where we relate several infinitesimal deformations by
means of exact sequences which allow to compute not only $\Ties_R$
effectively but also gives, on the tangent level, a geometric interpretation
of the related deformation functors.

The construction of $\uDefes_R$ as subfunctor of $\uDef_R$ is so explicit that
it leads to an algorithm for the computation of a semiuniversal equisingular
deformation in good characteristic (that is, of the $\mu$--constant stratum in
characteristic $0$). This has been implemented in the computer
algebra system {\sc Singular} (cf. \cite{CaGrLo} for a description of the
algorithm).     

In bad characteristic, however, new phenomena occur. There are deformations
which are not equisingular but become equisingular after some finite (and
dominant) base change. We call such deformations {\em weakly equisingular} and
show that the functor of weakly equisingular deformations of $R$ has a
semiuniversal object. Its base space is, in general not smooth but it becomes
smooth after a finite and purely inseparable base extension. The proof of this fact is
rather involved and occupies sections \ref{sec:equipolygonal} and
\ref{sec:8}. We prove this by constructing a {\em weak equisingularity stratum} in
the base space of any deformation of $R$ (with section) which has a certain
universal property (Theorem \ref{theo:6.3} and Theorem \ref{thm:6.3}).

To prove the existence and properties of the weak equisingularity stratum  we give explicit conditions defining a
subscheme in the base space of the given deformation such that the restriction
of the family to this subscheme can be simultaneously blown up and satisfies
additional conditions preserved under further blowing ups. All conditions
together define then the weak equisingularity stratum.

We like to stress that, keeping
the multiplicity constant along a section in each blowing up is equivalent to
keeping the Newton diagramme (with respect to generic adapted coordinates)
constant. Moreover, we have to consider an {\em adapted Jacobian ideal} taking care of the fact that
leading terms of $f$ may vanish after differentiation of $f$. This leads to
deformations which we call {\em equipolygonal deformations} and which we study  in detail. 

If we start with a versal deformation of $R$ with smooth base space then the
defining conditions for the weak equisingularity stratum become smooth after a
purely inseparable base change and we construct the weak equisingularity
stratum together with its smooth covering space at the same time. This
construction is functorial and, which is a key point, versality for
equipolygonal deformations is preserved under blowing ups.

Although the weak equisingularity stratum is not smooth, it is reduced and irreducible,
becomes smooth after a purely inseperable base extension and, has good
functorial properties. In good characteristic it is even smooth and coincides
with the (strong) equisingularity stratum considered above and, therefore, weak
and strong equisingular deformations are the same in that case. 

In bad
characteristic, however, a largest ''strong equisingularity stratum'' does in
general not exist. Indeed, we show that inside the base space of the
semiuniversal deformation of $R$ there may be (infinitely many) different
smooth subschemes, having all the same tangent space, such the restriction of
the given semiuniversal deformation to them is equisingular (even semiuniversal
equisingular). Each of these smooth substrata may be considered as a strong
equisingularity stratum, but the restriction of the semiuniversal deformation to the union of two such strata is not
(strongly) equisingular. Moreover, the Zariski closure of all these strongly
equisingular strata is the weak equisingularity stratum which is then of bigger
dimension. 

In the last section we study the geometry of the different equisingular
strata. We identify, inside the smooth covering space of the weak equisingularity
stratum, an intrinsically defined subspace $\Tisec_{\Rbar/R}$, being the
tangent space to deformations of the normalization $R\to\Rbar$ which leave $R$
fixed. This space can explicitely be computed, it is zero in good
characteristic, and in bad characteristic its vanishing gives a
sufficient and necessary condition that a largest (strong) equisingularity
stratum exists, and then automatically coincides with the weak equisingularity stratum.
\bigskip

{\bf Acknowledgements:} 
Our collaboration on this paper was supported
by the universities of Kaiserslautern and Valladolid and by a Research in
Paris stay at the Mathematisches Forschungsinstitut Oberwolfach. We like to
thank these institutions for their hospitality and support.

\bigskip
%\hrulefill Baustelle \hrulefill

% Assume we working over a field $K$ of {\em good characteristic},
%  that is, the characteristic of $K$ is either $0$ or it does not divide any of
%  the multiplicities of the branches of $R$. In this case, we have
%  \mbox{$\uDefes_{\Rbar\leftarrow P}\cong \uDefes_R$}, and the concepts of
%  strong and weak equisingular deformations (as introduced in
%  Section \ref{sec:6}) coincide, and we just speak of 
%  {\em equisingular deformations of the equation\/} (see Section \ref{}). 
% fix: section referenz - beweis dass strong=weak

%\bigskip

\newpage

\section{Deformations of the Parametrization and of the Normalization}%%%
\label{sec:DeformParametr}

%\noindent
In this section we fix the notations and state some basic facts about
deformations of the normalization and deformations of the parametrization for a
reduced plane curve singularity. We shall put special emphasis on deformations
with section.

$K$ denotes an algebraically closed field of characteristic \mbox{$p \ge 0$}.
If $A$ is a Noetherian complete local $K$-algebra with maximal ideal $\fm_A$,
we always assume that \mbox{$A/\fm_A = K$}. The category of these algebras is
denoted by $\sA_K$. $\Keps$ denotes the two-dimensional $K$-algebra with
\mbox{$\varepsilon^2 = 0$}. 

We consider reduced algebroid plane curve singularities $C$ over $K$, defined
by a formal power series \mbox{$f \in K[[x,y]]$}. Usually, we work with the 
complete local ring of \mbox{$C=\Spec(R)$},
$$
R = P/\langle f\rangle\,,\quad P = K[[x,y]]\,.
$$
If \mbox{$f=f_1\cdot\ldots\cdot f_r$} is an irreducible factorization of $f$
in $P$,
the rings 
$$
R_i = P/\langle f_i\rangle\,,\quad i=1,\dots,r\,,
$$
are the complete local rings of the branches of $C$. The lowest degree
\mbox{$\ord_{x,y} (f)$} of a monomial appearing in the power series development of
\mbox{$f\neq 0$} is called the {\em multiplicity} of $f$ and denoted by $\mt
(f)$; we set $\mt(0)=\infty$. Of course, $\mt (f)=\mt(f_1)+\cdots+
\mt(f_r)$. We say that the {\em characteristic of $K$ is good (for $R$)} if is
does not divide $\mt(f_i), \text{ for all } i=1, \ldots, r$. 

The normalization $\Rbar$ of $R$ is the integral closure of $R$ in its total
ring of fractions $\Quot(R)$. $\Rbar$ is the direct sum of the normalizations
$\Rbar_i$ of $R_i$, \mbox{$i=1,\dots,r$}, hence a semilocal ring. Each
$\Rbar_i$ is a discrete valuation ring, and we can choose uniformizing
parameters $t_i$ such that \mbox{$\Rbar_i\cong K [[t_i]]$}. After fixing the
parameters $t_i$, we identify $\Rbar_i$ with \mbox{$K [[t_i]]$} and 
get 
$$
\Rbar = \bigoplus^r_{i=1} \Rbar_i =  \bigoplus^r_{i=1} K [[t_i]]\,.
$$
The {\em normalization map\/} \mbox{$\nu: R\to \Rbar$} (induced by the
inclusion \mbox{$R\hookrightarrow \Quot(R)$}) is then given by the (primitive)
{\em parametrization of $R$\/} (or {\em of $C$\/}),
$$
\varphi=(\varphi_1,\dots,\varphi_r):P\lra \Rbar = \bigoplus^r_{i=1} K
[[t_i]]\,,$$ 
where \mbox{$\varphi_i(x) =x_i(t_i)$},  \mbox{$\varphi_i(y) =y_i(t_i)\in
  K[[t_i]]$}, \mbox{$i=1,\dots,r$}. Since \mbox{$\langle
  f\rangle=\Ker(\varphi)$}, $R$ may be recovered from $\varphi$. We call
\[
\ord\varphi_i:=\min\{\ord_{t_i}x_i, \ord_{t_i}y_i\}
\]
the {\em multiplicity} (or {\em order}) of $\varphi_i$ and the $r$--tupel
$\ord(\varphi)=(\ord \varphi_1, \ldots, \ord \varphi_r)$ the {\em
  multiplicity} (or {\em order) of the parametrization} $\varphi$. Note that
$\ord \varphi_i$ is the maximal integer $m_i$
s.t. $\varphi_i(\fm_P)\subset\langle t_i\rangle^{m_i}$. Moreover, we have
(cf. \cite{C1}) 
\[
mt(f)=\ord \varphi_1+\ldots+\ord\varphi_r.
\]

\begin{definition}\label{def:1.1}
  A {\em deformation with sections of the parametrization\/} of $R$ over
  \mbox{$A\in \sA_K$} is a commutative diagram with Cartesian squares  
$$
\UseComputerModernTips
\xymatrix@C=9pt@R=2pt@M=6pt{
{\Rbar}\ar@{}[ddrr]|{\Box} && \ar@{->>}[ll]
{\Rbar}_A\ar@/^3pc/[dddd]^{\osigma=\{\osigma_i\mid i=1,\dots,r
  \}}\\   
\\
P \ar[uu]^{\varphi}\ar@{}[ddrr]|{\Box} && \ar@{->>}[ll] P_A
\ar[uu]_{\varphi_A} \ar@/^/[dd]^{\sigma} \\ 
 \\
 K \ar@{^{(}->}[uu] && \ar@{->>}[ll] A
\ar[uu]
}
$$
with \mbox{$\Rbar_A=\bigoplus_{i=1}^r \Rbar_{A,i}$}, where $\Rbar_{A,i}$,
\mbox{$i=1,\dots,r$}, and $P_A$ are Noetherian complete local $K$-algebras
which are flat over $A$. $\sigma$ is a section of \mbox{$A\to P_A$}, and
$\osigma_i$ is a section of \mbox{$A\to \Rbar_{A,i}$}, \mbox{$i=1,\dots,r$}. We
denote such a deformation by \mbox{$\xi=(\varphi_A,\osigma,\sigma)$}.  

A morphism from \mbox{$\xi$} to another deformation
\mbox{$(P_B\!\xrightarrow{\varphi_B}\! {\Rbar}_B,\osigma_B,\sigma_B)$}
over \mbox{$B\in \sA_K$} is then given by morphisms of local $K$-algebras
\mbox{$A\to B$}, \mbox{$P_A\to P_B$} 
and \mbox{$\Rbar_{A,i}\to \Rbar_{B,i}$} such that the resulting diagram
commutes. The category of such deformations is denoted by
$\Defsec_{\Rbar\leftarrow P}$. If we consider only deformations over a fixed
base $A$, we obtain the (non-full) subcategory \mbox{$\Defsec_{\Rbar\leftarrow
    P}(A)$} with morphisms being the identity on $A$.  \mbox{$\Defsec_{\Rbar\leftarrow
    P}$} is a {\em fibred gruppoid\/} over $\sA_K$, in particular, each
morphism in  \mbox{$\Defsec_{\Rbar\leftarrow
    P}(A)$} is an isomorphism.

Giving \mbox{$\xi=(\varphi_A, \osigma, \sigma)$} and a morphism \mbox{$\psi:A\to B$} in $\sA_K$, the 
{\em induced deformation}\!
  \mbox{$\psi\xi= (\psi\varphi_A, \psi\bar{\sigma}, \psi\sigma)$} is an object in \mbox{$\Defsec_{\Rbar\leftarrow P}(B)$}, 
defined by $B\to P_A\widehat{\otimes}_A
  B\to {\Rbar}_A\widehat{\otimes}_A B$, 
%$$
%\UseComputerModernTips
%\xymatrix@C=9pt@R=2pt@M=6pt{
%{\Rbar}\ar@{}[ddrr]|{\Box} && \ar@{->>}[ll]
%{\Rbar}_A\widehat{\otimes}_A B
%\ar@/^3pc/@<1ex>[dddd]^{\psi\osigma=\{\psi\osigma_i\mid i=1,\dots,r 
%  \}}\\   
%\\
%P \ar[uu]^{\varphi}\ar@{}[ddrr]|{\Box} && \ar@{->>}[ll]
%P_A\widehat{\otimes}_A B 
%\ar[uu]_{\psi\varphi_A} \ar@/^/[dd]^{\psi\sigma} \\ 
% \\
% K \ar@{^{(}->}[uu] && \ar@{->>}[ll] B
%\ar[uu]^{\psi\iota}
%}
%$$
with morphisms \mbox{$b\mapsto 1\widehat{\otimes}\!\; b$},
\mbox{$\psi\varphi_A=\varphi_A\widehat{\otimes}\id_B$},
\mbox{$\psi\sigma:h\widehat{\otimes}\!\; b\mapsto \psi(\sigma(h))\cdot b$}, 
\mbox{$\psi\osigma_i:r\widehat{\otimes}\!\; b\mapsto \psi(\osigma_i(r))\cdot
  b$}. Here, $\widehat{\otimes}$ denotes the complete tensor product.

The set of isomorphism classes of objects in $\Defsec_{\Rbar\leftarrow P}(A)$
is denoted by $\uDefsec_{\Rbar\leftarrow P}(A)$, and 
\mbox{$  \uDefsec_{\Rbar\leftarrow P} : \sA_K \to \text{(Sets)} $}
denotes the corresponding {\em deformation functor\/} (which always
refers to isomorphism classes). Moreover, we denote by
\mbox{$\Tisec_{\Rbar\leftarrow P}:=\uDefsec_{\Rbar\leftarrow P}(\Keps)$} the
tangent space to this functor.
\end{definition}

\begin{remark}\label{rmk:remark1.2}
  Since $P$ and the $\Rbar_i$ are regular local rings, any deformation of $P$
  and of $\Rbar$ is trivial. That is, there are isomorphisms 
  \mbox{$P_A\cong A[[x,y]]$}  and \mbox{$\Rbar_A\cong \bigoplus_{i=1}^r
    A[[t_i]]$} over $A$, mapping 
  the sections $\sigma$ and $\osigma_i$ to the trivial sections. Hence, any
  object in  $\Defsec_{\Rbar\leftarrow P}(A)$ is isomorphic to a diagram of the
  form 
$$
\UseComputerModernTips
\xymatrix@C=9pt@R=10pt@M=6pt{
\bigoplus\limits_{i=1}^r K[[t_i]]\ar@{}[drr]|{\Box} && \ar@{->>}[ll]
\bigoplus\limits_{i=1}^r
A[[t_i]]\ar@/^3pc/[dd]^{\osigma=\{\osigma_i\,\mid i=1,\dots,r 
  \}}\\
K[[x,y]] \ar[u]^-{\varphi}\ar@{}[drr]|{\Box} && \ar@{->>}[ll] A[[x,y]]
\ar[u]_-{\varphi_A} \ar@/^/[d]^{\sigma} \\ 
 K \ar@{^{(}->}[u] && \ar@{->>}[ll] A
\ar[u]
}
$$
where $\varphi_A$ is the identity on $A$ and $\sigma$, $\osigma_i$ are the
trivial sections (that is, the canonical epimorphisms mod $x,y$, respectively
mod $t_i$). Hence, 
$\varphi_A$ is given by \mbox{$\varphi_A=(\varphi_{A,1},\dots,\varphi_{A,r})$},
where $\varphi_{A,i}$ is determined by 
$$ \varphi_{A,i}(x)=X_i(t_i)\,,\quad \varphi_{A,i}(y)=Y_i(t_i)\in
t_iA[[t_i]]\,,$$ 
\mbox{$i=1,\dots,r$}, such that \mbox{$X_i(t_i)\equiv x_i(t_i)$},
\mbox{$Y_i(t_i)\equiv y_i(t_i)$} mod $\fm_A$. 
\end{remark}

%\noindent
Similarly to Definition \ref{def:1.1}, by replacing $P$ by $R$ and $\varphi$
  by $\nu$ resp. $P_A$ by $R_A$ and $\varphi_A$ by $\nu_A$, we can define {\em deformations with
  section of the normalization} \mbox{$R\to \Rbar$}, and obtain
  the category \mbox{$\Defsec_{\Rbar\leftarrow R}$} resp. the deformation  functor $\uDefsec_{\Rbar\leftarrow R}$. Indeed, we are mainly interested
  in the latter functor, which a priori is more complicated than
  $\uDefsec_{\Rbar\leftarrow   P}$, since $R$ is not regular in
contrast to $P$. The following proposition shows that both functors are isomorphic. To prove this, we have to consider deformations with section of the
sequence of morphisms \mbox{$P\to R\to\Rbar$}, whose definition is analogous
to Definition \ref{def:1.1}. Such deformations may be called {\em deformations
  of the normalization with embedding}. The corresponding deformation
category, respectively the corresponding deformation functor of isomorphism classes, 
is denoted by  $\Defsec_{\Rbar\leftarrow R\leftarrow P}$, respectively by
$\uDefsec_{\Rbar\leftarrow R\leftarrow P}$. 

\begin{proposition}\label{prop:1.3}
The forgetful functor from $\Defsec_{\Rbar\leftarrow R\leftarrow P}(A)$ to
$\Defsec_{\Rbar\leftarrow P}(A)$, respectively to $\Defsec_{\Rbar\leftarrow
  R}(A)$, is an isomorphism, respectively smooth. Both induce an
isomorphism between the corresponding 
deformation functors. In particular, 
$$\uDefsec_{\Rbar\leftarrow P} \cong \uDefsec_{\Rbar\leftarrow R}\,.
$$  
Moreover, if \mbox{$(P_A\!\xrightarrow{\varphi_A}\!
  \Rbar_A,\,\osigma,\,\sigma)$} is an object of 
$\Defsec_{\Rbar\leftarrow P}(A)$, then $\ker(\varphi_A)$ is a principal ideal,
and the lifting of \mbox{$(\varphi_A,\osigma,\sigma)$} to an object of 
$\Defsec_{\Rbar\leftarrow R\leftarrow P}(A)$ is obtained by setting \mbox{$R_A
  =P_A/\ker(\varphi_A)$}.
\end{proposition}

\begin{proof}
  Let \mbox{$(P_A\!\xrightarrow{\varphi_A}\!
  \Rbar_A,\,\osigma,\,\sigma)$} be a deformation with sections of the
parametrization \mbox{$P \to {\Rbar}$} over \mbox{$A\in \sA_K$}. Since 
$\varphi_A$ is quasifinite, ${\Rbar}_A$ is a finite $P_A$-module, and we have a
minimal free resolution 
\begin{equation}
  \label{eq:eq1}
 0 \longleftarrow {\Rbar}_A \longleftarrow F_0 \stackrel{M}{\longleftarrow}
F_1 \longleftarrow F_2 \longleftarrow \ldots
\end{equation}
of ${\Rbar}_A$ as a $P_A$-module.  Since $P_A$ and $\Rbar_A$ are $A$-flat,
the exactness of the sequence \eqref{eq:eq1} is preserved when tensorizing with
$\otimes_A K$, obtaining in this way a minimal free resolution of ${\Rbar}$ as
$P$-module with presentation matrix \mbox{$M_0=M\otimes_A K$}. Since $P$
is regular of dimension two, and since $\Rbar$ has depth 
one, the Auslander-Buchsbaum formula gives that the minimal resolution of
${\Rbar}$ has length one. Thus, $M_0$ is injective. By the local criterion of
flatness (cf. [Ma, Theorem 22.5]), $M$ is injective too, and the free $P_A$-modules $F_0$ and $F_1$ have
the same rank. 

The ideal \mbox{$\langle \det(M)\rangle\subset P_A$} is independent of the
chosen resolution, and we set \mbox{$R_A:=P_A/\langle \det(M)\rangle$},
which is flat over $A$. Note
that the ideals $\langle f\rangle$ and $\langle \det(M_0)\rangle$ of $P$
have the same support and coincide in the generic points where
\mbox{$R=P/\langle 
f\rangle$} is regular. It follows that the two principal ideals \mbox{$\langle
f\rangle$} and $\langle
\det(M_0)\rangle$ of $P$ coincide.   

Since $\det(M)$ annihilates ${\Rbar}_A$ by Cramer's rule, and since the kernel of $\varphi_A$
is equal to the annihilator of ${\Rbar}_A$ as $P_A$-module, the
canonical projection \mbox{$P_A\twoheadrightarrow R'_A:=P_A/\ker
  \varphi_A$} induces a surjection \mbox{$R_A \twoheadrightarrow R'_A$}. The
kernel of this surjection is supported by the singular locus of the fibres and is zero after tensorizing with 
$\otimes_A K$.  Thus, by 
Nakayama's lemma, \mbox{$R_A = R'_A$}. It follows that $\varphi_A$ factors as
\mbox{$P_A\twoheadrightarrow 
  R_A\hookrightarrow \Rbar_A$}, defining in this way an object of
$\Defsec_{\Rbar\leftarrow R\leftarrow P}(A)$.
Moreover, if \mbox{$P_A\twoheadrightarrow 
  R_A''\hookrightarrow \Rbar_A$} is any lifting of \mbox{$P_A\to \Rbar_A$} to
an object of $\Defsec_{\Rbar\leftarrow R\leftarrow P}(A)$,
then \mbox{$P_A\twoheadrightarrow 
  R_A''$} is surjective and \mbox{$R_A''\hookrightarrow \Rbar_A$} is injective
(by Nakayama's lemma). Thus, as before, \mbox{$R_A''=R_A$}. As morphisms in 
$\Defsec_{\Rbar\leftarrow P}(A)$  can be uniquely lifted, too, this shows that 
the forgetful functor induces an isomorphism of categories 
\mbox{$\Defsec_{\Rbar\leftarrow
    R\leftarrow P}(A)\cong  
\Defsec_{\Rbar\leftarrow P}(A)$}. 

To get the statements for deformations (with section) of
the normalization, note that $P$ is a regular local ring. Hence, any
deformation (with section) of $R$ may be lifted to an ``embedded'' deformation,
that is, a deformation of \mbox{$P\to R$}. Thus, the forgetful
functor induces a surjection
\mbox{$\Defsec_{\Rbar\leftarrow R\leftarrow 
    P}(A)\twoheadrightarrow \Defsec_{\Rbar\leftarrow R}(A)$}. 
The fibre is a principal homogeneous space under isomorphisms of $P$ fixing
$R$, showing smoothness. Moreover, if two
deformations $R_A$ and $R'_A$ of $R$ over $A$ are isomorphic, the isomorphism
\mbox{$R_A\cong R'_A$} may be lifted to an isomorphism \mbox{$P_A\cong P'_A$}
(since \mbox{$P_A\cong A[[x,y]]$}). 
\end{proof}

\begin{remark}
  If we omit the sections in Definition \ref{def:1.1} and in the subsequent
  discussion, we get analogous results for deformations without sections. The
corresponding categories, respectively deformation functors are denoted by
$\Def_{\Rbar\leftarrow P}$, respectively by $\uDef_{\Rbar\leftarrow P}$,
etc.\,. Indeed, as the proof of Proposition \ref{prop:1.3} shows, the
sections do not affect the arguments at all. Hence, Proposition \ref{prop:1.3}
remains true with the superscript '$\sec$' being omitted.
\end{remark}

\section{Equisingular Deformations of the Parametrization}%%% section 2
\label{sec:section 2}

%\noindent
In this section, we define equisingular deformations of the
parametrization with section, and we discuss the uniqueness of the
sections.

In order to define equisingular deformations, we
consider the notion of infinitely near points.

Consider the natural diagram of graded $K$-algebras given by 
$$
\xymatrix@C=9pt@R=12pt
{
\gr_{\fm_i} (R_i) && \gr_{\fm}(R)\ar[ll]  && \gr_{\fm_P}(P)\ar[ll]\\
\Gr_{\fm_i} (R_i)\ar@{->>}[u] && \Gr_{\fm}(R)\ar[ll]\ar@{->>}[u] &&
\Gr_{\fm_P}(P)\ar[ll]\ar@{->>}[u]\\ 
}
$$
where $\fm,\fm_i,\fm_P$ are the respective maximal ideals of $R,R_i,P$ and
where, for an ideal \mbox{$I\subset S$}, \mbox{$\gr_I(S)=
  \bigoplus_{k=0}^\infty I^k/I^{k+1}$},  \mbox{$\Gr_I(S)=
  \bigoplus_{k=0}^\infty I^k$}. Applying $\Proj$ to the above diagram, we get
the blow up schemes \mbox{$\Bl_{\fm_i}(R_i)$}, \mbox{$\Bl_\fm(R)$}, 
\mbox{$\Bl_{\fm_P}(P)$} in the lower row, and the corresponding exceptional
divisors $E_{\fm_i}$,  $E_{\fm}$, and  $E_{P}$ in the upper row. $\Proj$
applied to the vertical maps in the diagram gives rise to natural embeddings of these
objects.  

\begin{definition}\label{def:2.1}
An {\em infinitely near point\/} $P'$ {\em in the first infinitesimal
  neighbourhood\/} of $P$ is the completion of the local ring of a closed point
$O$ on the exceptional divisor $E_P$ in \mbox{$\Bl_{\fm_P}(P)$}. We always use
the same notation $O$ for the local ring and for the point in
\mbox{$\Bl_{\fm_P}(P)$}. 
\end{definition}

%\noindent
Since exceptional divisors are projectivizations of tangent cones, $E_P$ is a
projective line, and the image of $E_{\fm_i}$ in $E_P$ is one point $O_i$,
counted $m_i$ times, where $m_i$ is the multiplicity of the branch $R_i$. Among
the infinitely near points in the first neighbourhood, those of type
\mbox{$P'=\widehat{O}_i$} for some \mbox{$1\leq i\leq r$} are called {\em
  infinitely near points (of $P$) on $R$}. For such a $P'$, we set
$$ \Lambda_{P'}:=\left\{ i\in \{1,\dots,r\} \:\left|\: P'=\widehat{O}_i
  \right.\right\}\,. $$
For each \mbox{$i\in \Lambda_{P'}$}, we also say that the branch {\em $R_i$
    passes through $P'$}.

Note that we refer to $P$ itself as an infinitely near point of $P$ (in the
$0$-th infinitesimal neighbourhood) on $R$.

\begin{remark}\label{rmk:2.4}
In analytical terms, we have \mbox{$P=K[[x,y]],\ \gr_{\fm_P} (P)=K[x,y]$},
  \mbox{$\gr_{\fm} (R)=K[x,y]/\langle J_mf\rangle$}, and \mbox{$\gr_{\fm_i}
  (R_i)=K[x,y]/\langle J_{m_i}f_{i}\rangle$},  
\mbox{$i=1,\dots,r$}, where $J_{m_i}f_{i}$ denotes the sum of terms of
smallest degree \mbox{$m_i=\mt(f_i)$} in the power series expansion of $f_i$, and
\mbox{$J_mf= \prod_{i=1}^r J_{m_i}f_{i}$}. 

Note that each $J_{m_i}f_{i}$ is the 
$m_i$-th power of a non--zero linear form \mbox{$\alpha_iy-\beta_ix$},
\mbox{$\alpha_i,\beta_i\in K$}. We also say that it {\em corresponds to the
  tangent direction \mbox{$(\alpha_i\!\!\::\!\!\:\beta_i)\in \P^1_K$}}. If the
infinitely near point \mbox{$O\in E_P$} corresponds to the tangent direction
\mbox{$(1:\beta)\in \P^1_K$}, then
\mbox{$O=P\left[\tfrac{y}{x}\right]_{\langle 
    x,y'\rangle}$}, \mbox{$y'=\tfrac{y}{x}-\beta$}, 
and \mbox{$P'=\widehat{O}= K[[x,y']]$}.  

Let \mbox{$\pi':P\to P'$} be the
blow--up map \mbox{$x\mapsto x,\ y\mapsto x(y'+\beta)$}. For $g\in P$ (of
multiplicity $m$) such that \mbox{$J_mg=c(y-\beta x)^m$} for some \mbox{$c\in
  K^\ast$} we set 
\[
g'=x^{-m} \pi' (g)\ ,\ \widehat{g}=xg'
\]
and call \mbox{$R'=P'/\langle g'\rangle$ } the {\em strict transform}
resp. \mbox{$\widehat{R}=P'/\langle \widehat{g}\rangle$} the
{\em reduced total transform} of \mbox{$R=P/\langle g\rangle$}. Moreover,
\mbox{$P'/\langle \pi' (g)\rangle$} resp. \mbox{$P'/\langle x\rangle$} are
the {\em total transform} of $R$ resp. the {\em exceptional curve} of the blow--up. In particular, if \mbox{$O=O_i$} (that is,
\mbox{$\alpha_i\neq 0$} and \mbox{$\beta=\frac{\beta_i}{\alpha_i}$}), then \mbox{$R'_i=P'\left/\left\langle f'_i\right\rangle \right.$}, respectively
\mbox{$R'=P'\big/\bigl\langle \prod_{i\in
      \Lambda_{P'}} f'_i \bigr\rangle$},
is the strict transform of the branch $R_i$, respectively of
$R$, at $P'$. 

For each infinitely near point $P'$ in the first infinitesimal neighbourhood
on $R$, and for each \mbox{$i\in \Lambda_{P'}$}, the normalization of the
strict transform $R_i'$ is $\Rbar_i$, and the normalization of $R'$ is 
\mbox{$\Rbar\!\,':=\bigoplus_{i\in \Lambda_{P'}}\Rbar_i$}. Further, notice that
the parametrization \mbox{$P'\to\Rbar\!\,'$} of $R'$ is given by
\mbox{$x_i(t_i), y'_i(t_i)$}, \mbox{$i\in \Lambda_{P'}$}, where 
$$
y'_i(t_i):=\frac{y_i(t_i)}{x_i(t_i)}-\frac{\beta_i}{\alpha_i}\in 
K[[t_i]]\,.$$ 
Here again, we assume that \mbox{$\alpha_i\neq 0$} (which is no restriction, since
after a general linear change of the coordinates $x,y$ 
it holds for all \mbox{$i=1,\dots,r$}).
\end{remark}

We extend the above definitions to higher infinitesimal neighbourhoods by induction.

\begin{definition}\label{def:2.3}
  Let \mbox{$k\geq 2$}, and assume that the infinitely near points $P'$ on $R$ in the
  \mbox{$(k\!\!\:-\!\!\:1)$}-th neighbourhood of $P$ are
    defined. Assume also that for each of these points $P'$ a set
    \mbox{$\Lambda_{P'}$}, the strict transform $R'_i$ of each branch $R_i$, 
\mbox{$i\in \Lambda_{P'}$}, and the strict transform $R'$ of $R$ at $P$ are
defined, as well as the reduced total transform $\widehat{R'}_i$
  resp. $\widehat{R'}$ and the exceptional curve $E'$.  Then, we call each infinitely near point $P''$ on $R$ in the first
infinitesimal neighbourhood of such a point $P'$ an 
{\em infinitely near point on $R$ in the $k$-th neighbourhood of
$P$}. We introduce
$$ \Lambda_{P''}:=\left\{i\in \Lambda_{P'}\:\big|\: R'_i \text{ passes through
  }P''\right\} $$
and define the {\em strict transform} of $R_i$ (respectively $R$) at $P''$ to
be the strict transform of $R'_i$ (respectively $R'$) at $P'$. Moreover, the
  reduced total transform of $R_i$ (resp. $R$) at $P''$ is the reduced total
  transform of $\widehat{R'_i}$ (resp. $\widehat{R'}$) and the exceptional
  curve at $P''$ is the reduced total transform of $E'$.

Given infinitely near points $P',P''$ as above, we call $P''$ {\em
  consecutive to $P'$}. According to Remark \ref{rmk:2.4}, if
\mbox{$P'=K[[u,v]]$}, \mbox{$P''=K[[w,z]]$}, then up to interchanging $u$ and
$v$, we can assume that \mbox{$w=u$}, \mbox{$z=\frac{v}{u}-\beta$} for some
\mbox{$\beta \in K$}. The map \mbox{$P'\to P''$} is called a {\em formal
  blow-up\/} (of the maximal ideal $\fm_{P'}$ in $P'$), as it satisfies the
following two properties:
\begin{enumerate}
\item[(i)] \mbox{$\fm_{P'} P'' = \langle u,v\rangle \cdot P''$} is a
  principal ideal, and  
\item[(ii)] there is no proper subalgebra \mbox{$S\in \sA_K$} of $P''$ with
  \mbox{$\fm_{P'}S$} being a principal ideal. 
\end{enumerate}
We call a point $P'$ an {\em infinitely near point\/} of $P$ on $R$ if it is
an infinitely near point on $R$ in the $k$-th neighbourhood of $P$ for some
\mbox{$k\geq 0$}. The above consideration shows that infinitely near points of
$P$ on $R$ are related to $P$ by compositions of formal blow-ups. An
infinitely near point $P'$ of $R$ is called {\em free\/} (resp. {\em satellite\/}) if
exactly one (resp. two) components of the exceptional curve $E'$ pass through
$P'$. The point $P$ itself is considered as free. 

%We say that an infinitely near point $P'$ of $P$ on $R$ is {\em essential for
%  $R$} iff it is not consecutive to an infinitely near point where the strict
%transform of $R$ is smooth. The set of such points is denoted by $\Ess(R)$.
\end{definition}

%\noindent
%Although the set of infinitely near points $P'$ on $R$ is infinite, for most
%purposes, we may restrict ourselves to considering only finitely many
%infinitely near points. Indeed, those 
%$P''$ %which are in the $k$-th neighbourhood (\mbox{$k\geq 1$}) of a point $P'$
%for which the strict transform of $R$ satisfies \mbox{$R'=R'_i=\Rbar_i$} for
%some index $i$ (and intersects the exceptional
%divisor transversally) usually do not have to be considered. By the theorem of
%resolution of singularities (see, e.g.\ \cite{C1,Lipman,Zariski}), this
%behaviour of the strict transform may be achieved after finitely many formal
%blow-ups. In particular, the set $\Ess(R)$ is finite. 

We say that an infinitely near point $P'\neq P$ is {\em essential\/} for $R$ if the reduced total transform of $R$ at $P'$ is not a
  node (i.e. a normal crossing of two smooth branches). $P$ itself is essential if
  $R$ is not smooth. The set of essential points for $R$ is denoted by Ess$(R)$. The set \mbox{Ess$(R)$} will be considered for an embedded (good)
  resolution of $R$. By the theorem of resolution of singularities
  (cf. e.g. \cite{C1}, \cite{Lipman}, \cite{Zariski}) Ess$(R)$ is finite.

\begin{definition} 
We define the {\em multiplicity} (or {\em order\/})
  {\em of a deformation with sections of the parametrization}
  \mbox{$(\varphi_A,\osigma,\sigma)$}, to be the $r$-tuple
  $\bord(\varphi_A\osigma, \sigma):=\bm=(m_1,\dots,m_r)$ such that
  \mbox{$\varphi_{A,i}(I_{\sigma})\subset I_{\osigma_i}^{m_i}$} and $m_i$ is
  the maximal integer with this property. Here, 
  \mbox{$I_\sigma=\ker \sigma\subset P_A$} and  \mbox{$I_{\osigma_i}=\ker
    \osigma_i\subset \Rbar_{A,i}$} are the ideals of the sections. A
  deformation with sections \mbox{$(\varphi_A, \osigma, \sigma)$} of $\varphi$
  is called 
  {\em equimultiple\/} (or, an {\em em-deformation\/}) if 
\[
\bord (\varphi_A, \osigma, \sigma)=\bord(\varphi).
\]
 We introduce the category, resp.\ deformation
  functor, of em-deformations with sections of the parametrization,
  $\Defem_{\Rbar\leftarrow P}$, resp.\   $\uDefem_{\Rbar\leftarrow P}$.
\end{definition}

\begin{remark}
A deformation \mbox{$\varphi_A: A[[x,y]] \to \bigoplus_{i=1}^r
    A[[t_i]]$} of the parametrization with trivial sections as in Remark
  \ref{rmk:remark1.2}, given by power series \mbox{$X_i(t_i),Y_i(t_i)\in A[[t_i]]$}  
 is equimultiple iff, for each \mbox{$1\leq
  i\leq r$}, the minimum of \mbox{$\ord_{t_i} X_i$} and \mbox{$\ord_{t_i} Y_i$}
coincides with the minimum of the $t_i$-orders of the residues $x_i,y_i$ mod
$\fm_A$. If this minimum is attained by, say $x_i(t_i)$, this means that the
coefficient of the term of smallest $t_i$-degree in $X_i$ is a unit in $A$. If
the deformation with trivial sections defined by $\varphi_A$ is 
equimultiple, then each generator \mbox{$F\in A[[x,y]]$} of $\ker \varphi_A$
(which is a principal ideal due to  Proposition \ref{prop:1.3}) defines an
em-deformation of \mbox{$R=K[[x,y]]/\langle f\rangle$}, i.e.
    \mbox{$\ord_{x,y}(F)= \ord_{x,y}(f)$} (cf. \cite{GrLoShu}, Lemma 2.26 for a
    proof for $k=\C$ which can be modified to work for arbitrary $K$).

\quad Notice, however,
  that the converse is not true if $A$ is not
  reduced. 
  For instance, consider the irreducible plane
  curve singularity \mbox{$R=K[[x,y]]/\langle x^5\!+y^3\rangle$}. The
  deformation of the parametrization  (with trivial sections) over \mbox{$A =
    K[\varepsilon]$} given by \mbox{$X(t) = t^3  -3\varepsilon t$},
  \mbox{$Y(t)= t^5+ 5\varepsilon t^3$} is obviously not equimultiple as a deformation of the
  parametrization. But, the corresponding deformation of $R$, which is given by
  \mbox{$F=x^5\!+y^3$}, is trivial, hence equimultiple.
\end{remark}

\begin{definition}\label{def:2.6}
%\begin{enumerate}
%\item 
(1) An {\em equisingular deformation of the
    parametrization \mbox{$P \to {\Rbar}$}} (or {\em
    es-deformation of \mbox{$P \to {\Rbar}$}}) over $A$ is a deformation with
  sections $(\varphi_A,\osigma,\sigma)$ of the parametrization which is
    equimultiple and which satisfies:  

For each infinitely near point $P'$ on $R$ there exists a deformation
  \mbox{$(\varphi'_A,\osigma',\sigma')$} of the parametrization \mbox{$P' \to {\Rbar}\!\,'$} over $A$ such that the  following diagram is commutative with Cartesian squares
$$
\UseComputerModernTips
\xymatrix@C=3pt@R=5pt{
\ {\Rbar}\phantom{\big|} \ar@{->>}[dr] &&&&&&&
{\Rbar}_A \ar@{->>}[dl]\ar@{->>}[lllllll] \ar@/^4pc/[lddddddd]^-{\osigma}
  \\ 
&{\Rbar}\!\,'\ar@{}[ddrrrrr]|{\Box} &&&&& \ar@{->>}[lllll]
{\Rbar}\!\,'_A
\ar@/^3pc/[dddddd]^-{\osigma'}|(0.2)\hole
\\ 
\\
&P'\phantom{\big|}\!\! \ar[uu]_{\varphi'}\ar@{}[ddrrrrr]|{\Box} &&&&&
\ar@{->>}[lllll]
P'_A
\ar@/^2pc/[dddd]^(0.3){\sigma'} |(0.2)\hole
\ar[uu]^{\varphi'_A} \\
\\ 
&P\phantom{\big|}\!\!
\ar[uu]_{\pi'}\ar@/^1pc/[uuuuul]^{\varphi}\ar@{}[ddrrrrr]|{\Box} 
&&&&& 
\ar@{->>}[lllll] 
P_A \ar@/^/[dd]^{\sigma} 
\ar@/_1pc/[uuuuur]_(0.6){\varphi_A}
\ar[uu]^{\pi'_A} \\ 
\\
& K \phantom{\big|}\!\! \ar@{^{(}->}[uu] &&&&& \ar@{->>}[lllll] A
\ar[uu]
}
$$
and the following conditions hold:

\begin{enumerate}
\item[(i)] \mbox{$\osigma'_i=\osigma_i:\Rbar_{A,i}\to A$} for all \mbox{$i\in
  \Lambda_{P'}$}. If \mbox{$P'=P$} then \mbox{$(\varphi'_A, \osigma',
  \sigma')=(\varphi_A, \osigma, \sigma)$}. 
\item[(ii)] \mbox{$(P'_A\xrightarrow{\varphi'_A}\Rbar'_A, \osigma', \sigma')$}
  is equimultiple, i.e. an object of \mbox{$\Defem_{\Rbar'\leftarrow P'}(A)$}.
\item[(iii)] The system of such diagrams is {\em compatible}: that is, if
  $P''$ on $R$
  is in some infinitesimal neighbourhood of $P'$ then there exists a
  morphism 
  \mbox{$P'_A \to P''_A$} such that the obvious diagram commutes.
\item[(iv)] If $P''$ is consecutive to $P'$, then \mbox{$P'_A \to P_A''$} is
  a {\em formal blow-up of the section \mbox{$\sigma': P'_A \to A$}}. That is,
  \mbox{$I_\sigma'P''_A$} is a principal ideal and $P''_A$ does not contain a
  proper complete local $A$--subalgebra $S$ such that \mbox{$I_\sigma'S$} is
  principal.
\end{enumerate}
The corresponding (full) subcategory of $\Defsec_{\Rbar\leftarrow P}$
is denoted by \mbox{$\Defes_{{\Rbar} \leftarrow 
    P}$}, and the  
subfunctor of $\uDefsec_{\Rbar\leftarrow P}$ 
of isomorphism classes of equisingular deformations
of the parametrization is denoted by \mbox{$\uDefes_{{\Rbar} \leftarrow
    P}$}.

(2) An object \mbox{$\xi=(\nu_A, \osigma, \sigma)\in
  \Defsec_{\Rbar\leftarrow R} (A)$} is called an {\em equisingular
  deformation of the normalization of $R$} if it is in the image of
  \mbox{$\Defes_{\Rbar\leftarrow P} (A)$} under the natural functor
  \mbox{$\Defsec_{\Rbar\leftarrow
  P}(A)\xleftarrow{\cong}\Defsec_{\Rbar\leftarrow R\leftarrow P} (A)\to
  \Defsec_{\Rbar\leftarrow R} (A)$} given by Proposition \ref{prop:1.3}. The
  corresponding category resp. deformation functor is denoted by
  \mbox{$\Defes_{\Rbar\leftarrow R}$} resp. \mbox{$\uDefes_{\Rbar\leftarrow
  R}$}.
%\end{enumerate}
\end{definition}

 Note that \mbox{$\uDefes_{\Rbar\leftarrow R}\cong
  \uDefes_{\Rbar\leftarrow P}$} by Proposition \ref{prop:1.3}.

\begin{remark}\label{rmk:2.5}
  \begin{Enumerate}
\item[(1)] \ Condition (i), together with the fact that
  $(\varphi'_A,\osigma',\sigma')$ defines a deformation with section of
the parametrization \mbox{$P'\to \Rbar\!\,'$},  implies that the composition
\mbox{$P'_A\to 
    \Rbar_{A,i}\xrightarrow{\osigma_i} A$} is independent of 
  \mbox{$i\in \Lambda_{P'}$}.  If $P''$ is an infinitely near point of $P'$,
  the image of \mbox{$I_{\sigma'}$} in $P''_A$ is  contained in $I_{\sigma''}$

\item[(2)] \ If \mbox{$\osigma_i(t_i)=a\in \fm_A$}, then replacing $t_i$ by
  \mbox{$t_i-a$} trivializes the section $\osigma_i$. In the same way, each
  section $\sigma'$ can be trivialized by choosing appropriate generators $u,v$
  for \mbox{$I_{\sigma'}\subset P'$} (which is possible
  due to Nakayama's lemma).  This choice
  corresponds to the choice of an 
  isomorphism, \mbox{$P'_A\cong A[[u,v]]$} (up to reordering the
  indeterminates). 

\item[(3)] \ Let $P''$ be consecutive to $P'$. Then the formal blow-up
  \mbox{$P'_A\to P''_A$} of the section  \mbox{$\sigma': P'_A \to A$} can be
  expressed in the same way as the formal blow-up of the maximal ideal of a point in Definition
  \ref{def:2.3}:  

\quad According to (2), we may assume that \mbox{$P'_A=A[[u,v]]$}
  and  \mbox{$P''_A=A[[w,z]]$} and that the sections $\sigma',\sigma''$ are
  trivial.  
Additionally, we assume that \mbox{$m(u)\leq m(v)$}, where $m$
    denotes the $\langle w,z\rangle$-order of the residue mod $\fm_A$ of
    elements in $A[[z,w]]$.  Since \mbox{$P'_A\to P''_A$}
  is a formal blow-up, one has 
  \mbox{$\langle u,v\rangle A[[w,z]]=\langle h\rangle A[[w,z]]$} for some
  \mbox{$h\in A[[w,z]]$}. Therefore, \mbox{$u=hp$}, \mbox{$v=hq$}, and
    \mbox{$h=ru+sv$} for some \mbox{$p,q,r,s\in A[[w,z]]$}. Hence, \mbox{$m(u)= m(h)+m(p)\geq
      m(u)+\nu(p)$}, which implies \mbox{$m(p)=0$}, that is, $p$ is a unit
    in $A[[w,z]]$. Replacing $h$ by \mbox{$hp=u$}, we get
    \mbox{$v=u(qp^{-1})=u(v'\!+\!\!\:\beta)$} for a unique \mbox{$v'\in \langle
      w,z\rangle_A$} and \mbox{$\beta\in A$}, and thus \mbox{$\langle
      u,v\rangle A[[u,v']]=\langle u\rangle A[[u,v']]$}. The minimality
    condition for the formal blow-up implies now that
    \mbox{$A[[u,v']]=A[[w,z]]$} and, therefore, \mbox{$\langle
      u,v'\rangle=\langle w,z\rangle=I_{\sigma''}$}. Thus, we may choose 
\mbox{$w=u$} and \mbox{$z=\frac{v}{u}-\beta$}, that is, the formal blow-up of $\sigma'$ is given by \mbox{$A[[u,v]]\to
  A[[w,z]]$}, \mbox{$u\mapsto w$}, \mbox{$v\mapsto w(z\!+\!\beta)$} for a
  unique \mbox{$\beta\in A$}.  
\item[(4)]\ Although the set of infinitely near points $P'$ on $R$ is infinite
  we need to consider in Definition \ref{def:2.6} only the subsest of Ess$(R)$
  consisting of those $P'$ for which the strict transform $R'$ of $R$ is
  singular (which is finite since Ess$(R)$ is finite).
\end{Enumerate}
\end{remark}

%\noindent
Before giving a proof for the existence of a
semiuniversal deformation for equisingular deformations of
the parametrization, we consider (versal) equimultiple deformations.

Recall the notion of
versality.   A deformation $\xi$ over some base $B$ is {\em versal} if
the following holds:  given two deformations $\widetilde{\eta}$ and $\eta$ over
$\widetilde{C}$, 
respectively $C$, such that $\eta$ is induced from $\widetilde{\eta}$ by a
surjective morphism \mbox{$\chi: \widetilde{C} \twoheadrightarrow C$} and from
$\xi$ by a morphism \mbox{$\psi:B\to C$}. Then $\widetilde{\eta}$ can be
induced from $\xi$ by a morphism 
\mbox{$\widetilde{\psi} : B \to \widetilde{C}$} satisfying \mbox{$ \chi
  \circ \widetilde{\psi}=\psi$}. 
$\xi$ is called {\em semiuniversal} if, moreover, the tangent map of
$\widetilde{\psi}$ is uniquely determined.

%\begin{theorem}\label{theo3.6}
%  Any reduced plane curve singularity  over an algebraically closed
%  field of arbitrary characteristic  has a semiuniversal equisingular
%  deformation of the parametrization with smooth base space.
%\end{theorem}
%
%\noindent
We shall identify explicitly a
semiuniversal deformation for \mbox{$\uDefes_{{\Rbar} \leftarrow P}$} 
  as a subfamily of a semiuniversal deformation for
\mbox{$\uDefsec_{{\Rbar} \leftarrow P}$}.

We begin by describing a semiuniversal deformation for
\mbox{$\uDefsec_{{\Rbar} \leftarrow P}$}. We do this in a slightly more general
context: 
Given an integer vector \mbox{$\bm = (m_1, \dots,
  m_r)$} such that either \mbox{$\bm=\bnull=(0, \ldots, 0)$}, or \mbox{$1 \le 
  m_i \le \ord \varphi_i$}, we
call a deformation of the parametrization over $A$,
\mbox{$(P_A\!\xrightarrow{\varphi_A}\!\Rbar_A,\osigma,\sigma)$}, 
an {\it $\bm$-multiple deformation\/} if 
\[
\varphi_A (I_\sigma) \subset I_{{\osigma}}^{\bm} :=
  \bigoplus^{r}_{i=1} I^{m_i}_{{\osigma}_i}\,.
\]
Here, \mbox{$I_{{\osigma_i}} := \ker {\osigma_i}\subset \Rbar_{A,i}$} and
\mbox{$I_\sigma = \ker {\sigma}\subset P_A$}.  The corresponding category
is denoted by \mbox{$\Def^{\bm}_{\Rbar \leftarrow P}(A)$}, and the deformation
functor by \mbox{$\uDef^{\bm}_{\Rbar
    \leftarrow P}$}. 

Note that \mbox{$\Def^{\bm}_{\Rbar \leftarrow P}$} coincides with
\mbox{$\Def_{\Rbar \leftarrow P}$} for \mbox{$\bm=\bnull=(0,\dots,0)$}, with 
\mbox{$\Defsec_{\Rbar
    \leftarrow P}$} for \mbox{$\bm=\bone=(1,\dots,1)$}, and with 
\mbox{$\Defem_{\Rbar \leftarrow P}$} for
\mbox{$\bm=(\ord \varphi_1,\dots, \ord \varphi_r)$}.

The following proposition gives an explicit description of the tangent space 
\mbox{$T^{1,\bm}_{\Rbar \leftarrow P}:=\uDef^{\bm}_{\Rbar
    \leftarrow P}(\Keps)$}  to the functor
of $\bm$-multiple deformations of the parametrization. We start by introducing
the notations
$$
\dot{\bx}:= 
%\left(\dot{x}_1,\dots,\dot{x}_r\right)
\left\lgroup
\begin{matrix}
\dot{x}_1\\[-0.2em]
:\\[-0.2em]
\dot{x}_r
\end{matrix}
\!\right\rgroup\!
\,, \quad 
\dot{\by}:= 
%\left(\dot{y}_1,\dots,\dot{y}_r\right)
\left\lgroup
\begin{matrix}
\dot{y}_1\\[-0.2em]
:\\[-0.2em]
\dot{y}_r
\end{matrix}
\!\right\rgroup \quad\raisebox{+1ex}{,}
$$
where $\dot{x}_i,\dot{y}_i$ denote the derivatives of $x_i,y_i$ with
respect to $t_i$. Moreover, we write \mbox{$\fm$} for the maximal
ideal of $R$, \mbox{$\fmbar$} for the Jacobson radical of $\Rbar$, and
$$ \fmbar^{\bm} := \bigoplus_{i=1}^r  t_i^{m_i}K[[t_i]]
\subset \bigoplus_{i=1}^r K[[t_i]]
 = \Rbar\, .$$

We define the following $R$---modules:

\begin{eqnarray*}
  M^{\bnull}_{{\Rbar} \leftarrow P} & :=& 
        \left(\Rbar\oplus\Rbar\right)/\left(\Rbar\cdot
          (\dot{\bx},\dot{\by})+(R\oplus R)\right)\,, \\
  M^{\bm}_{{\Rbar} \leftarrow P} & :=& 
        \left(\fmbar^{\bm}\!\oplus\fmbar^{\bm} \right)/\left(\fmbar\cdot
          (\dot{\bx},\dot{\by})+(\fm \oplus
        \fm)\right)\,, \quad 1\leq m_i\leq \ord \varphi_i\,.
\end{eqnarray*}

%In particular, 
%$\bm$-multiple deformations coincide with deformations with
%section if \mbox{$\bm=(1,\dots,1)$}, respectively with equimultiple
%deformations if \mbox{$m_i = \min\{\ord x_i(t_i), \ord y_i(t_i)\}$} for all
%$i$. 
%
%Let 
%\mbox{$T^{1,m}_{{\Rbar} \leftarrow P}$} denote the corresponding tangent
%space, that is, $m$--multiple deformations over the dual numbers
%$K[\varepsilon]$. 

  \begin{proposition}\label{prop3.7} With the above notations, we have
\begin{Enumerate}
\item[(1)] \  
$ T^{1,\bm}_{{\Rbar} \leftarrow P} \cong M^{\bm}_{{\Rbar} \leftarrow P} 
$.        
\item[(2)] \ 
Let \mbox{$\varphi_A:A[[x,y]]\to \bigoplus_{i=1}^r A[[t_i]]$} define an 
        $\bm$-multiple deformation of the parametrization (with trivial
        sections) over \mbox{$A=K[[\bT]]=K[[T_1, \dots, T_k]]$}, given by
        power series $X_i(\bT,t_i),\,Y_i(\bT,t_i)$, \mbox{$i=1,\dots,r$}. 
        Then this deformation is a versal (respectively semiuniversal)
$\bm$-multiple deformation  iff the column vectors
$$ \left(\frac{\partial X_i}{\partial T_j}(\bnull,t_i),\frac{\partial
    Y_i}{\partial T_j}(\bnull,t_i)\right)_{i=1,\dots,r} \in
\fmbar^{\bm}\! 
\oplus \fmbar^{\bm}\,,\quad j=1,\dots,k\,,$$
represent a system of generators (respectively a basis) for the vector space
$ M^{\bm}_{{\Rbar} \leftarrow P}$.  

\item[(3)] \
Let \mbox{$\ba^{(j)},\bb^{(j)}\in \fmbar^{\bm}=
      \bigoplus_{i=1}^r t_i^{m_i}K[[t_i]]$} be such that 
    $$ (\ba^{(j)},\bb^{(j)}) =\left(
    \left\lgroup 
\begin{matrix}
a_1^{(j)} \\[-0.2em]
 : \\[-0.2em]
a_r^{(j)} 
\end{matrix}
\right\rgroup ,
\left\lgroup
\begin{matrix}
b_1^{(j)} \\[-0.2em]
 : \\[-0.2em]
b_r^{(j)} 
\end{matrix}
\right\rgroup\right)\,,
\quad  j=1,\dots,k\,,$$ 
represent a basis (resp. a system of generators) for
    $M^{\bm}_{{\Rbar} \leftarrow P}$. Then the deformation (with trivial
    sections) of the
    parametrization over \mbox{$K[[\bT]]=K[[T_1,\dots,T_k]]$}  defined by  
%\mbox{$\varphi_A=(\varphi_{A,1}, \dots,\varphi_{A,r}) :A[[x,y]]\to
%  \bigoplus_{i=1}^r t_i^{m_i} A\{t_i\}$}, 
    \begin{eqnarray*}
%\varphi_{A,i}(x) =  
X_i(\bT,t_i) & = & x_i(t_i) + \sum_{j=1}^k a_i^{(j)}(t_i) \cdot T_j\,,\\
%\varphi_{A,i}(y) = 
Y_i(\bT,t_i) & = & y_i(t_i) + \sum_{j=1}^k b_i^{(j)}(t_i) \cdot T_j\,,
    \end{eqnarray*}
\mbox{$i=1,\dots,r$}, is a semiuniversal (resp. versal) $\bm$-multiple deformation of the
parame\-tri\-za\-tion.

\quad In particular, $\bm$-multiple
deformations of the parametrization are unobstructed, and they have a smooth
semiuniversal base space of dimension $\dim_K ( M^{\bm}_{{\Rbar} \leftarrow
  P})$.  
\end{Enumerate}
\end{proposition}

%\noindent
We omit the proof, since it is similar to (but simpler than) the proof of the
analogous statement for equisingular deformations (Theorem \ref{theo3.8}).

\medskip
%\noindent
At the end of this paragraph, we consider the problem of the uniqueness of the
sections. This depends on the characteristic of $K$. 

\begin{lemma}\label{lem:unique section}
If the characteristic
    of $K$ is good, then 
the forgetful natural transformation \mbox{$\Defem_{\Rbar\leftarrow P}\to 
  \Def_{\Rbar\leftarrow P}$} %(neglecting the sections) 
is injective in the following cases: 
\begin{quote}
  \begin{enumerate}
  \item[{\it Case\;1.}] All branches of $R$ are singular. 
  \item[{\it Case\;2.}] The branches of $R$ have pairwise
    different tangent directions.
  \end{enumerate}
\end{quote}
In other words, in these cases the equimultiple sections $\osigma$ and
$\sigma$ are uniquely determined (if they exist).  
\end{lemma}

\begin{proof}
 By \cite[1.1.4]{Wahl1}, it suffices to show the injectivity on the tangent
 level, that is, for \mbox{$A=\Keps$}.  
Assume first that $R$ is
irreducible,  \mbox{$\mult(R)=m\geq 2$}. Assume further that
\mbox{$X=x(t)+\varepsilon x_{\varepsilon}(t),Y=y(t)+\varepsilon
  y_{\varepsilon}(t)$}, \mbox{$ x_{\varepsilon},y_{\varepsilon}\in K[[t]]$},
 defines an em-deformation along the trivial sections, that is,
 \mbox{$\min\{\ord_t( x_{\varepsilon}), \ord_t( x_{\varepsilon})\} \geq m$}. If
 $X,Y$ 
 defines also an em-deformation with the
 sections $\osigma,\sigma$ given by \mbox{$I_{\osigma}=\langle
   t+\varepsilon \gamma_0\rangle$}, \mbox{$I_{\sigma}=\langle
   x+\varepsilon \alpha_0, y+\varepsilon \beta_0\rangle$},
 \mbox{$\alpha_0,\beta_0,\gamma_0\in K$}, then the equimultiplicity
 condition \mbox{$\varphi_A(I_{\sigma})\subset I^{m}_{\osigma}$} is
 equivalent to 
\setlength{\arraycolsep}{2pt}
\begin{eqnarray*}
   x(t)+\varepsilon x_{\varepsilon}(t) +\varepsilon \alpha_0 &=&
   \bigl(h(t)+\varepsilon h_{\varepsilon}(t)\bigr)\cdot (t^m+m\varepsilon 
   \gamma_0t^{m-1})\\
 y(t)+\varepsilon y_{\varepsilon}(t) +\varepsilon
   \beta_0 &=& \bigl(k(t)+\varepsilon k_{\varepsilon}(t)\bigr)\cdot
   (t^m+m\varepsilon 
   \gamma_0t^{m-1})\,,
 \end{eqnarray*}
for some \mbox{$h,h_{\varepsilon},k,k_{\varepsilon}\in K[[t]]$}. Comparing
coefficients, this implies \mbox{$\alpha_0=\beta_0=0$} and, if the
characteristic is either $0$ or if it does not divide $m$, then
\mbox{$\gamma_0=0$}, too. 

The injectivity in Case 1 follows immediately from these considerations, since
an em-deformation of the parametrization of $R$ induces by definition
em-deformations of the parametrizations of the branches $R_i$.

In Case 2 we may assume 
that at least one of the branches, say $R_1$ is non-singular (otherwise, Case 1
applies). Then, for each fixed
\mbox{$j\in \{2,\dots,r\}$}, we may choose coordinates $x,y$ such that $R_1$ is
parametrized by \mbox{$(x_1,y_1)=(t_1,0)$} and $R_j$ has the tangent direction
\mbox{$x=0$}. Let $X,Y$ define an em-deformation with trivial sections, and
also an em-deformation with the sections
$\osigma,\sigma$ given by \mbox{$I_{\osigma_i}=\langle
   t+\varepsilon \gamma_{i}\rangle$}, \mbox{$I_{\sigma}=\langle
   x+\varepsilon \alpha_0, y+\varepsilon \beta_0\rangle$},
 \mbox{$\alpha_0,\beta_0,\gamma_i\in K$}. Then, similar to the above,
\setlength{\arraycolsep}{2pt}
\begin{eqnarray*}
  t+\varepsilon t a_{1,\varepsilon}(t) +\varepsilon \alpha_0 &=&
   \bigl(h_1(t)+\varepsilon h_{1,\varepsilon}(t)\bigr)\cdot (t+\varepsilon 
   \gamma_1)\\
 0+\varepsilon  t b_{1,\varepsilon}(t) +\varepsilon
   \beta_0 &=& \bigl(k_1(t)+\varepsilon k_{1,\varepsilon}(t)\bigr)\cdot
   (t+\varepsilon \gamma_0)\,,
 \end{eqnarray*}
 which implies \mbox{$\alpha_0=\gamma_1$} and \mbox{$\beta_0=0$}. If the branch
 $R_j$ is singular, Case 1 shows that \mbox{$\gamma_j=\alpha_0=0$}, hence the
 uniqueness. Thus, we can assume that $R_j$ is smooth and parametrized by 
\mbox{$(x_j,y_j)=(0,t_j)$}. Then the same reasoning as above gives
\mbox{$\beta_0=\gamma_j$} and \mbox{$\alpha_0=0$}, thus the uniqueness of the
sections.   
\end{proof}

%\noindent
However, there are examples of em-deformations such that both 
sections $\osigma$ and $\sigma$ are {\em not\/} unique:
%
%For \mbox{$\Char(K)=p>0$}, however, there are deformations of irreducible
%singularities $R$, \mbox{$\mult(R)=p$}, such that the equimultiple section
%\mbox{$\sigma:A[[x,y]]\to A$} is {\em not\/} unique:

\begin{example}
\begin{Enumerate}
\item[(1)] \
 Let \mbox{$\Char(K)=p>0$}, and consider the irreducible
 singularity \mbox{$R=K[[x,y]]/\langle y^p\!-x^{2p+1}\rangle$}. Then 
  \mbox{$(X,Y)=\bigl(t^p,\,t^{p}(t-\alpha)^{p+1}\bigr)$}, \mbox{$\alpha\in \fm_A$}, 
  defines an em-deformation of the parametrization with trivial
  sections over $A$. But, it also defines an em-deformation with the
  sections $\osigma,\sigma$ given by \mbox{$I_{\osigma}=\langle
    t-\alpha\rangle$}, 
  \mbox{$I_{\sigma}=\langle x\!\!\:-\!\!\:\alpha^p,y \rangle$}.  

\item[(2)] \ Let %\mbox{$\Char(K)\neq 2$}, and let 
  \mbox{$R=K[[x,y]]/\langle 
  x^4\!-y^2\rangle$}, which decomposes into two smooth branches with the same
tangent direction. Then, for each \mbox{$\alpha\in \fm_A$},
 \mbox{$(X_1,Y_1)=(t_1,\,-t_1^2+\alpha t_1)$},
\mbox{$(X_2,Y_2)=(t_2, \,t_2^2\!-\alpha t_2)$}, 
defines an em-deformation with trivial 
sections. And, it defines an em-deformation with sections $\osigma,\sigma$
given by  \mbox{$I_{\osigma_1}=\langle t_1-\alpha
    \rangle$}, \mbox{$I_{\osigma_2}=\langle t_2-\alpha
    \rangle$}, \mbox{$I_\sigma=\langle x-\alpha ,y\rangle$}. 
\end{Enumerate}

\medskip%\noindent
Note that none of these examples defines an equisingular deformation. Indeed,
after formally blowing up the sections, we do not get an equimultiple
deformation of the strict transform. In the second example, this is caused by
the fact that the sections do not satisfy the compatibility condition of
Definition \ref{def:2.6} (iii).
%After the formal blow-up given by \mbox{$x\mapsto u$},
%\mbox{$y\mapsto uv$} of the trivial section, we obtain
%\mbox{$U_1=t_1$}, \mbox{$V_1=-t_1+2\alpha$}, \mbox{$U_2=t_2$},
%\mbox{$V_2=-t_2-2\alpha$}. For both branches individually, we find a unique
%section $\sigma'$: for the first branch, we have \mbox{$I_{\sigma'}=\langle
%  v-2\alpha\rangle$}, for the second \mbox{$I_{\sigma'}=\langle
%  v+2\alpha\rangle$}. Thus, if \mbox{$\Char(K)\neq 2$}, there is no section
%$\sigma'$ as above. In particular, $\varphi_A$ defines an equimultiple
%deformation with trivial sections of the parametrization which is not
%equisingular. 
\end{example}

\begin{proposition}\label{prop:sections unique}
If the characteristic of $K$ is good, and if $R$ is singular, then the 
forgetful functor \mbox{$\Defes_{\Rbar\leftarrow P}\!\to
  \Def_{\Rbar\leftarrow P}$} is injective. That is, the equisingular
sections $\osigma$ and $\sigma$ are unique.
\end{proposition}

\begin{proof}
If either all branches of $R$ are singular, or if they have pairwise
different tangent directions, the uniqueness of the sections is implied by
Lemma \ref{lem:unique section}. Moreover, the proof of Lemma
\ref{lem:unique section} shows that the section $\sigma$ is uniquely
determined if $R$ either has a singular branch $R_j$ or two smooth
branches $R_i,R_j$, \mbox{$i\neq j$}, intersecting transversally. Moreover, the
proof shows that if the branch $R_j$ is either singular, or smooth and
transversal to some other branch, the section $\osigma_j$ is uniquely
determined. It remains to consider smooth branches that are tangential to all 
of the other branches. After finitely many formal blow-ups the strict transform
of such a branch $R_i$ becomes transversal to one of the other branches. Thus,
the section $\osigma_i$ is uniquely defined, and hence also $\sigma$ (see the
proof of Lemma \ref{lem:unique section}).
\end{proof}

%\noindent
There are examples of equisingular deformations in bad characteristic
\mbox{$p>0$} where $\osigma$ is not unique. For
instance, the trivial deformation over $\Keps$ of the parametrization 
\mbox{$(t^{p+1},t^p)$} of \mbox{$R=K[[x,y]]/\langle
  x^p-y^{p+1}\rangle$} is equisingular along the trivial sections and along the
section defined by \mbox{$I_{\osigma}=\langle t+\varepsilon \rangle$}.  
However, the section $\sigma$ is always unique:

\begin{proposition}\label{prop:sigma unique}
Let \mbox{$A\in \sA_K$}, and let 
\mbox{$\varphi_A:A[[x,y]]\to \bigoplus_{i=1}^r A[[t_i]]$} 
define a deformation of \mbox{$P\to \Rbar$}. If $R$ is singular then for at most one 
\mbox{$\sigma:A[[x,y]]\to A$} there is a lifting $\osigma$ such that
$(\varphi_A,\osigma,\sigma)$ is equisingular.
\end{proposition}

\begin{proof}
From the considerations in the proofs of Lemma \ref{lem:unique section} and
Proposition \ref{prop:sections unique}, it is clear that it suffices to
consider  
the case of an irreducible singularity $R$ such that the multiplicity of $R$
and of all of its singular strict transforms are divisible by 
\mbox{$p=\Char(K)>0$}.  Indeed, it suffices to consider the last singular
strict transform of $R$. It has a parametrization
\mbox{$(x,y)=(t^{kp},t^{kp+1})+\,\text{higher terms}$},  
 \mbox{$k\geq 1$}. 

Let $X,Y$ define an em-deformation of $(x,y)$ with trivial sections
 over \mbox{$A\in \sA_K$}. 
Then, up to terms of $t$-order \mbox{$kp+1$}, respectively \mbox{$kp+2$}, we
have \mbox{$X(t)=(1+a)t^{kp}$}, \mbox{$Y(t)=b_1t^{kp}+(1+b_2)t^{kp+1}$} for
some 
\mbox{$a,b_1,b_2\in \fm_A$}. If $X,Y$ is also equimultiple along the sections
$\osigma, \,\sigma$ defined by  \mbox{$I_{\osigma}=\langle
   t+\gamma\rangle$}, \mbox{$I_{\sigma}=\langle
   x+ \alpha, y+ \beta\rangle$},
 \mbox{$\alpha,\beta,\gamma\in \fm_A$}, then we get (again up to terms of
 $t$-order  \mbox{$kp+1$}, respectively \mbox{$kp+2$})  
\setlength{\arraycolsep}{2pt}
\begin{eqnarray*}
   (1+a)t^{kp} + \alpha &=&
   \bigl(c_0+c_1t +\ldots+c_{kp}t^{kp}\bigr)\cdot (t^{p}+\gamma^{p})^k\\
   b_1t^{kp}+(1+b_2)t^{kp+1} + \beta &=&
   \bigl(d_0+d_1t +\ldots+d_{kp+1}t^{kp+1}\bigr)\cdot (t^{p}+\gamma^{p})^k\,,
 \end{eqnarray*}
for some \mbox{$c_j,d_j\in A$}. Comparing coefficients, we get 
\mbox{$\alpha=c_0\gamma^{kp}$}, \mbox{$\beta=d_0\gamma^{kp}$}, and the
conditions \mbox{$0=d_1\gamma^{kp}$}, \mbox{$1+b_2 =d_1+e \gamma^{p}$}
for some \mbox{$e\in A$}. 

Thus, \mbox{$0=(1+b_2-e\gamma^{p})\cdot \gamma^{kp}$}, which implies
\mbox{$\gamma^{kp}=0$}. Together with the above equalities, this yields
\mbox{$\alpha=\beta=0$} as claimed.
\end{proof}

\section{Versal Equisingular Deformations of the Parametrization}\label{sec:section 3}

In this section, we give a proof for the existence of a semiuniversal
equisingular deformation of the parametrization and show that it has an
algebraic representative with a smooth base. Moreover, we show that
equisingular versality is an open property.

Let \mbox{$\Ties_{{\Rbar} \leftarrow P}=\uDefes_{{\Rbar} \leftarrow R}(\Keps)$}
denote the tangent space to the functor \mbox{$\uDefes_{{\Rbar} \leftarrow
    R}$}.   It is a subspace of 
\mbox{$\Tisec_{{\Rbar} \leftarrow P} = \left(\fmbar \oplus
  \fmbar\right)/\left(\fmbar(\dot{\bx},\dot{\by}) + (\fm \oplus
  \fm)\right)$} 
(indeed, it is a subspace  of each \mbox{$T^{1,\bm}_{{\Rbar} \leftarrow P}$}
where \mbox{$1 \le m_i \le \ord \varphi_i$}).   Hence,
\[
  \Ties_{{\Rbar} \leftarrow P}  = \dfrac{\Ies_{{\Rbar} \leftarrow
      P}}{\fmbar \cdot(\dot{\bx},\dot{\by}) + (\fm \oplus \fm)}
\]
with 
\[
\Ies_{{\Rbar} \leftarrow P} := \left\{
(\ba,\bb) \in \fmbar\oplus \fmbar \:\left|\:
  \begin{array}{l}
\bigl\{(x_i(t_i) + \varepsilon a_i, y_i(t_i) + \varepsilon b_i) \:\big|\:
i=1,\dots,r \bigr\} \text{
   defines}\\
\text{an equisingular deformation of the parametriza-}\\
\text{tion }\varphi \text{ with trivial sections over } K[\varepsilon] 
  \end{array}\!
\right.\right\}.
\]
We call \mbox{$\Ies_{{\Rbar} \leftarrow P}$} the {\em equisingularity
  module of \mbox{$P\to \Rbar$}}. Below, we show that it is an
$R$-module (Corollary \ref{cor:Ies is submodule}).

The main theorem of this section states now that $\Defes_{{\Rbar} \leftarrow P}$ is a
``linear'' subfunctor of $\Defsec_{{\Rbar} \leftarrow P}$.
 As such it is already completely determined by its tangent
space:

\begin{theorem}\label{theo3.8}
  \begin{Enumerate}
  \item[(1)] \ Let $K[[\bT]] = K[[T_1, \dots, T_k]]$. Then an $r$-tuple of
    power 
    series \mbox{$X_i(\bT,t_i), Y_i(\bT,t_i)\in K[[\bT,t_i]]$},
    \mbox{$i=1,\dots,r$}, satisfying  
    \begin{align*}
      X_i(\bT,t_i) & \equiv x_i(t_i) + \sum_{j=1}^k a_i^{(j)}(t_i) \cdot T_j 
      \mod  \langle \bT\rangle^2 \,,\quad
      X_i(\bT,0)=0\\ 
      Y_i(\bT,t_i) & \equiv y_i(t_i) + \sum_{j=1}^k b_i^{(j)}(t_i) \cdot
      T_j \mod  \langle \bT\rangle^2\,,\quad Y_i(\bT,0)=0 \,,
    \end{align*}
defines an equisingular deformation with trivial
sections of the parametrization over $K[[\bT]]$ iff \mbox{$(\ba^{(j)},
  \bb^{(j)}) \in \Ies_{{\Rbar} \leftarrow P}$} for all \mbox{$j=1,\dots,k$}.

\item[(2)] \ Let $X_i(\bT,t_i),\,Y_i(\bT,t_i)$, \mbox{$i=1,\dots,r$}, define an
  equisingular deformation with trivial sections of the parametrization over
  $K[[\bT]]$.  Then this deformation is a 
  versal (respectively semiuniversal) object of $\uDefes_{\Rbar\leftarrow P}$
  iff 
$$
\left(\left\lgroup
\begin{matrix}
\tfrac{\partial X_1}{\partial T_j} (\bnull,t_1) \\[-0.2em]
:\\[-0.2em]
\tfrac{\partial X_r}{\partial T_j} (\bnull,t_r)
\end{matrix}
\!\right\rgroup ,
%\oplus
\left\lgroup
\begin{matrix}
\tfrac{\partial Y_1}{\partial T_j} (\bnull,t_1) \\[-0.2em]
:\\[-0.2em]
\tfrac{\partial Y_r}{\partial T_j} (\bnull,t_r)
\end{matrix}
\!\right\rgroup\right)
\,, \quad j=1,\dots,k\,,
$$
 represent a system of $K$-generators (respectively a $K$-basis) of
 $\Ties_{{\Rbar} \leftarrow P}$. 

\item[(3)] \  Let \mbox{$(\ba^{(j)},\bb^{(j)})\in \Ies_{\Rbar\leftarrow P}$}, 
\mbox{$j=1,\dots,k$}, represent a basis (respectively a system of generators)
of \mbox{$\Ties_{{\Rbar} \leftarrow P}$}. Then
\begin{eqnarray*}
 X_i(\bT,t_i)&=&x_i(t_i)+\sum_{j=1}^k a_i^{(j)}(t_i)\cdot T_j\,, \\   
 Y_i(\bT,t_i)&=&y_i(t_i)+\sum_{j=1}^k b_i^{(j)}(t_i)\cdot T_j\,, 
\end{eqnarray*}
\mbox{$i=1,\dots,r$}, define 
a semiuniversal (respectively versal) equisingular deformation with trivial
sections of the parametrization over $K[[\bT]]$. In particular, equisingular
deformations of the parametrization are unobstructed, and the semiuniversal
deformation has a smooth base of dimension 
\mbox{$\dim_K \Ties_{{\Rbar} \leftarrow P}$}.
\end{Enumerate}
\end{theorem}

%\noindent
 The proof of this theorem needs some preparations. It is based on considering
 {\em small extensions\/} in $\sA_K$, 
 that is, surjective morphisms \mbox{$\widetilde{A} \twoheadrightarrow A$} of
 Noetherian complete local $K$-algebras with a one-dimensional kernel. The
 generator of this kernel will be usually denoted by $\varepsilon$. Then, as
 $K$-vector spaces, \mbox{$\widetilde{A} = A \oplus \varepsilon K$}, and  
\mbox{$\varepsilon \fm_A = 0$}.  

In the following Proposition \ref{prop:2.9} we show that an em-deformation of
the parametrization \mbox{$P\to \Rbar$} together with a 
factorization through the deformation of an infinitely near point $P'$ can be
lifted to a small extension \mbox{$\widetilde{A}$} of $A$. Moreover, after
trivialization of the section, the lifting to the infinitely near point is
determined by the lifting of the deformation of \mbox{$P\to \Rbar$} to
$\widetilde{A}$. 

We fix some notation: let
\mbox{$\widetilde{A}\twoheadrightarrow A$} be a small extension in $\sA_K$, and let
$\varepsilon$ denote a generator for its kernel. Let 
$$\varphi_A=(\varphi_{A,1},\dots, \varphi_{A,r}): P_A=A[[x,y]] \to
  \bigoplus_{i=1}^r A[[t_i]]=\Rbar_A\,,$$
\mbox{$\varphi_{A,i}(x)=X_i(t_i)$}, \mbox{$\varphi_{A,i}(y)=Y_i(t_i)$},  
define a deformation with trivial sections of the
parametrization 
$$\varphi=(\varphi_{1},\dots, \varphi_{r}):P=K[[x,y]] \to \bigoplus_{i=1}^r
K[[t_i]]=\Rbar\,,$$
\mbox{$\varphi_i(x)=x_i$,  $\varphi_i(x)=y_i\in K[[t_i]]$}.
Moreover, let \mbox{$P'=K[[u,v]]$} be an
infinitely point on $R$ in the first infinitesimal neighbourhood of $P$ such
that \mbox{$\varphi_i$}, \mbox{$i\in \Lambda_{P'}$}, factors as
\mbox{$\varphi_i= \varphi'_{i}\circ \pi'$}, where \mbox{$\pi':K[[x,y]]\to
  K[[u,v]]$} is a formal blow-up 
of the maximal ideal in $P$.
Finally, let \mbox{$\pi'_A:A[[x,y]]\to A[[u,v]]$} be a formal blow-up of the
trivial section in $P_A$ extending $\pi'$. 
After a linear change of variables, we may
assume that \mbox{$\ord_{t_i} x_i\leq \ord_{t_i} y_i$} for all $i$, and 
\mbox{$\pi'_A(x)=u$}, \mbox{$\pi'_A(y)=u(v+\beta)$}, \mbox{$\beta\in A$} (see
Remark \ref{rmk:2.5}\,(5)).  

\begin{proposition} \label{prop:2.9}
 Assume that \mbox{$(\varphi_A,\osigma,\sigma)$} is 
  equimultiple and that the components
  \mbox{$\varphi_{A,i}$}, \mbox{$i\in \Lambda_{P'}$}, of $\varphi_A$ factor as
  \mbox{$\varphi_{A,i}= 
    \varphi'_{A,i}\circ \pi'_A$} such that
 $$\varphi'_A=\bigl(\varphi'_{A,i}\bigr)_{i\in \Lambda_{P'}}:A[[u,v]]\to
 \bigoplus_{i\in \Lambda_{P'}}  A[[t_i]] $$  
defines a deformation of \mbox{$P'\!=K[[u,v]]\to \bigoplus_{i\in \Lambda_{P'}}
    K[[t_i]]=:\Rbar\!\,'$} with trivial sections 
  over $A$. 
Then there exists an extension of $\varphi'_{A}$,
 $$\varphi'_{\widetilde{A}}=\bigl(\varphi'_{\widetilde{A},i}\bigr)_{i\in
   \Lambda_{P'}}:\widetilde{A}[[u,v]]\to \bigoplus_{i\in 
   \Lambda_{P'}}  \widetilde{A}[[t_i]]\,, $$  
and a formal blow-up
\mbox{$\pi'_{\widetilde{A}}:\widetilde{A}[[x,y]]\to \widetilde{A}[[u,v]]$}
of the trivial section $\sigma'$  in $\widetilde{A}[[x,y]]$ extending $\pi'_A$, such that
the following hold:
\begin{enumerate}
\item[(i)] $\varphi'_{\widetilde{A}}$ defines a deformation of \mbox{$P'\to
    \Rbar\!\,'$} with trivial sections over $\widetilde{A}$;
\item[(ii)] \mbox{$\varphi'_{\widetilde{A}}\circ \pi'_{\widetilde{A}}$} defines
  an em-deformation of \mbox{$P\to \Rbar\!\,'$} with trivial
  sections.
\end{enumerate}
Moreover, $\varphi'_{\widetilde{A}}$ and $\pi'_{\widetilde{A}}$ are uniquely
determined by $\varphi'_{A}$, $\pi'_{A}$, and %the composition
\mbox{$\varphi'_{\widetilde{A}}\circ 
  \pi'_{\widetilde{A}}$}.
 \end{proposition}

%\noindent
In other words, the proposition states that the following diagram of solid
arrows (and trivial sections) can be completed by the dotted arrows 
$$
\UseComputerModernTips
\xymatrix@C=9pt@R=1pt@M=6pt{
{\Rbar\!\,'}\ar@{}[ddrr]|{\Box} && 
{\Rbar}\!\,'_A \ar@{->>}[ll]\ar@{}[ddrr]|{\Box}&& 
{\Rbar}\!\,'_{\widetilde{A}}\ar@{->>}[ll]\ar@/^4pc/@<10pt>[dddddd]^{\osigma'}
\\   
\\
P' \ar[uu]^{\varphi}\ar@{}[ddrr]|{\Box} && \ar@{->>}[ll] P'_A
\ar[uu]_{\varphi_A} \ar@{}[ddrr]|{\Box}\ar@/^2pc/[dddd]^(0.4){\sigma'}&& \ar@{->>}[ll]
P'_{\widetilde{A}}\ar@{.>}[uu]^{\varphi_{\widetilde{A}}} \ar@/^2pc/@<6pt>[dddd]^{\widetilde{\sigma}'}
\\
 \\
P \ar[uu]^{\varphi}\ar@{}[ddrr]|{\Box} &&  P_A \ar@{->>}[ll]
\ar[uu]_{\pi'_A} \ar@/^0.51pc/@<1pt>[dd]^-{\sigma}\ar@{}[ddrr]|{\Box} && 
P_{\widetilde{A}}\ar@{->>}[ll]\ar@{.>}[uu]^{\pi'_{\widetilde{A}}}\ar@/^1pc/@<2pt>[dd]^(0.4){\widetilde{\sigma}} 
 \\ 
 \\
 K \ar@{^{(}->}[uu] && \ar@{->>}[ll] A\ar[uu] && \ar@{->>}[ll] {\widetilde{A}}
\ar[uu]
}
$$

and that the dotted arrows are uniquely determined by their composition.

\begin{proof} {\it Step 1. Uniqueness.} 
Assume we have extensions $\varphi'_{\widetilde{A}}$, $\pi'_{\widetilde{A}}$ of
$\varphi'_A$, $\pi'_A$ as in the proposition. Setting 
$$U_i:=\varphi'_{A,i}(u)\,,\ V_i:=\varphi'_{A,i}(v)\,,\ 
\widetilde{U}_i:=\varphi'_{\widetilde{A},i}(u)\,,\
\widetilde{V}_i:=\varphi'_{\widetilde{A},i}(v)\,,$$  
we have
\mbox{$\widetilde{U}_i=U_i +
  \varepsilon a'_i$}, \mbox{$\widetilde{V}_i=V_i +
  \varepsilon b'_i$} for \mbox{$a'_i,b'_i\in t_iK[[t_i]]$}, and
\mbox{$\pi'_{\widetilde{A}}(x)=u+\varepsilon \alpha_{\varepsilon}$},
\mbox{$\pi'_{\widetilde{A}}(y)=(u+\varepsilon
  \alpha_{\varepsilon})(v+\beta+\varepsilon\beta_{\varepsilon})$}  
where \mbox{$\alpha_{\varepsilon},\beta_{\varepsilon}\in K$}, with
\mbox{$\alpha_\varepsilon=\beta_\varepsilon=0$} if $\widetilde{\sigma}$ is
trivial. Here \mbox{$\pi'_A(x)=u,\ \pi'_A(y)=u(v+\beta),\ \beta\in A$} by Remark
2.7 (2).  
 By assumption, $\varphi_{A,i}$ factors as \mbox{$\varphi_{A,i}=
   \varphi'_{A,i}\circ \pi'_A$}. Hence,
\mbox{$X_i=\varphi_{A,i}(x)=U_i$},\, \mbox{$Y_i =\varphi_{A,i}(y)= U_i\cdot
\bigl(V_i+\beta\bigr)$}.

Setting \mbox{$\widetilde{X}_i:=\varphi'_{\widetilde{A},i}\circ
  \pi'_{\widetilde{A}}(x)$}, and  
\mbox{$\widetilde{Y}_i:=\varphi'_{\widetilde{A},i}\circ
  \pi'_{\widetilde{A}}(y)$}, we get from (ii) that
$$ 
\widetilde{X}_i=X_i+\varepsilon a_i\,,\quad 
\widetilde{Y}_i=Y_i+\varepsilon b_i
$$
for some \mbox{$a_i,b_i\in t_iK[[t_i]]$} and 
$$
\widetilde{X}_i =  U_i +\varepsilon a'_i+\varepsilon
\alpha_{\varepsilon}\,,\quad 
\widetilde{Y}_i = 
\left(U_i+\varepsilon a'_i+\varepsilon
\alpha_{\varepsilon}\right)\cdot  
\left(V_i+\varepsilon b'_i+\beta+\varepsilon\beta_{\varepsilon}\right)
\,.
$$
Comparing coefficients, this implies 
$$ a_i=a'_i+\alpha_{\varepsilon}\,, \quad
b_i=(a'_i+\alpha_{\varepsilon})(v_i+\beta_0)+(b_i'+\beta_{\varepsilon})u_i\,,$$
where \mbox{$u_i\!\!\:=\!\!\:x_i,\,v_i\in t_iK[[t_i]]$}, resp.\ 
\mbox{$\beta_0\in K$} are the residues of $U_i,V_i$, resp.\ $\beta$ mod
$\fm_A$.   
Since \mbox{$a_i(0)=a'_i(0)=0$} this implies \mbox{$\alpha_{\varepsilon}=0$},
hence 
\begin{equation}\label{eqn:a'i,b'i}
a'_i=a_i\,,
\quad 
b'_i= \frac{b_i-a_i(v_i+\beta_0)}{x_i} -\beta_{\varepsilon}\,.
\end{equation}
Since \mbox{$b'_i(0)=0$}, this implies that
\begin{equation}\label{eqn:beta0}
 \beta_{\varepsilon} = \left(\frac{b_i-a_i(v_i+\beta_0)}{x_i}\right)(0) =\frac{b_i}{x_i}(0)+\frac{a_i\beta_0}{x_i}(0)\,. 
\end{equation}
In particular, the expression on the right-hand side does not depend on the
choice of \mbox{$i\in \Lambda_{P'}$}. Moreover, \eqref{eqn:a'i,b'i} and
\eqref{eqn:beta0} show that $\widetilde{U}_i$, $\widetilde{V}_i$,
$\widetilde{X}_i$, $\widetilde{Y}_i$ determine $a'_i$, $b'_i$,
$\alpha_{\varepsilon}$, $\beta_{\varepsilon}$. 

\medskip%\noindent 
{\it Step 2. Existence.} Using the above notations, we choose for each
\mbox{$i\in \Lambda_{P'}$} power series \mbox{$a_i,b_i\in K[[t_i]]$}
satisfying 
\mbox{$
\ord_{t_i} a_i, \,\ord_{t_i} b_i\geq \ord \varphi_i=\ord_{t_i} x_i
$}
and the compatibility condition
$$
\left(\frac{b_i-a_i\beta_0}{x_i}\right)(0) =
\left(\frac{b_j-a_j\beta_0}{x_j}\right)(0) \ 
\text{ for all }i,j\in 
\Lambda_{P'}\,. 
$$
Since \mbox{$v_i(0)=0$}, this allows us to define the needed extensions
according to \eqref{eqn:beta0} and \eqref{eqn:a'i,b'i}.
%\mbox{$\ord_{t_i} X_i, \ord_{t_i} Y_i\geq \ord \varphi_i$}.
\end{proof}

\begin{remark}\label{rmk:section unique}%[Uniqueness of Sections]
The uniqueness statement of Proposition \ref{prop:2.9} can be reformulated as
follows: Let \mbox{$\widetilde{A}\to A$} be  a small extension in $\sA_K$,  and
let  \mbox{$\widetilde{\xi}=(\varphi_{\widetilde{A}},
\osigma,\sigma)\in
  \Defem_{\Rbar\leftarrow P}(\widetilde{A})$} be a lifting of 
 \mbox{$\xi\in \Defem_{\Rbar\leftarrow P}(A)$}. Further, let
\mbox{$\pi'_{\widetilde{A}}:P_{\widetilde{A}}\to P'_{\widetilde{A}}$} be a
formal blow-up of the section $\sigma$. Then there is a unique morphism
\mbox{$\varphi'_{\widetilde{A}}$} as in the above diagram and at most one
section \mbox{$\sigma':P'_{\widetilde{A}}\to A$} such that
\mbox{$(\varphi'_{\widetilde{A}},\osigma',\sigma')\in
  \Defsec_{\Rbar\!\,'\leftarrow P'}(\widetilde{A})$}. 
Indeed, the diagram shows that $\sigma'$ exists iff the composition
\mbox{$\osigma'_i\circ \varphi'_{\widetilde{A}}$} is independent of the choice
of \mbox{$i\in \Lambda_{P'}$}.
\end{remark}

%\noindent
As a corollary of the proof of Proposition \ref{prop:2.9} (applied to
\mbox{$A=K$}, \mbox{$\widetilde{A}=\Keps$}), we obtain the following lemma
which allows us to argue by induction:

\begin{lemma}\label{lemma3.10}
\begin{Enumerate}
\item[(1)] \ Given \mbox{$a'_i, b'_i \in t_iK[[t_i]]$}, \mbox{$i\in
  \Lambda_{P'}$}, let  
\begin{equation}\label{eqn:lemma3.10}
a_i=a'_i\,,\quad b_i = a'_i\bigl(\varphi'_i(y)\!\!\:+\!\!\:
\beta_0\bigr)\!\!\:+\!\!\: 
  (b'_i\!\!\:+\!\!\: \beta_{\varepsilon})x_i\,,
\end{equation}
where \mbox{$\beta_0\!\!\:=\!\!\:(\beta\:\text{mod}\;\fm_A)
  ,\,\beta_{\varepsilon}\in K$}.  
Then the \mbox{$a_i, b_i$}, \mbox{$i\in
  \Lambda_{P'}$}, define an element of \mbox{$\Ies_{{\Rbar}\!\,^\prime
    \leftarrow P}$} iff  \mbox{$\,\min\{\ord_{t_i}a_i ,\, \ord_{t_i}b_i\} \geq 
  \ord_{t_i} x_i$} and the \mbox{$a'_i, b'_i$}, \mbox{$i\in \Lambda_{P'}$}
define an element  of \mbox{$\Ies_{{\Rbar}\!\,^\prime \leftarrow P'}$}.

\item[(2)] \ Let \mbox{$a_i, b_i \in t_iK[[t_i]]$}, \mbox{$i\in
  \Lambda_{P'}$} define an element of \mbox{$\Ies_{{\Rbar}\!\,^\prime
    \leftarrow P}$}. Then there exists a unique \mbox{$\beta_{\varepsilon}\in
  K$} such that setting
\begin{equation}\label{eqn:lemma3.10b}
a'_i=a_i\,,\quad  b'_i= \frac{b_i-a_i(\varphi'_i(y)+\beta_0)}{x_i}
-\beta_{\varepsilon}\,. 
\end{equation}
defines an element of \mbox{$\Ies_{{\Rbar}\!\,^\prime
    \leftarrow P}$}  satisfying \eqref{eqn:lemma3.10}.
\end{Enumerate}
\end{lemma}

\smallskip%\noindent
Note that \mbox{$\Rbar=\bigoplus_{P'}\Rbar'$}, hence \mbox{$\Ies_{{\Rbar} \leftarrow
    P}=\bigoplus_{P'}\Ies_{{\Rbar}\!\,' 
    \leftarrow P}$} where the sum on the right-hand side is taken over all
points $P'$ on $R$ in the first neighbourhood of $P$. 

  \begin{corollary}\label{cor:Ies is submodule}
    $\Ies_{{\Rbar} \leftarrow P}$ is an $R$-submodule 
    of  \mbox{$\fmbar \oplus \fmbar$}.
  \end{corollary}

\begin{proof}
We argue by induction on the number of blow-ups needed to resolve the
singularity $R$. If $R$ is regular, \mbox{$\Ies_{{\Rbar} \leftarrow
  P}=\fmbar\oplus\fmbar$}. If $R$ is singular, we may blow up the maximal ideal
in $P$ and consider the strict transform $R'_i$ at each infinitely near point
$P'$ in the first neighbourhood of $P$ on $R$. By the induction hypothesis, for
each such $P'$, 
\mbox{$\Ies_{{\Rbar}\!\,^\prime \leftarrow 
        P'}$} is an $R$-submodule of  \mbox{$\fmbar\!\,' \oplus
      \fmbar\!\,'$}, where \mbox{$\fmbar\!\,'=\bigoplus_{i\in \Lambda_{P'}}\!
        \fm_{\Rbar_i}$}. 
Then Lemma \ref{lemma3.10} shows that the same holds for
\mbox{$\Ies_{{\Rbar}\!\,^\prime \leftarrow P}$}. Indeed, the $K$-vector space
structure is obvious. Let \mbox{$h\in R$} and \mbox{$(\ba,\bb)\in
  \Ies_{{\Rbar}\!\,^\prime \leftarrow P}$}. Then, according to
\eqref{eqn:a'i,b'i} and \eqref{eqn:beta0}, we have
\mbox{$(ha_i)'=ha_i$} and  
\mbox{$(hb_i)'-hb'_i=
\bigl(h(\varphi'_i(x),\varphi'_i(y))-h(0))\cdot \beta_{\varepsilon}\in
\fm_{R'_i}$}. 
Since \mbox{$\fm_{R'}\!\oplus \fm_{R'}\subset \Ies_{{\Rbar}\!\,^\prime
    \leftarrow P}$}, \mbox{$\fm_{R'}=\bigoplus_{i\in
    \Lambda_{P'}}\!\fm_{R'_i}$}, and since 
\mbox{$h(b'_i)_{i\in \Lambda_{P'}}\!\in \Ies_{{\Rbar}\!\,^\prime \leftarrow
    P}$} by the induction hypothesis, we have \mbox{$\bigl((hb_i)')_{i\in
    \Lambda_{P'}}\!\in \Ies_{{\Rbar}\!\,^\prime \leftarrow P}$}. Then
Lemma \ref{lemma3.10} implies that \mbox{$h\bb\in
  \Ies_{{\Rbar}\leftarrow  P}$}, which proves the claim.
\end{proof}

\begin{lemma}\label{lem:2.11}
  Let the deformation  with trivial sections over $A$ defined by
  $\varphi_A$ be equisingular, and let \mbox{$(\ba,\bb)\in \fmbar\oplus
    \fmbar$}. Then the deformation with trivial
  sections \mbox{$(\varphi_{\widetilde{A}},\osigma,\sigma)$}  over
  $\widetilde{A}$ defined by   
\mbox{$\widetilde{X}_i:=X_i+\varepsilon a_i$},
\mbox{$\widetilde{Y}_i:=X_i+\varepsilon b_i$}, \mbox{$i=1,\dots,r$},
is equisingular iff \mbox{$(\ba,\bb)\in \Ies_{\Rbar\leftarrow P}$}.
\end{lemma}

\begin{proof}
%Let \mbox{$\varphi_{\widetilde{A}}:\widetilde{A}[[x,y]]\to \bigoplus_{i=1}^r
%  \widetilde{A}[[t_i]]$} be the defining map of the deformation, that is, 
%\mbox{$\varphi_{\widetilde{A}}(x)=\widetilde{X}_i$},
%\mbox{$\varphi_{\widetilde{A}}(y)=\widetilde{Y}_i$}, and
%\mbox{$\varphi_{\widetilde{A}}|_A=\id_A$}. 
We argue again by induction on the number of blow-ups needed to resolve the
singularity $R$. If $R$ is regular, the statement is obvious. Now, assume
that $R$ is singular. 

If the deformation \mbox{$(\varphi_{\widetilde{A}},\osigma,\sigma)$} is
equisingular, we find a formal blow-up
\mbox{$\pi'_{\widetilde{A}}$} of the 
trivial section $\sigma$ in \mbox{$P_A=A[[x,y]]$} such that, for each point
\mbox{$P'=K[[u,v]]$} in the first neighbourhood of $P$ on $R$, and for each  
\mbox{$i\in \Lambda_{P'}$}, the morphism 
\mbox{$\varphi_{\widetilde{A},i}$} factors as
\mbox{$\varphi_{\widetilde{A},i}=\varphi'_{\widetilde{A},i}\circ
  \pi'_{\widetilde{A}}$}. 
Moreover, we can assume that the \mbox{$\varphi'_{\widetilde{A},i}:A[[u,v]]\to
  K[[t_i]]$}, \mbox{$i\in \Lambda_{P'}$} define an equisingular deformation
with trivial sections of the parametrization \mbox{$P'\to \Rbar\!\,'$}.
Setting, as before, \mbox{$\varphi_{\widetilde{A},i}(u)=:U_i+\varepsilon
  a'_i$}, \mbox{$\varphi_{\widetilde{A},i}(v)=:V_i+\varepsilon b'_i$}, the
induction hypothesis gives that the $a'_i,b'_i$, \mbox{$i\in \Lambda_{P'}$}
define an element of \mbox{$\Ies_{\Rbar\!\,'\leftarrow P'}$}.
Since they necessarily satisfy the equality \eqref{eqn:lemma3.10}, Lemma
\ref{lemma3.10}\:(a) implies 
that \mbox{$(\ba,\bb)\in \Ies_{\Rbar\leftarrow P}$}. Note that the condition on
the order of $a_i,b_i$ is satisfied, since
\mbox{$\widetilde{X}_i,\widetilde{Y}_i$} 
defines an equisingular deformation along the trivial section.

Conversely, let us assume that \mbox{$(\ba,\bb)\in \Ies_{\Rbar\leftarrow P}$}.
Then \mbox{$(\varphi_{\widetilde{A}},\osigma,\sigma)$} is obviously
equimultiple. Moreover,
choosing \mbox{$a_i',b'_i\in K[[t_i]]$} according to
\eqref{eqn:lemma3.10b},
\mbox{$U_i+\varepsilon a'_i$}, \mbox{$V_i+\varepsilon b'_i$}
defines a deformation with trivial sections of the
parametrization \mbox{$P'\to \Rbar\!\,'$} which is obtained from
$\varphi_{\widetilde{A}}$ via a formal blow-up (see the proof of Proposition
\ref{prop:2.9}). Lemma \ref{lemma3.10}\:(b) together with the induction
hypothesis give that this deformation is, indeed, equisingular. 
\end{proof}

\begin{lemma}\label{lem:2.12}
\mbox{$(\varphi_A,\osigma,\sigma)\in 
  \Defsec_{\Rbar\leftarrow P}(A)$} is equisingular iff it is {\em formally
  equisingular}, that is, iff, for each 
\mbox{$N\geq 1$}, the induced deformation with sections over $A/\fm_A^N$
is equisingular.
\end{lemma}

\begin{proof}
Straightforward, by induction on the number of formal blow-ups.
%Since the necessity of formal equisingularity is obvious, let 
% \mbox{$(\varphi_A,\osigma,\sigma)$} be formally equisingular.
%As before, we may assume that \mbox{$\varphi_{A,i}:A[[x,y]]\to A[[t_i]]$}, 
%given by \mbox{$X_i,Y_i\in A[[t_i]]$} and
%that $\sigma, \osigma_i$ are trivial sections. Since, for all \mbox{$N\geq
%  1$}, the induced deformation $\varphi_{A/\fm_A^N}$ is equimultiple along
%$\sigma$, there can be no term in $X_i,Y_i$ of order less than the minimum of
%the order of the residues. That is, $\varphi_A$ is equimultiple along
%$\sigma$. Applying a formal blow-up as above, we get for each $P'$ on $R$ in
%the first neighbourhood of $P$ a deformation with trivial sections which is
%formally equisingular. Thus, proceeding by induction, we conclude that
%\mbox{$(\varphi_A,\osigma,\sigma)$} is equisingular. 
\end{proof}

\begin{proof}[Proof of Theorem \ref{theo3.8}]
Let \mbox{$A=K[[\bT]]$}, and let
    \mbox{$(\varphi_A,\osigma,\sigma)$} denote the deformation with trivial
    sections defined by \mbox{$X_i,Y_i\in t_iA[[t_i]]$}.

\medskip%\noindent
(1)  If \mbox{$(\varphi_A,\osigma,\sigma)$} is equisingular, then,
    for each \mbox{$j=1,\dots,k$}, the deformation 
    \mbox{$(\psi\varphi_A,\psi\osigma,\psi\sigma)$} induced by
    the projection \mbox{$\psi:K[[\bT]]\twoheadrightarrow K[[\bT]]/(\langle
      T_j^2\rangle+\langle T_\ell\mid \ell\neq j\rangle)$} is
    equisingular, too. Thus, Lemma \ref{lem:2.11} gives \mbox{$(\ba^{(j)},
  \bb^{(j)}) \in \Ies_{{\Rbar} \leftarrow P}$}.

 Conversely, let \mbox{$(\ba^{(j)}, \bb^{(j)}) \in \Ies_{{\Rbar}
    \leftarrow P}$} 
for all $j$. Since each extension of Artinian local rings factors through small
extensions,  Lemma \ref{lem:2.11} implies that 
\mbox{$(\varphi_A,\osigma,\sigma)$} is formally equisingular, hence
equisingular due to Lemma \ref{lem:2.12}.

\smallskip\medskip%\noindent
Since (3) is an immediate consequence of (2), it remains to prove (2):
 Let \mbox{$(\varphi_A,\osigma,\sigma)$} be a versal (respectively
  semiuniversal) object of $\Defes_{\Rbar\leftarrow P}$. Then, for any
  \mbox{$(\ba, \bb) \in \Ies_{{\Rbar} \leftarrow P}$}, the
    equisingular deformation with trivial sections defined by
    \mbox{$x_i+\varepsilon a_i$}, \mbox{$y_i +\varepsilon b_i$},
    \mbox{$i=1,\dots,r$}, can be induced (respectively uniquely induced) from
    \mbox{$(\varphi_A,\osigma,\sigma)$} via a morphism in $\sA_K$,
    \mbox{$K[[\bT]]\to \Keps$}, \mbox{$T_j\mapsto \beta_j\varepsilon$}. But
    this means 
    that \mbox{$x_i+\varepsilon a_i = X_i(\bbeta\varepsilon,t_i)$}, \mbox{$
y_i +\varepsilon b_i = Y_i(\bbeta\varepsilon,t_i)$} for some (respectively for
a unique) \mbox{$\bbeta=(\beta_1,\dots,\beta_k)$}. Expanding and comparing 
coefficients, we get 
$$(a_i,b_i)=\sum_{j=1}^k \beta_j
  \left(\frac{\partial X_i}{\partial T_j}(\bnull,t_i),\frac{\partial
      Y_i}{\partial T_j}(\bnull,t_i)\right)\,,$$ 
thus the necessity of the condition.   

\quad To show that the given deformation is versal (semiuniversal) along the
trivial sections, it suffices to show that it is formally versal
(semiuniversal), according to 
\cite[(5.2) Satz]{Fl}. Thus, it is sufficient to consider a small extension
\mbox{$\chi:\widetilde{C}\twoheadrightarrow C$} in $\sA_K$ with kernel
$\varepsilon K$, and equisingular deformations $\eta,\widetilde{\eta}$ over $C,
\widetilde{C}$  respectively. We assume that $\eta$ is induced from
$\widetilde{\eta}$ by $\chi$ and from $\xi$ by a morphism \mbox{$\psi:A\to C$}
in $\sA_K$,
and have to show that there exists a morphism \mbox{$\widetilde{\psi}:A\to
\widetilde{C}$} such that \mbox{$\chi\circ \widetilde{\psi}=\psi$} and that
\mbox{$\widetilde{\eta}$} can be induced from $\xi$ via $\widetilde{\psi}$. 
To show the semiuniversality, we have to show additionally that the tangent map
of $\widetilde{\psi}$ is uniquely determined.

We introduce the following notations:
\begin{itemize}
\item \ $\eta$ is given by $W_i(t_i),Z_i(t_i)\in C[[t_i]]$, 
\item \ $\widetilde{\eta}$ is given 
  by \mbox{$\widetilde{W}_i(t_i)=W_i(t_i)+\varepsilon
  w_i^{\varepsilon}$}, \mbox{$\widetilde{Z}_i(t_i)+\varepsilon
  z_i^{\varepsilon}$}, where \mbox{$w_i^{\varepsilon},z_i^{\varepsilon} \in
    K[[t_i]]$},   
\item \ $\psi(T_j)=\widetilde{\psi}(T_j)+\varepsilon \beta_j^{\varepsilon}$,
  where \mbox{$\beta_j^{\varepsilon}\in K$}.
\end{itemize}
Then the assumption is that there are a $C$-automorphism of $C[[x,y]]$, mapping
\mbox{$x\mapsto H_1(x,y)$}, \mbox{$y\mapsto H_2(x,y)$}, \mbox{$H_1,H_2\in
  \langle x,y\rangle C[[x,y]]$}, and
$C$-automorphisms of $C[[t_i]]$ mapping \mbox{$t_i\mapsto s_i\in t_iC[[t_i]]$},
  \mbox{$i=1,\dots,r$}, such that 
\begin{equation}\label{eqn:mod mC}
 H_1\equiv x \;\text{mod}\; \fm_C\,, \quad H_2\equiv y \;\text{mod}\; \fm_C ,
 \quad  s_i \equiv t_i \;\text{mod}\;   \fm_C 
\end{equation}
and \mbox{$
  X_i\bigl(\psi(\bT),t_i\bigr)=H_1\bigl(W_i(s_i),Z_i(s_i)\bigr)$},  \mbox{$
 Y_i\bigl(\psi(\bT),t_i\bigr)  =   H_2\bigl(W_i(s_i),Z_i(s_i)\bigr)$}.

We show that these automorphisms can be extended to 
$\widetilde{C}$-automorphisms
\begin{itemize}
\item \ \mbox{$x\mapsto H_1+\varepsilon h_1^\varepsilon$}, \ 
\mbox{$y\mapsto H_2+\varepsilon h_2^\varepsilon$}, \ 
\mbox{$h_1^\varepsilon,h_2^\varepsilon \in \langle x,y\rangle K[[x,y]]$},
\item \ \mbox{$t_i\mapsto s_i+\varepsilon s_i^\varepsilon$}, \  
\mbox{$s_i^\varepsilon \in \langle t_i\rangle K[[t_i]]$},
\,\mbox{$i=1,\dots,r$}, 
\end{itemize}
such that
\begin{eqnarray}
  X_i\bigl(\psi(\bT)+\varepsilon
  \bbeta^{\varepsilon},t_i\bigr)&=&(H_1+\varepsilon h_1^\varepsilon)
\Bigl(\widetilde{W}_i(s_i+\varepsilon
s_i^\varepsilon),\widetilde{Z}_i(s_i+\varepsilon 
s_i^\varepsilon)\Bigr), \label{eqn:rhs1}\\
 Y_i\bigl(\psi(\bT)+\varepsilon
  \bbeta^{\varepsilon},t_i\bigr) & = &  (H_2+\varepsilon h_2^\varepsilon)
\Bigl(\widetilde{W}_i(s_i+\varepsilon
s_i^\varepsilon),\widetilde{Z}_i(s_i+\varepsilon 
s_i^\varepsilon)\Bigr)\,. \label{eqn:rhs2}
\end{eqnarray}
Expanding the left-hand sides as power series in $\varepsilon$ (and using
\mbox{$\varepsilon \fm_{\widetilde{C}} =0$}), we get
%\begin{eqnarray*}
$$
  X_i\bigl(\psi(\bT)+\varepsilon
  \bbeta^{\varepsilon},t_i\bigr) = X_i\bigl(\psi(\bT),t_i\bigr)+ \varepsilon 
\sum_{j=1}^k \frac{\partial X_i}{\partial T_j}\bigl(\bnull,t_i\bigr) \cdot
\beta_j^{\varepsilon} 
$$
%\\
%  Y_i\bigl(\psi(\bT)+\varepsilon
%  \bbeta^{\varepsilon},t_i\bigr) &=& Y_i\bigl(\psi(\bT),t_i\bigr)+ \varepsilon
%\sum_{j=1}^k \frac{\partial Y_i}{\partial T_j}\bigl(\bnull,t_i\bigr) \cdot
%\beta_j^{\varepsilon} \,.
%\end{eqnarray*}
and similarly for $Y_i$.
%Next, we expand the right-hand sides of \eqref{eqn:rhs1},
%\eqref{eqn:rhs2}. First, 
Note that \mbox{$\widetilde{W}_i(s_i+\varepsilon s_i^\varepsilon) =  
 W_i(s_i) + \varepsilon \bigl(
 \dot{x}_i(t_i) \cdot s_i^\varepsilon + w_i^{\varepsilon}\bigr)$},
% $$
%\widetilde{W}_i(s_i+\varepsilon s_i^\varepsilon) = \widetilde{W}_i(s_i) +
%\varepsilon \frac{\partial 
%   \widetilde{W}_i}{\partial t_i}(s_i)\cdot s_i^\varepsilon =
% W_i(s_i) + 
% \varepsilon \bigl(
% \dot{x}_i(t_i) \cdot s_i^\varepsilon + w_i^{\varepsilon}\bigr)\,,
%$$
and similarly for $Z_i$. Expanding the right-hand side of \eqref{eqn:rhs1},
we get (taking into account \eqref{eqn:mod mC}, which implies that
\mbox{$\frac{\partial H_1}{\partial x}\equiv 1$\;mod$\;\fm_C$} and 
\mbox{$\frac{\partial H_1}{\partial y}\equiv 0$\;mod$\;\fm_C$})
\begin{eqnarray*}
 \lefteqn{ (H_1+\varepsilon h_1^\varepsilon)
\Bigl(W_i(s_i) + 
 \varepsilon \bigl(
 \dot{x}_i \cdot s_i^\varepsilon + w_i^{\varepsilon}\bigr),\:
Z_i(s_i) + 
 \varepsilon \bigl(
 \dot{y}_i(t_i) \cdot s_i^\varepsilon + z_i^{\varepsilon}\bigr)
\Bigr) }\hspace{1cm}\\
%& = & H_1\bigl(W_i(s_i),Z_i(s_i)\bigr) + \varepsilon \Bigl(
%h_1^\varepsilon\bigl(W_i(s_i),Z_i(s_i)\bigr) + 1\cdot  \bigl(
% \dot{x}_i \cdot s_i^\varepsilon + w_i^{\varepsilon}\bigr)
%\Bigr)\\
& = &  H_1\bigl(W_i(s_i),Z_i(s_i)\bigr) + \varepsilon
\Bigl(h_1^\varepsilon\bigl(x_i,y_i\bigr)+ \dot{x}_i \cdot
s_i^\varepsilon + w_i^{\varepsilon}\Bigr)\,.
\end{eqnarray*}
Comparing this with the expansion of the left-hand side of \eqref{eqn:rhs1}
shows that it suffices to find $\bbeta^\varepsilon, h_1^\varepsilon,
h_2^\varepsilon, \bs^\varepsilon$ satisfying the equality
$$
  \bigl(w_i^\varepsilon,z_i^\varepsilon\bigr)  =  \sum_{j=1}^k
\beta_j^{\varepsilon} \cdot
  \left(\frac{\partial   X_i}{\partial T_j}\bigl(\bnull,t_i\bigr),
    \frac{\partial   X_i}{\partial T_j}\bigl(\bnull,t_i\bigr)\right)
  -s_i^\varepsilon \cdot \bigl(\dot{x}_i, \dot{y}_i\bigr)
 -  \Bigl(h_1^\varepsilon\bigl(x_i,y_i\bigr),
  h_2^\varepsilon\bigl(x_i,y_i\bigr)\Bigr)  \,.
$$
From Lemma \ref{lem:2.11}, we know that
\mbox{$\bw^\varepsilon,\bz^\varepsilon\in \Ies_{\Rbar\leftarrow P}$}. Hence,
the assumptions imply that we find a solution \mbox{$\bbeta^\varepsilon,
  h_1^\varepsilon, h_2^\varepsilon, \bs^\varepsilon$} for the above equation 
(respectively a solution with uniquely determined \mbox{$\bbeta^\varepsilon$}).
\end{proof}

%\noindent
As a corollary of Theorem \ref{theo3.8} and Corollary \ref{cor:Ies is submodule}, 
we get the {\em ``openness
  of versality''} for equisingular deformations with sections of the
parametrization: let $(S,\ko_S)$ be an algebraic $K$-scheme, then we call a 
family \mbox{$\ko_S \to \ko_S[[x,y]] \to 
\bigoplus_{i=1}^r\ko_S[[t_i]]$} of parametrizations of reduced plane curve
singularities over $S$ equisingular, if for any (closed) point \mbox{$s \in 
  S$}, the induced family \mbox{$\widehat{\ko}_{S,s} \to \widehat{\ko}_{S,s}
  [[x,y]]  \to \bigoplus_{i=1}^r\widehat{\ko}_{S,s}[[t_i]]$} defines an
equisingular deformation of the parametrization with trivial sections. 
We say that the family is {\em  equisingular-versal at
  $s$}, if the induced 
family over the complete local ring $\widehat{\ko}_{S,s}$ is a versal
equisingular deformation of the special fibre.

\begin{corollary}\label{coro3.13}
  Let \mbox{$\ko_S \to \ko_S[[x,y]] \to \bigoplus_i\ko_S[[t_i]]$} be an
  equisingular family of parametrizations of reduced plane curve singularities.
  Then the set of (closed) points \mbox{$s \in S$} such that the family is
   equisingular-versal  at $s$
  is open in $S$.
\end{corollary}

\begin{proof}
  Let
  \mbox{$\kj^{\es}=\kj^{\es}_{\bigoplus_i\ko_S[[t_i]]\leftarrow\ko_S[[x,y]]}$}
  be the subsheaf 
  of \mbox{$\bigoplus_i \bigl(t_i\ko_S[[t_i]]\oplus t_i\ko_S[[t_i]]\bigr)$} of
  elements $(\ba,\bb)$ such 
  that \mbox{$X_i(t_i) + a_i(t_i)$},  \mbox{$Y_i(t_i) + b_i(t_i)$} defines 
an equisingular family over $S$. Here, as usual,
  \mbox{$X_i, Y_i \in \ko_S[[t_i]]$} denote the images of $x$ and $y$
  in $\ko_S[[t_i]]$.  Let $\dot{X}_i,\dot{Y}_i$ denote 
  the partials of $X_i,Y_i$ (with respect to $t_i$). Then
  the quotient sheaf  
\[
\kT^{1,\es}=\kT^{1,\es}_{\bigoplus_i\ko_S[[t_i]]\leftarrow \ko_S[[x,y]]} =
\kj^{\es}\!\!\left/ \left( \bigoplus_i
  \bigl(\dot{X}_i,\dot{Y}_i\bigr)t_i\ko_S [[t_i]]\right. 
 + \langle x,y\rangle \ko_S [[x,y]]\right),
\]
is a coherent $\ko_S$-sheaf.   Moreover, we have the Kodaira-Spencer
map \mbox{$ \Theta_S \to \kT^{1,\es}$}
which maps $\delta$ to the class of \mbox{$\bigl(\delta(X_i),
\delta(Y_i)\bigr)_{i=1}^r$}.
By Theorem \ref{theo3.8} the equisingular-versal locus is the complement of
the support of the cokernel of this map. Thus, the equisingular-versal locus
is open in $S$. 
\end{proof}

\section{Equisingular Deformations of the Equation}\label{sec:section 4}

%\noindent
We now turn to deformations of the singularity itself, that is, to deformations
of \mbox{$R=P/\langle f\rangle$} with or without section. Analogous to Definition
\ref{def:1.1}, they are defined by flat morphisms \mbox{$A\to R_A$} in $\sA_K$
(with section $\sigma$), together with a surjection \mbox{$R_A\twoheadrightarrow
  R$} over \mbox{$A\twoheadrightarrow
  K$} such that the corresponding diagram is Cartesian,
\[
\xymatrix@R=3pt@C=14pt{
R &  R_A\ar[l] \ar@/^1pc/[dd]^\sigma \\
\\
K \ar@{^{(}->>}@<2pt>[uu] &  A\ar@{->>}[l] \ar@{}[uul]|{\Box}\ar[uu] 
}\raisebox{-6ex}{.}
\]

Morphisms are morphisms of diagrams. We denote the corresponding
category by $\Defsec_R$, resp.\ $\Def_R$, and the deformation functors
by $\uDefsec_R$, resp.\ $\uDef_R$. In order to distinguish them from
deformations of the parametrization, we refer to the objects of
$\Def_R$ (resp.\ $\Defsec_R$)
also as {\em deformations of the equation\/} (resp.\ with section) since, basically, we deform $f$. From Proposition
\ref{prop:1.3}, we deduce the following 
commutative diagram
$$
  \xymatrix@R=3pt@C=14pt{
 & \Defsec_{\Rbar \leftarrow R \leftarrow P} \ar[dl]_(0.55)\cong
\ar@{->>}[dr]^(0.5){\quad \text{ smooth}} & \\
\Defsec_{\Rbar\leftarrow P}\ar@{-->}[dr] &  
& \Defsec_{\Rbar\leftarrow
  R}\ar[dl] \\
&   \Defsec_R\ar[dd] \\
\\
& \Def_R
}
$$
where the solid arrows are the natural forgetful functors and the
dashed arrow is defined by making the diagram commutative. For the 
deformation functors of isomorphism classes we have 

$$\uDefsec_{\Rbar\leftarrow P} \xleftarrow{\cong } \uDefsec_{\Rbar \leftarrow R
  \leftarrow P} \xrightarrow{\cong } \uDefsec_{\Rbar\leftarrow R} \to
\uDefsec_R \to \uDef_R\,.$$\\

We turn now to equisingular deformations.

\begin{definition}\label{def:3.1}
A deformation of $R$ is called {\em (strongly) equisingular\/} (or,
an {\em es-deformation\/}) if it
is induced by an equisingular deformation of the parametrization of
$R$. That is, we define the category $\Defes_R$ to be the full
subcategory of $\Def_R$,  
$$ \Defes_R = \text{image}\,\left(\Defes_{\Rbar\leftarrow P}\to
 \Def_R\right) \,.$$
\mbox{$\uDefes_R$} denotes the
corresponding subfunctor of $\uDef_R$. In particular, we introduce 
$$ \Ties_R = \uDefes_R (\Keps)\,. $$
\end{definition}

Similarly we define 
\[
\Defessec_R=\text{ image }(\Defes_{\Rbar\leftarrow P}\to \Defsec_R)
\]
as full subcategory of $\Defsec_R$ and call objects of $\Defessec_R$ {\em
  equisingular deformatons of $R$ with section}. The corresponding functor
  \mbox{$\uDefessec_R$} of isomorphism classes is called the {\em equisingular
  deformation functor with section}.
\medskip

Since, by Proposition
\ref{prop:sigma unique}, every equisingular deformation of
\mbox{$P\to \Rbar$} has a unique section $\sigma$, the forgetful
functors from the image of
$\Defes_{\Rbar\leftarrow P}$ in 
$\Defsec_R$ to $\Def_R$ is injective on objects, if $R$ is singular. Moreover, since the
section is singular, any isomorphism in $\Def_R$ must respect the
section and hence lifts to an isomorphism in $\Defsec_R$. That is,
\mbox{$\Defes_R \text{ and } \Def^{\es, \sec}_R$} are equivalent categories
and the functors \mbox{$\uDefes_R \text{ and }  \uDef^{\es, \sec}_R$} are
isomorphic. In particular, the vector spaces $\Ties_R$ and $\Tiesec_R$ are
isomorphic (but not equal: the first is a subspace of $R/J$, the second of
$\fm/\fm J$, cf. section \ref{sec:exact seq}).
\bigskip

Theorem \ref{theo3.8} yields immediately the first main result of this paper:

\begin{theorem}\label{theo:def of eqn}
(1)  The natural transformations \mbox{$\uDefes_{\Rbar\leftarrow P} \to
  \uDefessec_R \cong  \uDefes_R$} is smooth. In particular, equisingular
  deformations of $R$ are unobstructed and have
  a semiuniversal deformation with smooth base of dimension $\dim_K
  \Ties_R$. 

%\noindent
(2) \mbox{$\uDefes_{\Rbar\leftarrow
  P}\!\to \uDefes_R$} is an isomorphism of functors if and only if \mbox{$\dim_K
  \Ties_{\Rbar\leftarrow P} =\dim_K \Ties_R$}. If this holds, a semiuniversal
equisingular deformation of $R$ is obtained from the semiuniversal
equisingular deformation of the parametrization (as given in Theorem
\ref{theo3.8}\:(3)) by elimination of the parametrizing variables. 

(3) If the characteristic of $K$ is good, then
\mbox{$\uDefes_{\Rbar\leftarrow R} \cong \uDefes_R$}, and the semiuniversal object of
$\uDefes_R$ has an algebraic representative.
\end{theorem}

\begin{proof}
(1) Since $\uDefes_R$ is the image functor under \mbox{$\uDefes_{\Rbar\leftarrow
      P}\to \uDef_R$}, it follows that any versal object of
  $\uDefes_{\Rbar\leftarrow P}$ induces  a versal object
  of $\uDefes_R$. Since $\uDefes_{\Rbar\leftarrow P}$ is unobstructed by
  Theorem  \ref{theo3.8}, the same follows for $\uDefes_R$. Since
  $\uDefes_{\Rbar\leftarrow P}$ and $\uDef_R$ satisfy Schlessinger's
  conditions \cite{Schl} for the existence of a formal versal deformation (the
  first by Theorem \ref{theo3.8}, for the second this is well-known), the same
  holds for $\uDefes_R$.
Now, it follows
  from \cite[(5.2) Satz]{Fl} that $\uDefes_R$ has a semiuniversal deformation
  with smooth base, and that \mbox{$\uDefes_{\Rbar\leftarrow P}\to \uDefes_R$} is
  smooth. 

%\noindent
(2)   Since the functors
$\uDefes_R$ and $\uDefes_{\Rbar\leftarrow R}$ are both unobstructed,
they are equivalent iff the surjection \mbox{$\Ties_{\Rbar\leftarrow
    P} \to \Ties_R$} is an isomorphism, showing (2).

(3) In good characteristic, Proposition \ref{5.1} below says that
\mbox{$\Ties_R\cong \Ties_{\Rbar\leftarrow R}$}.
%\medskip

From Proposition \ref{prop:1.3} and its proof, we deduce that a semiuniversal
object of $\uDefes_R$ is obtained from the semiuniversal object for
$\uDefes_{\Rbar\leftarrow P}$ given in Theorem \ref{theo3.8}\:(3) by
eliminating the uniformizing parameters $t_i$ from the ideal generated
by  \mbox{$x-X_i(\bT,t_i)$} and \mbox{$y-Y_i(\bT,t_i)$},
\mbox{$i=1,\dots,r$}.  
The resulting power series
\mbox{$F_i(\bT,x,y)$}, \mbox{$i=1,\dots,r$}, respectively their product
\mbox{$F=F_1\cdot \ldots\cdot F_r$}, define the ideal of the total
space of the semiuniversal deformations of the branches $R_i$,
respectively of $R$.

Since $R$ has an isolated singularity, it is finitely determined,
hence we may assume that $x_i$ and $y_i$ are
polynomials in $t_i$, see \cite[Thm.\ B]{Hironaka}.  
Since $\Ties_{\Rbar\leftarrow P}$ is a finite dimensional $K$-vector
space, we can also choose the $a_i^{(j)}$ and the $b_i^{(j)}$ in
Theorem \ref{theo3.8}\:(3) to be polynomials in $t_i$.
Then \mbox{$X_i,Y_i\in K[\bT,t_i]$}, and we have to eliminate $t_i$
from \mbox{$x-X_i,\, y-Y_i$}, that is, we have to compute a generator
$F_i$ for the ideal \mbox{$\langle x-X_i,y-Y_i\rangle \cdot
  K[[x,y,\bT,t_i]]\cap K[[x,y,\bT]]$}. 

If \mbox{$p=\Char(K)$} is
good, we can compute $m_i$-th roots of units in $K[[t_i]]$
(\mbox{$m_i=\mult(R_i)$}), and hence we may assume that the
parametrization is of the form \mbox{$x_i(t_i)=t_i^{m_i}$}, and
\mbox{$\ord_{t_i} y_i(t_i)>m_i$}. But then \mbox{$\{x-X_i=y-Y_i=0\}$}
intersects the $t_i$-axis only in \mbox{$t_i=0$}. Thus, we can
eliminate in the polynomial ring, that is, we get \mbox{$\langle
  F_i\rangle=\langle x-X_i,y-Y_i\rangle \cap K[x,y,\bT]$} (cf.\
\cite{GrPf}). The product 
\mbox{$F=F_1\cdot \ldots\cdot F_r$} defines via \mbox{$K[\bT]\to
  K[x,y,\bT]/\langle F\rangle$} an algebraic representative of the
semiuniversal equisingular deformation of $R$.
\end{proof}

 % fix: to be completed
%\hspace{3cm}

\begin{remark}\label{rem:4.3}
\begin{Enumerate}
\item[(1)] \ Let \mbox{$\xi$} be an object of $\Defes_{\Rbar\leftarrow P}(A)$.
  Then Proposition \ref{prop:1.3} shows that there exists a unique lifting
to $\Defes_{\Rbar\leftarrow R\leftarrow P}(A)$, and this lifting
induces a deformation of the equation, \mbox{$\eta\in \Defes_{R}(A)$}. 
If $\xi$ is 
versal, then $\eta$ is versal too. If $\xi$ is semiuniversal then, however, 
$\eta$ need not be semiuniversal.

\quad  More precisely, if $\xi$ is semiuniversal and if \mbox{$\eta_{s}\in
   \Defes_R(B)$} is  
semiuniversal equisingular then,
by Theorem \ref{theo:def of eqn}, $\eta$ can be induced from $\eta_{s}$ by a map
\mbox{$B\to A$}, where 
\mbox{$A\cong B[[z_1,\dots,z_\ell]]$}. Here, \mbox{$\ell=\dim_K
  \ker(\Ties_{\Rbar\leftarrow P}\to 
  \Ties_R)$}, which can be computed by using the exact sequences in 
Section \ref{sec:exact seq}.

\item[(2)] \ For $K$ a field of characteristic $0$, J.\ Wahl \cite{Wa}
  introduced 
in a different way a functor $\ES$ of equisingular deformations of $R$ (over
Artinian rings). He considered an embedded resolution 
(by finitely many successive blow-ups) of the singularity $R$ . Then a deformation of $R$ is
equisingular in the sense of Wahl, if it is an em-deformation of the
equation along some section. Further he requires that, after blowing up the
section, there 
exist sections through the infinitely near points on $R$ in the first
neighbourhood of $P$ along which the blown up family induces an equisingular
deformation of the equation of the reduced total transform. Thus, the
definition is by induction on the number of blow-ups needed to resolve the
singularity. Any em-deformation of a node is equisingular.

\quad Notice that equimultiplicity for the parametrization (as in our
definition) differs from equimultiplicity for the equation (as in
Wahl's definition): an em-deformation of the parametrization induces an
em-deformation of the equation, but not conversely. For instance, it
can be easily seen that the em-deformation along the trivial section
of the cusp \mbox{$R=K[[x,y]]/\langle x^2\!-y^3\rangle $} given by
\mbox{$R_A=A[[x,y]]/\langle x^2\!-y^3-sy^2\rangle$} can be lifted to
$\Defsec_{\Rbar\leftarrow P}(A)$, but not to $\Defem_{\Rbar\leftarrow P}(A)$.
Indeed, the unique lifting is given by the parametrization
\mbox{$X(t)=t^3-s^2t$}, \mbox{$Y(t)=t^2-s^2$} which is not
equimultiple along any section.

\quad Hence, the relation between Wahl's $\overline{\ES}$, the functor
of isomorphism classes of $\ES$, and our $\uDefes_R$
is not completely obvious. Nevertheless, equisingular deformations in
Wahl's sense lift to equisingular deformations of the
parametrization and it can be shown that 
$\uDefes_R$ is isomorphic to $\overline{\ES}$, (cf. \cite{GrLoShu}).

\item[(3)] \ The main result of Wahl's paper \cite{Wa} is that $\overline{\ES}$
  is unobstructed. Theorem \ref{theo:def of eqn} provides a new proof. Wahl's
  proof appears to be more involved than ours, as it uses deformation theory of
  global objects, namely of divisors supported on the exceptional divisor of
  the (embedded) resolution of $R$. Moreover, as Wahl shows, there is
  no easy description of $\overline{\ES}$. For example, it is in
  general {\em not\/} a linear subfunctor of $\uDef_R$, the functor of
  deformations of $R$.  

\quad On the other hand, our functor $\uDefes_{\Rbar\leftarrow P}$ is {\em very
  easy to describe\/} as a 
linear subfunctor  of $\uDefsec_{\Rbar\leftarrow P}$. Hence, even in
characteristic $0$, our approach to equisingular deformations of $R$ via
deformations of the parametrization provides an easy understanding of
the objects in $\Defes_R$.

\quad By Theorem \ref{theo:def of eqn}\:(2), a semiuniversal base for
$\uDefes_R$ is obtained from a semiuniversal base for
$\uDefes_{\Rbar\leftarrow P}$ by elimination. As elimination is highly non-linear,
this ``explains'' why we cannot expect $\uDefes_R$ to be a linear
subfunctor of $\uDef_R$.  

%\item[(4)] \ Starting from the semiuniversal object of
%  $\uDefes_{\Rbar\leftarrow 
%    P}$ as given in Theorem \ref{theo3.8}, we get (in good characteristic) the
%  semiuniversal object of $\uDefes_R$ basically by eliminating the parameters
%  $t_i$ from the ideal generated by \mbox{$x-X_i(\bT,t_i)$},
%  \mbox{$y-Y_i(\bT,t_i)$}, \mbox{$i=1,\dots,r$}.  The resulting power series
%  \mbox{$F_i(\bT,x,y)$}, \mbox{$i=1,\dots,r$}, respectively their product
%  \mbox{$F=F_1\cdot \ldots\cdot F_r$}, define versal equisingular 
%  deformations of the branches $R_i$, respectively of $R$. This follows from
%  the proof of Proposition \ref{prop:1.3}. As elimination is highly non-linear,
%  this ``explains'' why we cannot expect $\uDefes_R$ to be a linear subfunctor
%  of $\uDefsec_R$. 

\quad On the other hand, Wahl introduced special equisingular deformations
$\text{ES}'\subset \text{ ES}$ such that $\overline{\text{ES}}'$ is a linear
subfunctor of $\uDef_R$. These special deformations lift to equisingular
deformations of the parametrization for which the elimination is linear.

\item[(4)] \ Let \mbox{$K=\C$}. All results proved so far are valid (with the
  same proofs) for convergent instead of formal power series. Then we use the
  geometric language of deformations of reduced plane curve singularities
  \mbox{$(C,0)=\bigl(\{f\!\!\:=\!\!\:0\}, 0\bigr)\subset (\C^2\!,0)$}
  over complex space germs $(T,0)$, where \mbox{$f\in \C\{x,y\}$}. Let
  \mbox{$\bigl(\overline{C},\overline{0}\bigr)\to (C,0)$} denote the
  normalization of $(C,0)$, and \mbox{$
    \bigl(\overline{C},\overline{0}\bigr)\to (C,0)\hookrightarrow (\C^2\!,0)$}
  the parametrization. We denote by \mbox{$(S^{\es}\!,0)$} the base space of
the semiuniversal equisingular deformation of $(C,0)$ (in this case isomorphic
to the base space of the semiuniversal equisingular deformation of the
  parametrization). Note that Theorem \ref{theo3.8} gives a
  convergent semiuniversal equisingular deformation (even an algebraic
  representative) and not only a formal object. \mbox{$(S^{\es}\!,0)$}
is a closed subgerm of the base space $(S,0)$ of the
  semiuniversal deformation \mbox{$(\sC,0)\to (S,0)$} of $(C,0)$. As is
  well-known (cf. \cite{GrLoShu}), \mbox{$(S^{\es}\!,0)$} coincides with the {\em $\mu$-constant
    stratum\/} of $(C,0)$, 
$$
(S^{\es}\!,0) = (S^\mu\!,0) = \Bigl(\bigl\{s\in S\:\big|\: \mu(\sC_s)=\mu(C,0)
\bigr\},\, 0\Bigr)\,, 
$$ 
where \mbox{$\mu(C,0)=\dim_{\C}\C\{x,y\}/\langle \frac{\partial f}{\partial x},\frac{\partial f}{\partial y}\rangle$} is the Milnor
number of $(C,0)$, \mbox{$\phi:\sC\to S$} a sufficiently small
representative of \mbox{$(\sC,0)\to (S,0)$}, and \mbox{$\mu(\sC_s)=\sum_{p\in
    \sC_s} \mu(\sC_s,p)$} is the total Milnor number of the fibre of $\phi$
over $s$. If \mbox{$\mu(\sC_s)=\mu(C,0)>0$}, then it is also known that there
is exactly one singular point \mbox{$p\in \sC_s$} (satisfying 
\mbox{$\mu(\sC_s,p)=\mu(\sC_s)$}), and the restriction of $\phi$ to
\mbox{$\phi^{-1} (S^\mu)$} admits a unique section, picking up the singular
point over \mbox{$s\in S^{\mu}$} (see \cite{Lazzeri,Teissier}).

\quad Hence, Theorem \ref{theo:def of eqn} gives a new proof that the $\mu$-constant
stratum \mbox{$(S^\mu\!,0)$} is smooth. 

\item[(5)] \ In good characteristic, the $\mu$-constant stratum can be
  generalized to the (strong) equisingularity stratum
  \mbox{$S^{\es}=\Spec(\Bes)\subset S=\Spec(B)$}. Here, $B$ is the
  base ring of the semiuniversal deformation \mbox{$B\to R_B$} of $R$,
  and \mbox{$\fp\subset B$} is a prime ideal defining a smooth
  subscheme such that the restriction \mbox{$\Bes:=B/\fp\to
  R_B\widehat{\otimes}_B B/\fp$} is equisingular. \mbox{$S^{\es}$} is the unique maximal closed subscheme
of $S$ having the following universal property: if  
a strongly equisingular deformation  \mbox{$A\to R_A$} of $R$ is induced 
from \mbox{$B\to R_B$} by some map \mbox{$\psi:B\to A$}, then 
\mbox{$\Spec(\psi): \Spec(A)\to S$} factors through
$S^{\es}$. 

\item[(6)] \ In bad characteristic, a ``strongly equisingular stratum''
 which generalizes the $\mu$-constant stratum does, in general,  {\em not exist}.
The situation is a follows: let \mbox{$\xi^{\es}\in \Defes_R(\Bes)$}
be a semiuniversal es-deformation, and let \mbox{$\eta\in \Defsec_R(B)$} be a
semiuniversal deformation with section. Then, by semiuniversality of
$\eta$, there exists a morphism \mbox{$\psi:B\to \Bes$} in $\sA_K$, which is
unique on the tangent level, such that \mbox{$\psi\eta\cong \xi^{\es}$}. 
By definition, $\uDefes_R$ is a subfunctor of $\uDef_R$, and we
have a unique subspace \mbox{$\Ties_R\subset  T^1_R$}. 
Since the tangent map of $\psi$ is
injective, the dual map \mbox{$\fm_B/\fm_B^2\to
  \fm_{\Bes}/\fm_{\Bes}^2$} is surjective. Thus, $\psi$ induces an isomorphism
\mbox{$B/\Ker\psi \cong \Bes$}, 
and the deformation over \mbox{$B/\ker \psi$} induced by $\eta$ is
isomorphic to $\xi^{\es}$, hence strongly equisingular semiuniversal.

\quad Thus \mbox{$\Spec(B/\ker\psi)\subset \Spec(B)$} could be considered as a
(strongly) 
equisingular stratum. However, it does not have an intrinsic meaning. Indeed, as Example \ref{exa:3.4} shows, $\psi$ is not unique and
\mbox{$\Spec(B/\ker\psi)\subset \Spec(B)$} may vary for different
choices of $\psi$ (only the tangent space 
is fixed). Moreover the restriction of
$\eta$ over the union for different choices is not strongly equisingular.

\quad However, we can define in arbitrary
characteristic an intrinsic {\em weak equisingularity stratum} in $\Spec(B)$
(see Section \ref{sec:6}) which coincides with the stratum $S^\es$ defined in
(5) if the characteristic is good.
\end{Enumerate}
\end{remark}

\begin{example}\label{exa:3.4}

%\noindent
%\hrulefill not yet fixed \hrulefill

%Let \mbox{$\Char(K)=2$} and consider the irreducible plane curve singularity
%given by \mbox{$f=y^4+x^6+x^7=(y^2\!+x^3)^2+x^7$}. One has 
%\[
%\langle f,x\frac{\partial f}{\partial x},x\frac{\partial f}{\partial y},y\frac{%\partial f}{\partial x},y\frac{\partial f}{\partial y}\rangle = \langle y^4\!+x%^6,\, x^7,\,
%x^6y\rangle\,. 
%\] 
%Thus, a basis for the $K$-vector space $\Tisec_R$ is given by the
%classes of the monomials with exponents in the set \mbox{$D=\{ (i,j)
%  \mid 0\leq i\leq 5, \, 0\leq j\leq 3, \, i+j\neq 0 \}$}, and the
%monomial $x^6$. In particular, one has \mbox{$\dim_K
%  \Tisec_{R}=24$}, and a semiuniversal deformation with
%section is given by
%\[
%F:= y^4+x^6+x^7 + \sum_{(i,j)\in D} u_{i,j} x^iy^j + u_{6,0}x^6\,.
%\]

We give now several examples which show that all possible pathologies in bad
characteristic do actually occur. Consider the plane curve singularities  (in bad characteristic) given by the equation $f$ in the
following cases:
\begin{enumerate}
\item [(1)] \mbox{$\Char (K)=p>0 \text{ and } f=y^{2p}+x^{2p+1}+x^py^{p+1}$}
\item [(2)] \mbox{$\Char (K)=2 \text{ and } f= y^4+x^6+x^7$}
\item [(3)] \mbox{$\Char (K)= p\geq 3, l=\frac{p+1}{2} \text{ and } f=y^p-x^{p+2}+x^ly^l$}
\item [(4)] \mbox{$\Char (K)= p>0 \text{ and } f=y(y^p-x^{p+1})$}
\end{enumerate}

The Tjurina ideal \mbox{$\langle f, x\frac{\partial f}{\partial x},
  x\frac{\partial f}{\partial y}, y\frac{\partial f}{\partial x},
  y\frac{\partial f}{\partial y}\rangle$} is given by the respective ideals
\begin{enumerate}
\item [(1)] \mbox{$\langle x^{2p+1}, x^{2p}y, x^{p+1}y^p, x^p y^{p+1},  y^{2p}\rangle $}

\item [(2)] \mbox{$\langle x^7, x^6y, y^4+x^6\rangle$}
\item [(3)] \mbox{$\langle x^{p+2}, 2x^{p+1} y-lx^{l-1}y^{l+1}, x^{l+1}y^{l-1},
    x^ly^l, y^p\rangle$}
\item [(4)] \mbox{$\langle x^{p+2}-xy^p, x^{p+1}y, x^py^2, y^{p+1} \rangle$}
\end{enumerate}

Thus, a basis for the $K$--vector space \mbox{$\Tisec_R$} is given by the
classes of the monomials with exponents in the set $D$ given by
\begin{enumerate}
\item [(1)] \mbox{$D=D_0\cup D_1\cup D_2$} where \\\mbox{$D_0=\{ (i, j)\ |\
  0<i+j\leq 2p, j<2p\},$}\\ \mbox{$ D_1=\{(i,j)\ |\ 1<i<p, p+2<j<2p,
  i+j>2p\},$} \\\mbox{$D_2=   \{(i,j)\ |\ p+2<i<2p, 1<j<p, i+j>2p\}$}
\item [(2)] \mbox{$D=\{(i,j)\ |\ 0\leq i\leq 5, 0\leq j\leq 3, i+j>0\}\cup\{(6,0)\}$}
\item [(3)] \mbox{$D=D_0\cup D_1\cup D_2\cup\{(p+1, 0)\}$} where\\
 \mbox{$D_0= \{ (i, j)\ |\
  0<i+j\leq p\} $}, \\ \mbox{$D_1=\{(i,j)\ |\ 1<i<l, l<j<p, i+j>p\},$}\\
\mbox{$ D_2=\{(i,j)\ |\  l+1<i<p+1, 1<j<l-1, i+j>p\}$}
\item [(4)] \mbox{$D=\{(i, j)\ |\ 0<i+j, i<p, j<p+1\}\cup \{(p,0), (p,1), (p+1,0)\}$}
\end{enumerate}

In particular, the dimension of \mbox{$\Tisec_R$} is given by $3p^2+1$ in case
  (1), 24 in case (2), \mbox{$2l^2+(l-1)^2$} in case (3), and
  \mbox{$p^2+p-2$} in case (4). The semiuniversal
  deformation with section $\eta$ of \mbox{$R=K[[x,y]]/\langle f\rangle$} is given by
\[
F:= f+\sum\limits_{(i,j)\in D} u_{i,j}x^iy^j
\]
over the base \mbox{$S=\Spec(B)\text{ where }B=K[[u_{i,j}\ |\ (i,j)\in D]]$}.

\textbf{Case (1)}: The subscheme \mbox{$S'\subset S$} given by the
smooth conditions
\[
\begin{array}{lcl}
u_{ij} & = & 0\ ,\ (i,j)\in D\smallsetminus\{(p,p), (2p,0)\}, \text{ and }
2u_{2p,0}-u_{p,p}^2=0, \text{ if } p\neq 2;\\
u_{ij} & = & 0\ ,\ (i,j)\in D\smallsetminus\{(4,0)\}, \text{ if } p=2
\end{array}
\]
can be viewed as an equisingular locus for $\eta$, and one has \mbox{$\dim
  S'=1+(p-1)(p-2)$}. Now, for every \mbox{$h\in\langle u_{ij}\ |\ (i,j)\in
  D\rangle$} one has the smooth \mbox{$(p-1)(p-2)$}--dimensional subscheme
  $S_h$ of $S'$ given by adding, to the equations defining $S'$, the following equation
\[
u_{p,p}-2h^2=0, \text{ if } p\neq 2;\  u_{4,0}-h^4=0, \text{ if } p=2.
\]
Then, for each $h$ the deformation induced by $\eta$ on $S_h$ is a semiuniversal
deformation for $\uDefes_R$, so that \mbox{$S_h=\Spec(B/\Ker\psi_h)$} for some
$\psi_h$ as in \ref{rem:4.3} (6). Notice that all the smooth subschemes $S_h$
share the same tangent space, namely the subspace $\Ties_R$ of
$\Tisec_R$. Also notice that $S'$ coincides with the Zariski closure of
the union of all $S_h$.

\textbf{Case (2)}: Take any pair \mbox{$h, h'\in \langle u_{ij}\ |\ (i,j)\in
  D\rangle$} and consider the unqiue \mbox{$h''\in \langle u_{ij}\ |\ (i,j)\in
  D\rangle$} which satisfies the equality
\[
(1+h''^2)u_{42}=h'^2+h''+(1+h''^2)u_{33}h^3+(1+h'')(u_{4,3}h^3+u_{5,2}h^2).
\]
$h''$ is the solution to an implicit equation in the power series
ring $B$. Consider the $5$--dimensional subscheme \mbox{$S_{h,
    h'}\subset S$} given by the smooth conditions
\[
\begin{array}{l}
u_{ij}=0, (i,j)\in D\smallsetminus\{(4,0), (4,2), (5,1), (6,0), (3,3), (4,3),
(5,2), (5,3)\}\\
u_{4,0}+h^4=0\\
u_{5,1}+h^2u_{3,3}=0\\
u_{6,0}+h''^2+u_{5,1}h +u_{4,2}h^2+u_{3,3}h^3=0
\end{array}
\]
%Again, the deformation induced by $\eta$ on $S_{h,h'}$ is a semiuniversal one
%for $\uDefes_R$, one has \mbox{$S_{h,h'}=\Spec(B/\Psi_{h,h'})$} for some
%$\Psi_{h,h'}$ as in \ref{rem:4.3} (6), and all $S_{h, h'}$ share as tangent
%space the subspace $\Ties_R$ of $\Tisec_R$. Finally, notice that the Zariski
%closure of \mbox{$\underset{h,h'}{\cup} S_{h,h'}$} is the 7--dimensional non
%smooth subscheme \mbox{$S\subset\Spec(B)$} given by
%\[
%\begin{array}{l}
%u_{ij}=0, (i,j)\in D\smallsetminus\{(4,0), (4,2), (5,1), (6,0), (3,3), (4,3),
%(5,2), (5,3)\}\\
%u_{4,0}+h^4=0\\
%u_{5,1}+h^2u_{3,3}=0\\
%u_{6,0}+h''^2+u_{5,1}h+u_{4,2}h^2+u_{3,3}h^3=0
%\end{array}
%\]
Again, the deformation induced by $\eta$ on \mbox{$S_{h, h'}$}   is
  semiuniversal for \mbox{$\uDefes_R$}, one has \mbox{$S_{h, h'}=\Spec(B/\Ker\psi_{h, h'})$} for some $\psi_{h,h'}$, as in \ref{rem:4.3} (6),
  and all \mbox{$S_{h, h'}$} share as tangent space the subspace \mbox{$\Ties_R$} of \mbox{$\Tisec_R$}. Finally, notice that the Zariski closure of
  \mbox{$\underset{h,h'}{\cup} S_{h, h'}$} is the $7$--dimensional non smooth
  subscheme \mbox{$S'\subset S$} given by 
\[
\begin{array}{l}
u_{ij}=0, (i,j)\in D\smallsetminus\{(4,0), (4,2), (5,1), (6,0), (3,3), (4,3),
(5,2), (5,3)\}\\[0.5ex]
u_{5,1}^2+u_{4,0} u_{3,3}^2=0.
\end{array}
\]

\textbf{Case (3)}: The subscheme \mbox{$S'\subset S$} given by the smooth
conditions
\[
u_{i,j}=0, \ (i,j)\in D, i+j\leq p \text{ or } i+j= p+1 \text{ and } j\leq
l-1,
\]
can be viewed as the equisingular locus for $\eta$, and one has
\mbox{$\dim S'=(l-2)^2$}. The deformaton induced by $\eta$ on $S'$ is a semiuniversal
deformation for $\uDefes_R$, therefore $S'$ is the only subscheme of type
\mbox{$\Spec(B/\Ker\psi)$} for $\psi$ as in \ref{rem:4.3} (6).

\textbf{Case (4)}: The subscheme \mbox{$S'\subset S$} given by the smooth
conditions
\[
u_{i,j}=0,\ (i,j)\in D, i+j\leq p +1\ ,
\]
again, can be viewed as the equisingular locus for $\eta$, and one has
\mbox{$\dim S'=\frac{1}{2}(p-1)(p-2)$}. The deformation induced by $\eta$ on $S'$ is
semiuniversal for $\uDefes_R$, and the subscheme $S$ is the only one of
type \mbox{$\Spec(B/\Ker\psi)$} with $\psi$ as in  \ref{rem:4.3} (6).

Conclusion: In all four cases  the subscheme $S'$ will become the weak equisingularity stratum
(section \ref{sec:6}). 

In case (1) (resp. (2)) $S'$ is the Zariski closure of infinitely many substrata
$S_h$ (resp. $S_{h, h'}$) which are smooth with the same tangent space. For each
$h$ (resp. $(h, h')$) the restriction of $\eta$ to $S_h$ (resp. $S_{h, h'}$)
is a seminuniversal strongly equisingular deformation of $R$. However,
the restriction of $\eta$ to the union of two different strata of $\{S_h\}$
(resp. $\{S_{h, h'}\}$) is not (strongly) equisingular. Hence, an intrinsic
  largest strong equisingularity stratum in $S'$ does not exist. In case (1)
  $S'$ is smooth while in case (2) $S'$ is singular.

In cases (3) and (4) $S'$ is a strong equisingularity stratum. It exists (and
coincides with the weak one) although the characteristic is bad.
\end{example}
\bigskip

%\newpage

\section{Exact Sequences of Infinitesimal Deformations}%%5
\label{sec:exact seq}

%\noindent
In this section, we consider maps between the tangent spaces of the
deformation functors introduced so far. We consider additionally 
deformations of the normalization (without sections) which fix
$R$. They are given by morphisms $R_A\to \Rbar_A$ over $A$ such that
\mbox{$R_A=R\,\widehat{\otimes}\, A$} (and each
morphism between them induces the identity on $R_A$). The
corresponding category, resp.\ deformation functor, is denoted by
$\Def_{\Rbar/R}$, resp.\  $\uDef_{\Rbar/R}$.

We denote by $T^0_{R}(\Rbar),\;
T^0_R,\; T^0_{{\Rbar}},\; T^0_{\Quot(R)}$ the $K$-vector space of derivations
of $R$ in ${\Rbar}$, $R$ in $R$, ${\Rbar}$ in ${\Rbar}$, and $\Quot(R)$ in
$\Quot(R)$, respectively.  Because all the considered derivations can be extended
to $\Quot(R)$, we can regard $T^0_{R}(\Rbar),\; T^0_R,\; T^0_{{\Rbar}}$ as 
vector subspaces of $T^0_{\Quot(R)}$.

%One has respective tangent spaces $T^1_{\Rbar/R},\;
%T^1_{{\Rbar} \leftarrow R},\; T^1_R$ to the functors $\Def_{{\Rbar}/R},\;
%\Def_{{\Rbar} \leftarrow R},\; \Def_R$, that is, the respective sets of
%deformation classes over the $K$--algebra $K[\varepsilon],\; \varepsilon^2=0$,
%endowed with a natural structure of $K$-vector spaces.   
From the obvious relations among the deformation functors, we deduce linear maps
\mbox{$T^1_{{\Rbar}/R} \to T^1_{{\Rbar} \leftarrow R} \to T^1_R$}.
The elements of \mbox{$T^1_{{\Rbar}/R}=\uDef_{{\Rbar}/R}(\Keps)$} can
be interpreted as derivation classes in the 
following sense: each deformation of \mbox{$R\to \Rbar$} which fixes $R$ is 
represented by a deformation \mbox{$R \otimes_K\!
K[\varepsilon] \hookrightarrow {\Rbar}\otimes_K\!
K[\varepsilon]$}, given by an injective morphism \mbox{$R+\varepsilon
R\to \Rbar+\varepsilon \Rbar$} of $K$-algebras mapping an element
$g$ of \mbox{$R+\varepsilon 
R$} to \mbox{$g+\varepsilon \partial g$} for some fixed \mbox{$\partial\in
T^0_{R}(\Rbar)$}.  
Since two such morphisms define isomorphic deformations
iff their derivatives are equal modulo
\mbox{$T^0_{{\Rbar}} \cap T^0_{R}(\Rbar)$}, we can identify
$T^1_{{\Rbar}/R }$ with the quotient
\mbox{$T^0_{R}(\Rbar)\big/\bigl(T^0_{{\Rbar}} \cap T^0_{R}(\Rbar)\bigr)$}. 
%(d.h.\, $T^0_{R,{\Rbar}}/T^0_{{\Rbar}} \cap T^0_{R,{\Rbar}}$).
The kernel of the map \mbox{$T^1_{{\Rbar}/R} \to T^1_{{\Rbar} \leftarrow R}$}
consists of the deformation classes determined by derivations in \mbox{$T^0_R
  \subset T^0_{R}(\Rbar)$}. Thus, it is identified with the $K$-vector space 
$$M_R :=T^0_R\big/\bigl(T^0_{{\Rbar}} \cap T^0_R\bigr)\,.$$
Note that, in characteristic zero, we
have \mbox{$M_R = 0$} as every derivation in $T^0_R$ can be extended to one
in $T^0_{{\Rbar}}$ (see \cite{Deligne}).   However, if \mbox{$\Char(K)=p > 0$},
this is not true. For instance, if \mbox{$p \nmid q$}, the derivative
$\tfrac{\partial}{\partial y}$ is tangent to the curve \mbox{$\{y^p +
  x^q=0\}$} and the induced 
derivation is in \mbox{$T^0_R$} but not in $T^0_{{\Rbar}}$\,.

The involved $K$-vector spaces can be described in terms of the parametrization
of the curve, too.  In fact, denote by 
%
%$\ord_t$ the $t$--order of the multi-power
%series in 
%${\Rbar}$ (that is, the $r$--tuple of the $t_i$--orders of the components), by
%$\dot{x}, \dot{y}$ the $t$--derivatives of \mbox{$x = x(t)$}, \mbox{$y =
%  y(t)$}, and by 
\mbox{$\bd =(d_1, \dots, d_r)$} the {\em differential
multi-exponent}, that is, \mbox{$d_i = \min\bigl\{\ord_{t_i}(\dot{x}_i),
  \ord_t(\dot{y}_i)\bigr\}$}.  Then, deformations of the parametrization
which give rise to classes in $T^1_{{\Rbar}/R}$ are precisely those
defined by power   series \mbox{$X_i,Y_i\in A[[t_i]]$} of type 
\mbox{$X_i = x_i
  + \varepsilon h_i t_i^{-d_i} \dot{x}_i$},  \mbox{$Y_i =
  y_i + \varepsilon h_i t_i^{-d_i} \dot{y}_i$}
  where \mbox{$h_i \in K[[t_i]]$}.  Deformations leading to elements 
  in $T^1_{{\Rbar} \leftarrow R}$ are precisely those given by 
\mbox{$X_i = x_i + \varepsilon a_i$}, \mbox{$Y_i(t_i) =
  y_i+ \varepsilon b_i$}, where \mbox{$a_i,b_i \in
    K[[t_i]]$}. 

Taking into account those deformations leading to trivial deformations for
the respective functors, we get the following lemma:

\begin{lemma}\label{lem:4.1} We have the following isomorphisms of $K$-vector spaces:
\setlength{\arraycolsep}{3pt}
\[
\begin{array}{rclrcl}
  M_R & \cong &\dfrac{\bt^{-\bd} {\Rbar}\cdot (\dot{\bx}, \dot{\by}) \cap (R \oplus
    R)}{{\Rbar}\cdot (\dot{\bx}, \dot {\by}) \cap (R \oplus R)}, &\qquad 
T^1_{{\Rbar}/R} & \cong & \dfrac{\bt^{-\bd}
  {\Rbar}\cdot (\dot{\bx},\dot{\by})}{{\Rbar}\cdot (\dot{\bx},\dot{\by})} \\[4.0ex]
T^1_{{\Rbar}\leftarrow R} & \cong & \dfrac{{\Rbar} \oplus
  {\Rbar}}{\Rbar\cdot (\dot{\bx},
  \dot{\by}) + (R \oplus R)}, &\qquad 
 T^1_R & \cong &R/J\,,
\end{array} 
\]
where $J$ is the Jacobian ideal of the curve, that is, the ideal of
\mbox{$R=K[[x,y]]/\langle f\rangle$} generated by the partials
$\frac{\partial f}{\partial x},\,\frac{\partial f}{\partial y}$.
In particular, the vector spaces are all equipped with a natural
$R$-module structure.
\end{lemma}

%\noindent
Here, $\bt^{-\bd} {\Rbar}$ is a short-hand notation for
$$\bigoplus_{i=1}^r t_i^{-d_i} K[[t_i]]\subset
  \Quot(R)=\bigoplus_{i=1}^r \Quot(\Rbar_i) =\bigoplus_{i=1}^r
  K((t_i))\,.$$ 
Note that, in good characteristic, we have \mbox{$d_i = m_i-1$}.   

Altogether, we get an exact sequence of $R$-modules, which are finite
dimensional $K$-vector spaces (see also \cite{Bu} and \cite{GrLoShu}).
$$
0 \lra M_R \lra T^1_{{\Rbar}/R} \lra T^1_{{\Rbar} \leftarrow R} \lra T^1_R \lra
R\big/{\Rbar}J  \lra 0\,.
$$
All maps are obvious, except for the map \mbox{$T^1_{{\Rbar}
    \leftarrow R} \!\to T^1_R$},  which takes the 
class of \mbox{$(\ba,\bb) \in {\Rbar} \oplus {\Rbar}$} to the class mod
$J$ of the element \mbox{$\ba\frac{\partial f}{\partial x} + 
\bb\frac{\partial f}{\partial y}$}.
Note that ${\Rbar}J$ is an ideal of $R$ as $J$ is contained in the
conductor $\kc$ of $R$. 

For deformations with sections we have an analogous exact sequence of $R$-modules:
% where $\fm$
%(respectively $\fmbar$) denotes the maximal ideal of $R$ (respectively 
%the radical of ${\Rbar}$). 
\begin{equation}\label{eqn:exact sequence of T1}
0 \lra \Msec_R \lra \Tisec_{{\Rbar}/R} \lra
\Tisec_{{\Rbar} \leftarrow R} \lra \Tisec_R \lra
\fm/\fmbar J \lra 0\,.
\end{equation}
Here, $\fm$ denotes the maximal ideal of $R$, $\fmbar$ the Jacobson radical
of $\Rbar$, and
$$\Msec_R :=
  \Tnsec_R\big/\bigl(\Tnsec_{{\Rbar}} 
\cap \Tnsec_R\bigr)\,,$$
where
\mbox{$\Tnsec_R =
  \{\partial\in T^0_R \mid \partial(\fm) \subset \fm\}$},
\mbox{$\Tnsec_{{\Rbar}} = \{\partial \in  
T^0_{{\Rbar}} \mid \partial(\fmbar) \subset \fmbar\}$}. 

\begin{lemma}\label{lem:4.2} We have the following isomorphisms of $R$-modules:
\setlength{\arraycolsep}{3pt}
\[
\begin{array}{rclrcl}
\Msec_R & \cong & \dfrac{\bt^{-\bd+\bone} {\Rbar}\cdot
  (\dot{\bx},\dot{\by}) \cap (\fm \oplus 
  \fm)}{{\fmbar}\cdot (\dot{\bx}, \dot{\by}) \cap (\fm \oplus \fm)}, &\qquad
\Tisec_{{\Rbar}/R} & \cong & \dfrac{\bt^{-\bd+\bone}
  {\Rbar}\cdot (\dot{\bx},\dot{\by})}{{\fmbar}\cdot (\dot{\bx},\dot{\by})}\\[4.0ex]
\Tisec_{{\Rbar}\leftarrow R} & \cong
  & \dfrac{(\fmbar \oplus \fmbar)}{\fmbar\cdot (\dot{\bx},\dot{\by}) + (\fm \oplus
      \fm)}, &\qquad \Tisec_R & \cong
    &\fm/\fm J\,. 
\end{array} 
\]
\end{lemma}

\begin{remark}
In the description of $M_R, \Msec_R$ as well as of
  $T^1_{{\Rbar}/R}, \Tisec_{{\Rbar}/R}$ in Lemma \ref{lem:4.1} and \ref{lem:4.2},
  each involved derivation 
  $\partial$ is represented by the tuple $(\partial x, \partial y)$ in
  \mbox{$R\oplus R$}, resp.\ in \mbox{$\Rbar\oplus\Rbar$}.
\end{remark}

%\noindent
From the above exact sequences, we deduce the following lemma:
\begin{lemma} Let \mbox{$\delta = \dim_K {\Rbar}/R$}. Then the
  following holds:
\begin{eqnarray*}
  \dim_K T^1_R  &=& \dim_K \Tisec_{R} -\dim_K(J/\fm J)+1\,, \\
  \dim_K T^1_R  &=& \dim_K T^1_{{\Rbar} \leftarrow R} + \delta + \dim_K M_R\,, \\
\dim_K \Tisec_R  &=& \dim_K \Tisec_{{\Rbar}\leftarrow R} + \delta + r-1
+ \dim_K \Msec_R\,,
\end{eqnarray*}
where \mbox{$1\leq \dim_K(J/\fm J) \leq 2$}. If
\mbox{$\Char(K)=0$}, then \mbox{$\dim_K M_R= \dim_K \Msec_R=0$} and,
if $R$ is not regular, \mbox{$\dim_K(J/\fm J) = 2$}. 
\end{lemma}
\begin{proof} The first equality follows from the
  definitions. Moreover, we have
$$\dim_K T^1_{{\Rbar}/R} = \dim_K
  \Tisec_{{\Rbar}/R} = |\bd|:= d_1+\ldots+d_r\,.$$
%Note that, in good
%characteristic, we have \mbox{$|\bd| = \mult(R)-r$}.
From the equality \mbox{$J{\Rbar} = t^{\bd}\kc$} and
the fact that \mbox{$\dim_K {\Rbar}/\kc = 2\delta$}, we get
$$\dim_K R\big/\bigl(J{\Rbar}\bigr) 
  = \delta + |\bd| \,, \qquad  \dim_K \fm/ (\fmbar \cdot \Rbar J) =
  \delta+r-1 + |\bd|\,.$$
%as one has \mbox{$\dim_K {\Rbar}/\fmbar = r$} and
%  \mbox{$\dim_K R/\fm = 1$}.   
%The above two formul{\ae} now follow from the
%  exact sequences. 
The fact that \mbox{$\dim_K(J/\fm J)=2$} in characteristic $0$ follows
e.g.\ from \cite{Casas}. %fix: (kappa,tau) im Vergleich
\end{proof}

%\noindent
%Now, we consider the tangent spaces \mbox{$\Ties_{{\Rbar}/R} \subset
%\Tisec_{{\Rbar}/R}$}, \mbox{$\Ties_{{\Rbar} \leftarrow R}
%\subset \Tisec_{{\Rbar} \leftarrow R}$} and \mbox{$ \Ties_R\!
%\subset \Tisec_R$} to the respective functors
%\mbox{$\uDefes_{{\Rbar}/R}$}, \mbox{$\uDefes_{{\Rbar} \leftarrow R}$} and
%$\uDefes_R$. 

\begin{proposition}\label{5.1} The exact
sequence \eqref{eqn:exact sequence of T1} induces an exact sequence of
$R$-modules
$$0 \lra \Msec_R \lra \Ties_{{\Rbar}/R} \lra \Ties_{{\Rbar}
  \leftarrow R} \lra \Ties_R \lra 0\,,
$$
where \mbox{$\Msec_R =\Ties_{{\Rbar}/R}=0$} if the characteristic of
$K$ is good.
\end{proposition}

\begin{proof}
To check the 
exactness, first note that the map \mbox{$\Ties_{{\Rbar}\leftarrow R}\! \to
\Ties_R$} is surjective by definition of $\uDefes_R$, and that elements in
its kernel are nothing but classes of equisingular deformations keeping 
$R$ fixed.  Second, the
  image of $\Msec_R$ in $\Tisec_{{\Rbar}/R}$ consists of
    the classes of trivial deformations of the parametrization keeping
    $R$ fixed. So, $\Msec_R$ is contained in $\Ties_{{\Rbar}/R}$.

In good characteristic, a deformation of \mbox{$\Rbar\leftarrow P$} given by
\mbox{$X_i = x_i + \varepsilon h_i 
  t_i^{-d_i+1} \dot{x}_i$}, \mbox{$Y_i = y_i + \varepsilon h_i
  t_i^{-d_i+1} \dot{y}_i$} is
equimultiple along the trivial section iff \mbox{$\ord_{t_i}(h_i) \ge
  d_i$}.  Thus, in this 
case, \mbox{$\Ties_{{\Rbar}/R}  = \{0\}=  \Msec_R$}, and
the exact sequence 
simply states that the map \mbox{$\Ties_{{\Rbar} \leftarrow R} \!\to
    \Ties_R$} is an isomorphism.  
\end{proof}

%\noindent
In bad characteristic, however, each term in the sequence of
Proposition \ref{5.1} can be
nonzero, as the following examples show:

\begin{example}\label{exa:5.6}
Let \mbox{$\Char(K)=p>0$}. 

\medskip%\noindent

(1) \: The irreducible plane curve singularity given by 
\[
f=y^{2p} +x^{2p+1}+x^py^{p+1}
\]
satisfies $\dim_K(M^\sec_R)=0, \dim_K(\Ties_{\Rbar/R})=1$. Indeed, it is
parametrized by
\[
x(t)=\frac{-t^{2p}}{1+t^{p+1}}\ ,\ y(t)=\frac{-t^{2p+1}}{1+t^{p+1}},
\]
and under the isomorphism of lemma \ref{lem:4.2}, $\Ties_{\Rbar/R}$ is
generated by the class of
\mbox{$(\dot{x},\dot{y})=\left(\frac{t^{3p}}{(1+t^{p+1})^2},
  -\frac{t^{2p}}{(1+t^{p+1})^2}\right)$}, however,
\mbox{$\Msec_R=(0)$} as $\dot{x}\not\in \fm$. Note that, in this case,
\mbox{$\Ties_{\Rbar\leftarrow R}$} and \mbox{$\Ties_R$} are not isomorphic; in
fact one has \mbox{$\dim_K(\Ties_{\Rbar\leftarrow R})=1$} and
\mbox{$\dim_K(\Ties_R)=0$}.
\medskip

(2) \: The irreducible plane curve singularity
given by
$$f=y^{p^2}-x^{p^2+p}-x^{p^2+p+1}=(y^p\!-x^{p+1})^p-x^{p^2+p+1}$$
satisfies \mbox{$\dim_K(\Msec_R)=\dim_K(\Ties_{\Rbar/R})=2$}. Indeed,
it is parametrized by $$x(t)=t^{p^2}\,,\quad
y(t)=t^{p^2+p}+t^{p^2+p+1}\,,$$
 and under the isomorphism of
Lemma \ref{lem:4.2}, $\Msec_R$ is the two-dimensional $K$-vector space
generated by the classes of $(0,t^{p^2})$ and $(0,t^{p^2+p})$, i.e. by the
classes modulo \mbox{$T^{0,\sec}_{\Rbar} \cap T^{0,\sec}_R$} of the elements of
\mbox{$T^{0,\sec}_R$} given by \mbox{$x\frac{\partial}{\partial y}$} and \mbox{$y\frac{\partial}{\partial y}$} . 
Note that, in this case, \mbox{$\Ties_{{\Rbar} \leftarrow R} \cong \Ties_R$}.
\medskip%\noindent

(3) \: The irreducible plane curve singularity 
\[
f=y^p-x^{p+2}+x^ly^l
\]
for $p\geq 3$ and \mbox{$l=\frac{1}{2}(p+1)$}  satisfies
  \mbox{$\dim_K\Msec_R=\dim_K\Ties_R=0$}.

Indeed, it is parametrized by
\[
x(t)=\frac{t^p}{(1+t)^{l-1}}\ ,\ y=\frac{t^{p+2}}{(1+t)^l}
\]
and, since $d=p$, under the isomorphism of lemma 5.2,
  \mbox{$\Ties_{\Rbar/R}=(0)$}, as one has
  \mbox{$(\dot{x},\dot{y})=\left(\frac{(1-l)t^p}{(1+t)^l},
  \frac{2t^{p+1}-lt^{p+2}}{(1+t)^{l+1}}\right)$}. It follows that the
  deformations of type
\[
\begin{array}{lcl}
X(t) & = & x(t)+\varepsilon s t^{-d+1}\dot{x}\\
Y(t) & = & y(t)+\varepsilon s t^{-d+1}\dot{y}
\end{array}
\]

\noindent
which are equisingular must satisfy $s\in t^d\Rbar$, so
  \mbox{$\Ties_{\Rbar/R}=(0)$}. Notice, that in this case, one has \mbox{$\Ties_{\Rbar\leftarrow R}\cong \Ties_R$}.
  %fix to be inserted

\medskip

(4) \: The plane curve singulartiy with $r\geq 2$ branches given by
\[
f=y(y-x^2)(y-x^3)\ldots (y-x^{r-1})(y^p-x^{p+1})
\]
satisfies \mbox{$\dim_K\Msec_R=\dim_K\Ties_R=0$}. Indeed, the differential
multi--exponent is given by the $r$--tuple \mbox{$d=(0, \ldots, 0, p)$}, so the elements in
\mbox{$t^{-d+1}\Rbar(\dot{x},\dot{y})$} which give rise to equisingular
deformations of the parametrization over $K[\varepsilon]$ need to be in
\mbox{$t\Rbar(\dot{x},\dot{y})$}. This shows $\Ties_{\Rbar/R}=(0)$. Also, in
this case, one has \mbox{$\Ties_{\Rbar\leftarrow R}\cong\Ties_R$}.
\end{example}

%\noinden
 Wahl showed in \cite{Wa} that the tangent space to his functor
 $\ES$ is an ideal \mbox{$I\subset P= K[[x,y]]$}, the important {\em
   equisingularity ideal\/} of $R$. Let us show how this generalizes
 in our context to arbitrary characteristic.

Let us denote by $\Def_{R/P}$ the category of deformations of
\mbox{$P\to R$} inducing the product deformation of $P$, also denoted
by (embedded) {\em deformations of $R/P$}. The forgetful functor
\mbox{$\Def_{R/P}\to \Def_R$} is smooth, and we denote by
$\Defes_{R/P}$ the preimage in $\Def_{R/P}$ of $\Defes_R$, and by
$\Defes_{R/P,\fix}$ the objects in $\Defes_{R/P}$ which are
equisingular along the trivial section.

\begin{proposition}
  The tangent space $\Defes_{R/P}(\Keps)$ of $\Defes_{R/P}$,
$$ \Ies := \bigl\{ g\in P \:\big|\: f+\varepsilon g \text{ defines an
  element of }
\Defes_{R/P}(\Keps) \bigr\}
$$
is an ideal containing the Tjurina ideal \mbox{$\langle 
  f,\frac{\partial f}{\partial x},\frac{\partial f}{\partial
    y}\rangle$}. Likewise, the tangent space 
$\Defes_{R/P,\fix}(\Keps)$ of $\Defes_{R/P,\fix}$,
$$ \Iesf := \bigl\{ g\in P \:\big|\: f+\varepsilon g \text{ defines an
  element of }
\Defes_{R/P,\fix}(\Keps) \bigr\}
$$
is an ideal containing \mbox{$\langle 
  f \rangle + \fm_P \big\langle\frac{\partial f}{\partial
    x},\frac{\partial f}{\partial y}\big\rangle$}. 

Moreover, the canonical map \mbox{$\Ies_{R/P,\fix}\to \Ies_{R/P}$}
induces an isomorphism 
$$ \Ties_R \cong \Ies \big/ \langle 
  f,\tfrac{\partial f}{\partial x},\tfrac{\partial f}{\partial
    y}\rangle \cong  \Iesf \big/ \bigl(\langle 
  f \rangle + \fm_P \big\langle\tfrac{\partial f}{\partial
    x},\tfrac{\partial f}{\partial y}\big\rangle\bigr)\,.$$
\end{proposition}

\begin{proof}
  Consider the image of the $R$-module $\Ies_{\Rbar\leftarrow P}$
  (Corollary \ref{cor:Ies is submodule}) under the map \mbox{$\Rbar\oplus \Rbar \to R$}
  given by 
\mbox{$(\ba,\bb) \mapsto \ba\frac{\partial f}{\partial
    x} +\bb\frac{\partial f}{\partial
    y}$}.
$\Iesf$ is the preimage in $P$ of this image, hence an ideal in $P$.

A slight modification of the proof of Corollary \ref{cor:Ies is submodule} shows that 
$$
\Ies_{{\Rbar} \leftarrow P} := \left\{
(\ba,\bb) \in \fmbar\oplus \fmbar \:\left|\:
  \begin{array}{l}
\bigl\{(x_i(t_i) + \varepsilon a_i, y_i(t_i) + \varepsilon b_i) \:\big|\:
i=1,\dots,r \bigr\} \\
\text{defines an es-deformation of $P\to \Rbar$ }\\
\text{over $\Keps$ along some sections $\osigma,\sigma$} 
  \end{array}\!
\right.\right\}
$$
is an $R$-submodule of \mbox{$\Rbar \oplus \Rbar$}. Taking the image
of this submodule in $R$ under the same map as above, and then the
preimage in $P$, gives $\Ies$. 

Note that \mbox{$\Iesf \big/ \bigl(\langle 
  f \rangle + \fm_P \big\langle\frac{\partial f}{\partial
    x},\frac{\partial f}{\partial y}\big\rangle\bigr)$} is the tangent space to
the image of $\uDefes_{\Rbar\leftarrow P}$ in $\uDefsec_R$, which is
isomorphic to $\Ties_R$ by the remark after Definition \ref{def:3.1}. 
\end{proof}

\medskip
\section{Weakly Equisingular and Weakly Trivial Deformations}%%6
\label{sec:6}

%\noindent
Our definition of (strongly) equisingular deformation of the equation is the canonical extension
of equisingularity 
to fields of arbitrary characteristic $p$. However, there is a new phenomenon in
characteristic \mbox{$p>0$}, which does not appear in characteristic $0$. For
example, the deformation over $K[[a]]$ given by
\mbox{$y^p\!+ax^p +x^{p+1}$} is not trivial, but it becomes trivial after the base change \mbox{$a\mapsto a^p$}. Similarly, 
\mbox{$y^{2p}\!+2ax^p y^p +a^2x^{2p} +x^{2p+1}+x^p y^{p+1}$} does not define an
equisingular deformation, but after the base change \mbox{$a\mapsto a^p$}
it does (see Example \ref{exa:5.6}(1)).
Since this is a very natural phenomenon in positive characteristic, we
introduce the new concepts of weak triviality and weak equisingularity.
In this section we just state the main results about weakly equisingular
deformations. The proofs (which rely on the results of section
\ref{sec:equipolygonal}) are deferred to section \ref{sec:8}.

%deformations which become trivial (resp.\ strongly equisingular) after
%some base change weakly trivial (weakly equisingular).   

We start by considering weakly equisingular deformations which are induced from
a fixed deformation of $R$.

For \mbox{$C\in \sA_K$}, let $\sA_C$ denote the category of Noetherian
complete local $C$-algebras. For the following definition, we fix
\mbox{$C\in \sA_K$} and an object \mbox{$\eta\in \Defsec_R(C)$}.%, that is, a
%deformation \mbox{$C\to R_C$} of $R$ with section $\tau$.   

\begin{definition}
  A {\em weakly equisingular deformation\/} ({\em wes-deformation\/}) of $R$
  {\em based in  \mbox{$\eta=(C\to R_C,\tau)$}} over \mbox{$A\in \sA_C$}
  is a commutative diagram with Cartesian squares
$$
\UseComputerModernTips
  \xymatrix@C=11pt@R=4pt{
{\Rbar}\ar@{}[ddrr]|{\Box}
%\ar@/_2pc/[dddd] 
&& \ar[ll]
{\Rbar}_{A}\ar@/^2pc/[dddd]^(0.3){{\osigma}}  \\ 
\\
R \ar[uu]\ar@{}[ddrr]|{\Box}  && \ar[ll]
R_{A}\ar@/^/[dd]^(0.45){\!\sigma}\ar@{}[ddrr]|{\Box} \ar[uu] &&
R_C\ar[ll]|\hole\ar@/^/[dd]^\tau\\ 
\\
K \ar[uu]&& \ar[ll] A \ar[uu] && C \ar[uu] \ar[ll]
}\qquad\qquad\raisebox{-5ex}{$(\ast)$}
$$
such that \mbox{$\xi=\bigl(R_A\!\to\!\!\:\Rbar_A,\osigma,\sigma\bigr)$} is an
object in \mbox{$\Defes_{\Rbar\leftarrow R}(A)$}. We refer to such a diagram by
writing $\xi/\eta$. If $\xi$ is trivial, we call
the deformation \mbox{$\xi/\eta$} {\em weakly trivial}. 
\end{definition}

%\noindent
Here, $\Defes_{\Rbar\leftarrow R}$ denotes the category  of (strongly) equisingular
  deformations of the normalization (Definition \ref{def:2.6}). We could have worked with
  \mbox{$\Defes_{\Rbar\leftarrow P}$} as well, but at this point we prefer
   deformations of the normalization (which have the induced deformation of
  $R$ explicitly as part of their data). Recall (Proposition \ref{prop:1.3})
  that the deformaton functors \mbox{$\uDefes_{\Rbar\leftarrow R}$} and
  \mbox{$\uDefes_{\Rbar\leftarrow P}$} are isomorphic.

A {\em morphism of wes-deformations\/} based in $\eta$ is given in an obvious way by a
commutative diagram (inducing the identity on $\eta$). The corresponding
category is denoted by $\Defwes_{R,\eta}$, while $\Defwes_{R,\eta}(A)$ denotes
the (non-full) subcategory of deformations over $A$ with morphisms being the
identity on $A$. 
If \mbox{$\psi:A\to B$} is a morphism in $\sA_C$, then the induced deformation
$\psi\xi$ is an object 
in $\Defes_{\Rbar\leftarrow R}(B)$. Hence, together with $\eta$, it defines an
object \mbox{$\psi(\xi/\eta)\in\Defwes_{R,\eta}(B)$}.

Similarly, we define the category \mbox{$\Def^{wtr}_{R, \eta}$} of weakly
trivial deformatons of $R$ based in $\eta$.

%Here we omit a formal definition of weak equisingular deformation of $R$
%(i.e. of the equation) without reference to a given deformation $\eta$ since
%we are mainly interested in the strong and weak equisingularity strata and the
%relation between them.

%\noindent
\medskip
We can now formulate the main results about weakly equisingular deformations of this paper:

\begin{theorem}\label{theo:6.3}
Let \mbox{$\eta=(C\to R_{C},\tau)\in \Defsec_R (C)$}  be a fixed deformation
with section of $R$. Then the following holds:

\begin{enumerate}
\item [(1)] There exist an algebra \mbox{$C_\eta\in \sA_C$} with structure
  morphism \mbox{$\Psi_\eta:C\to C_\eta$} and an object \mbox{$\zeta/\eta\in
  \Def^\wes_{R, \eta}(C_\eta)$} which has the following universal property:

If \mbox{$\xi/\eta\in \Def^\wes_{R,\eta}(A),\ A\in \sA_C$}, then there is a
unique morphism \mbox{$\Psi: C_\eta\to A$} in $\sA_C$ such  that
$\Psi(\zeta/\eta)$ is isomorphic to $\xi/\eta$. In particular, $C_\eta$ is
unique up to a unique isomorphism.
\item [(2)] The ideal \mbox{$I^\wes_\eta:=\Ker(\Psi_\eta)\subset C$} is
  uniquely determined by $\eta$ and satisfies:
\begin{enumerate}
\item [(i)] \mbox{$C^\wes_\eta=C/I^\wes_\eta \hookrightarrow C_\eta$} is finite
\item [(ii)] for each $z\in C_\eta$ there exists a \mbox{$q=p^l$} for some $l\geq 0$
  s.t. \mbox{$z^q\in C/I_\eta$}, where \mbox{$p=\Char(K)$} and $q=1$ if
  $p$ is good.
\end{enumerate}
\item [(3)] The construction of $\zeta/\eta$ is functorial in $\eta$ which
  means the following: Let
  \mbox{$\eta\in \Defsec_R(C)\ ,\ \eta'\in \Defsec_R(C')$} and
  \mbox{$\varphi:C\to C'$} a morphism s.t. \mbox{$\varphi(I^\wes_\eta)=I^\wes_{\eta'}$} and
  \mbox{$\varphi(\eta)\cong \eta'$}. Let \mbox{$\zeta/\eta\in
  \Def^\wes_{R, \eta}(C_\eta)$} resp. \mbox{$\zeta'/\eta'\in \Def^\wes_{R,
  \eta'}(C_{\eta'})$} be the universal objects as in (1). Then there exists a
  morphism  \mbox{$\Psi:C_\eta\to C_{\eta'}$} in $\sA_C$ (where
  \mbox{$C_{\eta'}\in \sA_C   \text{ via } C\xrightarrow{\varphi} C'\to
  C_{\eta'})$} such that   \mbox{$\Psi(\zeta/\eta)\cong
  \zeta'/\eta'$}. 
\item [(4)] Let $\Ess(R)$ denote the set of essential infinitely near points
  on $R$, $m_Q$ the multiplicity of the strict transform $R_Q$ of $R$ at $Q$
  and 
\[
\con^\wes_R:=\sum\limits_{Q\in \Ess(R)}\frac{m_Q(m_Q+1)}{2}- \ef_R
\]
where $\ef_R$ is the number of essential and free points on $R$. Then
\[
\dim C_\eta=\dim C_\eta^\wes\geq \dim C-\con^\wes_R.
\]

\item [(5)] If, moreover, $\eta$ is versal in $\Defsec_R$, then $C_\eta$ is a
  regular local ring, $I^\wes_\eta$ is a prime ideal and 
\[
\dim C^\wes_\eta = \dim C-\con^\wes_R.
\]
\end{enumerate}
\end{theorem}

We call $\zeta/\eta$ from (1) the {\em universal weakly equisingular
  deformation of $R$ based in $\eta$} and \mbox{$\con^\wes_R$} the {\em
  number of conditions defining weakly equisingular deformations of $R$}.

\begin{definition}\label{def:6.3}
The morphism $\Spec(C_\eta)\to \Spec(C)$ is finite and we call the image \mbox{$S^{\wes}_{\eta}:= \Spec(C^\wes_{\eta})\subset
  \Spec(C)$} with the scheme structure defined by the prime ideal $\Ies_{\eta}$ the {\em weak
equisingularity stratum of $\eta$\/}. In good characteristic strong and weak
  equisingularity coincide and we call \mbox{$S^\wes_R=S^\es_\eta$} the {\em
  equisingularity stratum}.

If \mbox{$\eta\in \Defsec_R(B^\sec_R)$} is the semiuniversal deformation of
$R$ we write \mbox{$S^{\wes,\sec}_R$} instead of \mbox{$S^\wes_\eta$}. It is a
  subscheme of the base space \mbox{$S^\sec_R:=\Spec (B^\sec_R)$} of $\eta$ and called the {\em
    weak equisingularity stratum of $R$\/}. 
\end{definition} 

We define now weakly equisingular deformations without reference to a given deformation.

\begin{definition} A deformation \mbox{$\eta\in \Defsec_R(C)$} is called a
  {\em weakly equisingular deformation (with section) of $R$} if there
  exists an injective ring map \mbox{$\varphi:C\hookrightarrow A$} such that the
  induced deformation \mbox{$\varphi\eta\in \Defsec_R(A)$} is strongly
  equisingular, that is, an object of \mbox{$\Defessec_R (A)$}. We denote the corresponding full subcategory of
  $\Defsec_R$ by $\Def^{\wes, \sec}_R$ and $\uDef^{\wes,\sec}_R$ the functor
  of {\em weakly equisingular deformations of $R$ with section}. If
  $\varphi\eta$ is trivial we call $\eta$ {\em weakly trivial}. 

As a corollary of Theorem \ref{theo:6.3} we show that \mbox{$\uDef^{\wes,
    \sec}_R$} has a semiuniversal deformation which can be identified with the
    restriction of the semiuniversal deformation for $\uDefsec_R$ to the weak
    equisingularity stratum \mbox{$S^{\wes, \sec}_R \subset S^\sec_R$}.
\end{definition}

\begin{theorem}\label{theo:6.5}
Let \mbox{$\eta\in \Defsec_R(C)$} be a versal (resp. semiuniversal) deformation
with section of $R$ and \mbox{$\pi: C\twoheadrightarrow C_\eta^\wes$} the
canonical surjection. Then \mbox{$\pi\eta\in \Defsec_R(C^\wes_\eta)$} is a
versal (resp. semiuniversal) weakly equisingular deformation of $R$ with section.
\end{theorem}

\begin{proof}
First let $\eta$ be arbitrary in \mbox{$\Defsec_R(C)$} and \mbox{$\Psi': C\to
  C'$} a morphism in $\sA_K$ such that \mbox{$\Psi'\eta\in
  \Def^{\wes,\sec}_R(C')$}. We claim that $\Psi'$ factors as
  \mbox{$C\overset{\pi}{\twoheadrightarrow}C_\eta^\wes\to C'$}.

By definition there exists a ring extension \mbox{$\varphi: C'\hookrightarrow
  A$} s.t. \mbox{$(\varphi\circ \Psi')\eta$} is induced by an equisingular
  deformation of the normalization \mbox{$\xi\in \Defes_{\Rbar\leftarrow R}(A)
  $}. By Theorem \ref{theo:6.3} (1) we have \mbox{$\xi/\eta\cong \Psi(\zeta/\eta)$}
  for a unique morphism \mbox{$\Psi:C_\eta\to A$} with $\zeta/\eta$ the
  universal object such that the following diagram commutes
\[
\UseComputerModernTips
  \xymatrix@C=11pt@R=4pt{
C_\eta\ar[dd]_\Psi & & \ C_\eta^\wes\ar@{_{(}->}[ll]\ar@{-->}[dd] & &C\ar@{->>}[ll]\ar[ddll]^{\Psi'}\\
\\
A & & \ C'\ar@{_{(}->}[ll]_\varphi 
}\raisebox{-7ex}{.}
\]
The dotted arrow exists since $\varphi$ is injective, which proves the claim.

Now let $\eta$ be versal and consider two objects $\Theta'$ resp. $\Theta''$
in \mbox{$\Def^{\wes, \sec}_R$} over $C'$ resp. $C''$ where
\mbox{$\chi:C''\twoheadrightarrow C'$} is a surjection and
\mbox{$\Theta'\cong\chi\Theta''$}. We have a commutative diagram
\[
\UseComputerModernTips
  \xymatrix@C=11pt@R=4pt{
& &  C\ar@/_1pc/_{\Psi'}[dddll] \ar@/^1pc/^{\Psi''}[dddrr]\ar@{->>}^\pi[dd]&&\\
\\
&&  C^\wes_\eta\ar@{-->}[dll]\ar@{-->}[drr]&\\
C' & & &&C''\ar@{->>}^\chi[llll]
}\quad\raisebox{-10ex}{,}
\]

\noindent
such that \mbox{$\Psi'\eta\cong\Theta'$}\ ,\ \mbox{$\Psi''\eta\cong\Theta''$}
by the versality property of $\eta$ where the dotted arrows exist by the first
part of the proof. This proves the versality of $\pi\eta$.

If $\eta$ is semiuniversal then the tangent map of $\Psi'$ is unique and hence
also the tangent map of \mbox{$C^\wes_\eta\to C'$}. Hence $\pi\eta$ is
semiuniversal too.
\end{proof}

\begin{remark} Let \mbox{$\eta\in\Defsec_R(B^\sec_R)$} be a semiuniversal
  deformation of $R$ and \mbox{$\zeta/\eta\in \Def^\wes_{R, \eta}(B_\eta)$} 
  the universal weakly equisingular deformation of $R$ based in $\eta$.

(1) From Example \ref{exa:3.4} the base space of the semiuniversal
weakly equisingular deformation \mbox{$S^{\wes, \sec}_R\subset S^\sec_R = 
  Spec(B^\sec_R)$}  is in general not smooth in bad characteristic. It is,
however, the image of the smooth space $\Spec(B_\eta)$ under the finite and
surjective map \mbox{$\Spec(\Psi_\eta): \Spec (B_\eta)\twoheadrightarrow
  S^{\wes, \sec}_R$} which is actually a homeomorphism by Theorem
\ref{theo:6.3} (2). In particular, $S^{\wes,
  \sec}_R$ is always irreducible of dimension 
\[
\dim S^{\wes, \sec}_R= \tau^\sec(R)-\con^\wes(R),
\]
and we can say that the conditions defining the weak equisingularity stratum  are independent.

In good characteristic and, as the examples in \ref{exa:3.4} suggest, also in many cases
of bad characteristic, the weak equisingularity stratum $S^{\wes, \sec}_R$ is
smooth. 

Here we use the well known fact  that if a map \mbox{$\varphi:A\to B$} in $\sA_K$ is injective then the map
\mbox{$\Spec(\varphi):\Spec (B)\to \Spec(A)$} is dominant (cf. \cite{GrPf},
Proposition A.3.8). If $\varphi$ is
finite then $\Spec(\varphi)$ is closed and hence $\Spec(\varphi)$ is
surjective.

(2) We comment the situation for (weakly) trivial
  deformations: If \mbox{$\Char(K)=0$} then it is known that
  if a deformation \mbox{$\eta\in \Def_R(C)$} becomes trivial after some base
  change \mbox{$\varphi: C\to C'$} then $\varphi$ factors as
  \mbox{$C\overset{\pi}{\twoheadrightarrow} C^{tr}\to C'$} where $C^{tr}$ is a
  unique factor algebra of $C$  such that $\pi\eta$ is
  trivial (cf. \cite{GrKa}, Lemma 1.4). Hence \mbox{$\Spec(C^{tr})\subset
  \Spec(C)$} is the unique maximal substratum over which $\eta$ is trivial (and
  a family is weakly trivial iff it is trivial). The proof uses Schlessinger's
  theory of functors of Artin rings. If $\Char(K)>0$, however, we do not know
  whether Schlessinger's conditions are satisfied and hence we do not know
  whether unique weakly trivial strata exist.
\end{remark}

Everything as above can be formulated for deformations of the
parametrization. To do so, fix any object \mbox{$\zeta=(C\to R_C\to
  \overline{R}_C,\tau, \overline{\tau})\in \Defsec_{\overline{R}\leftarrow R}(C)$} with
\mbox{$C\in \sA_K$}.

Let \mbox{$\Defes_{\overline{R}\leftarrow R,\zeta}$} denote the category whose
  objects are Cartesian diagrams of type
$$
\UseComputerModernTips
  \xymatrix@C=11pt@R=4pt{
{\Rbar}\ar@{}[ddrr]|{\Box}
%\ar@/_2pc/[dddd] 
&& \ar[ll]
{\Rbar}_{A}\ar@/^2pc/[dddd]^(0.41){{\osigma}} \ar@{}[ddrr]|{\Box}&& \ar[ll]
{\Rbar}_{C}\ar@/^2pc/[dddd]^(0.3){\overline{\tau}} \\ 
\\
R \ar[uu]\ar@{}[ddrr]|{\Box}  && \ar[ll]
R_{A}\ar@/^/[dd]^(0.45){\!\sigma}\ar@{}[ddrr]|{\Box} \ar[uu] &&
R_C\ar[ll]|\hole\ar@/^/[dd]^\tau\ar[uu]\\ 
\\
K \ar[uu]&& \ar[ll] A \ar[uu] && C \ar[uu] \ar[ll]
}
$$
where \mbox{$\xi=(A\to R_A\to \Rbar_A,\sigma, \osigma)$} is an object in
$\Defes_{\overline{R}\leftarrow R}(A)$, and whose morphisms are the obvious
morphisms of diagrams inducing the identity on $\zeta$. We write $\xi/\zeta$ for
such objects.

For \mbox{$A\in \sA_C$} denote by \mbox{$\Defes_{\overline{R}\leftarrow
    R,\zeta}(A)$} the (non-full) subcategory of deformations of the
parametrization over $A$ with morphisms being the identity on $A$. If
\mbox{$\psi:A\to B$} is a morphism in $\sA_C$ and if $\xi/\zeta$ is an object
in \mbox{$\Defes_{\overline{R}\leftarrow R,\zeta}(A)$}, then
\mbox{$\psi(\xi/\zeta)=(\psi\xi)/\zeta$} is the induced object in
\mbox{$\Defes_{\overline{R}\leftarrow R,\zeta}(B)$}.\\

The next theorem is the parametric analogue of Theorem \ref{theo:6.3}. Here
the situation is, however, simpler: 

\begin{theorem}\label{thm:6.3}
For each deformation \mbox{$\zeta=(C\to R_C\to \Rbar_C,\tau, \otau)$} in
\mbox{$\Defsec_{\Rbar \leftarrow R}$} there is a 
universal object \mbox{$\xi/\zeta$} for the category $\Defes_{\Rbar \leftarrow R,\zeta}$
defined over an algebra \mbox{$C_{\zeta}\in \sA_C$} (with structure morphism
\mbox{$\psi_{\zeta}:C\to C_{\zeta}$}) which is functorial in $\zeta$ and has
the following properties:

  \begin{enumerate}
  \item[(1)] 
    \mbox{$C_{\zeta}=C^\es_\zeta:=C/\Ies_{\zeta}$}, where \mbox{$\Ies_{\zeta}=\ker
        \psi_{\zeta}$}. 
\item[(2)] The codimension of \mbox{$\Spec(C^\es_\zeta)\subset\Spec(C)$}
  satiesfies 
\[
\dim C - \dim C_{\zeta}^\es \leq \con^\es_{\Rbar\leftarrow P}:= \sum_{Q\in \Ess(R)} m_Q - ef_R-(r\!\!\:-\!\!\:1),
\]
 where  $\Ess(R)$ is the set of essential
  infinitely near points of $R$, $m_Q$ is the multiplicity of the strict
  transform $R_Q$ of $R$ at $Q$, $ef_R$ is the number of free essential points of $R$ and $r$ the number of branches of $R$.
\item[(3)] If $\zeta$ is a versal deformation for \mbox{$\uDefsec_{\Rbar
      \leftarrow R}$}, then $C_{\zeta}^\es$ 
  is smooth and satisfies 
\[
\dim C - \dim C_{\zeta}^\es =
      \con^\es_{\Rbar\leftarrow P}.
\]
Moreover, the induced deformation
  $\zeta^{\es}$ of $\zeta$ on $C_{\zeta}^\es$ is a versal deformation for the
  functor \mbox{$\uDefes_{\Rbar \leftarrow R}$}, and $\zeta^{\es}$ is
  semiuniversal for  \mbox{$\uDefes_{\Rbar \leftarrow R}$} if $\zeta$ is
  semiuniversal for  \mbox{$\uDefsec_{\Rbar \leftarrow R}$}. 
  \end{enumerate}
\end{theorem}

%\noindent
Again, universal in the statement means that if $\xi'/\zeta$ is an object in
$\Defes_{\Rbar \leftarrow R,\zeta} (A)$, then there is a unique morphism
\mbox{$\psi:C_{\zeta}^\es \to A$} such that \mbox{$\psi(\xi/\zeta)\cong\xi'/\zeta$}. In
  particular, $C_\zeta^\es$ is unique up to isomorphism, and it is given by
\mbox{$\Ies_{\zeta}$} which is a uniquely defined ideal of $C$, depending
functorial on
$\zeta$. We call the subscheme \mbox{$S^{\es}_{\zeta}:=
  \Spec(C_{\zeta}^\es)=\Spec(C/\Ies_{\zeta})$} the
{\em equisingularity stratum of $\zeta$\/}. 

Notice that in the parametric case it is not necessary 
to consider the analogue of weak equisingular deformations since the universal
object for it is the subscheme \mbox{$S^{\es}_{\zeta}$} of $\Spec(C)$. That
is, if a deformation of the parametrization over $A$ becomes equisingular over $B$  after a
finite base change \mbox{$A\hookrightarrow B$}, then it was already
equisingular before over $A$. This corresponds to the fact that in Theorem
\ref{thm:6.3} (1) no base change \mbox{$C^\es_\zeta\hookrightarrow C_\zeta$} as in Theorem
\ref{theo:6.3} (2) (i) is required.

Sections \ref{sec:equipolygonal} and \ref{sec:8} will provide proofs for Theorems
\ref{theo:6.3} and \ref{thm:6.3}.

%It is the image of $\Spec(C_{\eta})$ in
%$\Spec(C)$ and, therefore, has the following universal property:
%If \mbox{$\eta'\in \Defsec_R(C')$} becomes equisingular
%after some dominant %fix
%base change, then any morphism \mbox{$\psi:C\to C'$} such that
%\mbox{$\psi\eta\cong \eta'$} satisfies that 
%\mbox{$\Spec(\psi): \Spec(C')\to \Spec (C)$} factors through $S^{\es}_{\eta}$.

\section{Equipolygonal Deformations}
\label{sec:equipolygonal}

In this section, we introduce equipolygonal deformations for embedded plane
curve singularities as well as for their parametrizations. Such deformations
are auxiliary tools for describing equisingular strata. We also show the
relationship between equipolygonal and equisingular deformations.

Throughout the following, we consider a fixed (embedded) plane curve
singularity \mbox{$P\to R=P/\langle f\rangle$}. \\

\begin{center}
{\em Equipolygonal deformations of the equation}
\end{center}
Let $Q$ be an infinitely near point of $P$ on $R$, and let \mbox{$Q\to
  R_Q=Q/\langle g\rangle$} be the (embedded) strict transform of $R$ at $Q$. If
\mbox{$Q\neq P$}, the reduced total transform of $R$ at $Q$ consists of either
one or two additional smooth exceptional branches, depending on whether $Q$
is a {\em free point\/} or a {\em satellite point}. Namely, the branch
\mbox{$Q\to   E_Q$} given by the exceptional divisor of the blowing up
  creating $Q$ and, in the satellite case, another exceptional branch
  \mbox{$Q\to D_Q$} which is the strict transform of $E_{P'}$ for some
  infinitely near point $P'$ of $P$ (such that $Q$ is infinitely near to
  $P'$). We denote by \mbox{$Q\to H_Q$} the curve singularity consisting of
  the exceptional branches at $Q$ (for \mbox{$Q=P$}, $H_Q$ is defined to be
  the zero ring). We set \mbox{$\kp_{P,R}:=\{P\}$}, and, for  \mbox{$Q\neq P$},
$$ \kp_{Q,R}:= \bigl\{\,\text{infinitely near points of $Q$ on $R$ which are on
  $H_Q$}\,\bigr\}. $$
Moreover, we define \mbox{$e=e_{Q,R}$} (resp.\
\mbox{$d=d_{Q,R}$}) to be the number of points in $\kp_{Q,R}$ which are on
$E_Q$ (resp.\ on $D_Q$).  If $E_Q$, resp.\ $D_Q$, does not exist, then
we set \mbox{$e=1$}, (resp.\ \mbox{$d=1$}). 

We say that two elements \mbox{$u,v\in \fm_Q\subset Q$} are {\em adapted to
  $H_Q$}  (or {\em adapted coordinates of $Q$\/}) if 
\mbox{$Q=K[[u,v]]$} and if all points in $\kp_{Q,R}$ are on \mbox{$Q/\langle
uv\rangle$}. That is, up to a permutation of $u,v$, we have intersection
multiplicities \mbox{$i(E_Q,Q/\langle
u\rangle)\geq e$} and \mbox{$i(D_Q,Q/\langle v\rangle)\geq d$} (if $E_Q,D_Q$
exist). In this case, we say that 
$u$ is {\em adapted to $E_Q$} and $v$ is {\em adapted to $D_Q$}. In particular,
if \mbox{$E_Q=Q/\langle u\rangle$} and \mbox{$D_Q=Q/\langle v\rangle$}, then
$u$ is adapted to $E_Q$ and $v$ is adapted to $D_Q$.
If $D_Q$ does not exist, we call each element \mbox{$v\in \fm_Q$} which is
transversal to $u$ adapted to $D_Q$. For \mbox{$Q=P$}, any two elements $u,v$
which generate $\fm_Q$ are adapted to $H_Q$.

We call \mbox{$u,v\in \fm_Q$} {\em generic adapted\/} elements
if the 
set of infinitely near points of $R$ on \mbox{$Q/\langle
uv\rangle$} coincides with $\kp_{P,R}$. In other words, two adapted elements
$u,v$ are generic adapted if the intersection multiplicity of $Q/\langle
uv\rangle$ with  $R_Q$ is minimal, that is, \mbox{$\ord_u g(u,0)=i(E_Q, R_Q)$}
and \mbox{$\ord_v g(0,v)=i(D_Q, R_Q)$.}

Up to a permutation of $u$ and $v$, the following objects depend only on $R$
and $Q$ (but not on the choice of adapted elements \mbox{$u,v\in Q$}): 
\begin{enumerate}
\item The Newton polygon \mbox{$N=N_{Q,R}$} of $g$ with respect to generic
  adapted elements.
\item The ideal \mbox{$\Iep=\Iep_{R_Q}\subset Q$} generated by the monomials
  with exponents in \mbox{$N+\Z_{\geq 0}^2$}.
\item The {\em adapted Jacobian ideal\/} \mbox{$J_{Q,R}\subset Q$} generated by
  $g,u\frac{\partial  
    g}{\partial u},u^e\frac{\partial
    g}{\partial v}, v^d\frac{\partial g}{\partial u},v\frac{\partial
    g}{\partial v}$ (assuming that \mbox{$E_Q=Q/\langle u\rangle$} and 
 \mbox{$D_Q=Q/\langle v\rangle$}). We have \mbox{$J_{Q, R}\subset I^{\ep}$}.
\item The finite dimensional $K$-vector space 
\mbox{$\Tep=\Tep_{Q,R}=\Iep/J_{Q,R}$}.
\end{enumerate}

%\noindent
It is easy to see that $J_{Q,R}$ is the tangent space to the group of adapted
automorphisms of $Q$. Here, an automorphism $\psi$ of  
\mbox{$Q$} is called {\em adapted}, if \mbox{$\psi(u)\in \langle
  u,v^e\rangle$}, \mbox{$\psi(v)\in \langle u^d,v\rangle$} (for any adapted
$u,v$ with $u$ adapted to $E_Q$ and $v$ adapted to $D_Q$). Note that adapted
automorphisms map adapted elements to adapted 
elements. 

\begin{definition}\label{def:7.1}
  Let \mbox{$\eta=(Q_A\!\to\!\!\:\Hbar_{Q,A},\osigma,\sigma)\in
    \Defsec_{\Hbar_Q\leftarrow Q}(A)$} be a deformation with section of the
    parametrization of $H_Q$, and let \mbox{$\xi=(A\to R_{Q,A}, \sigma')$}
    define an object of \mbox{$\Defsec_{R_Q}(A)$}. Then we say that $\xi$ is 
{\em adapted to $\eta$\/} if $\xi$ fits into a commutative diagram 

$$
\UseComputerModernTips
\xymatrix@C=12pt@R=4pt{
 && \Hbar_{Q,A}\ar@/^3pc/[dddd]^-{\osigma}\\
\\
R_{Q,A} \ar@/_1pc/[ddrr]_(0.4){\sigma'} & &Q_A \ar[uu] \ar@{->>}[ll]\ar@/^1pc/[dd]^\sigma\\
\\
 & &A\ar[uu]\ar[uull]
}
$$

%$$
%\UseComputerModernTips
%\xymatrix@C=12pt@R=4pt{
%R_{Q,A} \ar@/_/[ddddr] && 
%{\Hbar}_{Q,A}\ar@/^/[ddddl]^{\osigma}\\   
%& Q_A \ar[ur] \ar@{->>}[ul]\ar@/^/[ddd]^(0.45){\!\sigma}\\
% \\
%\\
%& \,A\,, \ar[uuu]
%}
%$$
\noindent
that is, $\eta$ and $\xi$ are deformations with basically the same section, namely
\mbox{$\sigma=Q_A\twoheadrightarrow R_{Q, A}\xrightarrow{\sigma'} A$}.
\end{definition}

%\noindent
Let \mbox{$u,v\in Q$} be generic adapted elements, and let \mbox{$U,V\in
  I_\sigma\subset Q_A$} be such that \mbox{$u\equiv U\!\!\mod \fm_A$} and 
\mbox{$v\equiv V\!\!\mod \fm_A$}. Then Nakayama's Lemma implies that $U,V$
generate $I_\sigma/I_\sigma^2$, thus, \mbox{$Q_A=A[[U,V]]$} and
\mbox{$I_\sigma=\langle U,V\rangle$}. 

We call  \mbox{$U,V$} {\em adapted\/} to
$\eta$, if \mbox{$\ord_U E(U,0)=e$} and \mbox{$\ord_V D(0,V)=d$}. 
Here, \mbox{$E(U,V)\in A[[U,V]]$} (resp.\ \mbox{$D(U,V)\in A[[U,V]]$}) is
the equation for the deformation of $E_Q$, (resp.\ $D_Q$) induced by $\eta$.

%
%Given generic adapted elements \mbox{$u,v\in Q$}, there are adapted elements
%$U,V\in I_\sigma\subset Q_A$ such that \mbox{$u\equiv U\!\!\mod \fm_A$} and
%\mbox{$v\equiv V\!\!\mod \fm_A$}. By Nakayama's Lemma, $U,V$ generate
%$I_\sigma$, that is, \mbox{$Q_A=A[[U,V]]$} and \mbox{$I_\sigma=\langle
%  U,V\rangle$}. For such $U,V$ we have an equation \mbox{$G(U,V)\in
%  A[[U,V]]$} for the adapted deformation $\xi$ of \mbox{$Q/\langle g\rangle$}.%

  \begin{definition}\label{def:6.5alt} Let
    \mbox{$\eta=(\phi_A,\osigma,\sigma)\in 
    \Defsec_{\Hbar_Q\leftarrow Q}(A)$}, and let 
 \mbox{$\xi\in \Defsec_{R_Q}(A)$} and 
{$U,V\in I_\sigma$} be adapted
to $\eta$. Moreover, let \mbox{$G(U,V)\in A[[U,V]]$} be
an equation for $\xi$.  Then $\xi$ is called 
\begin{enumerate}
\item {\em equiadapted\/} if \mbox{$\ord_U G(U,0)=\ord_u g(u,0)$},
  \mbox{$\ord_V G(0,V)=\ord_v g(0,v)$}. 
\item {\em equipolygonal\/} if \mbox{$G(U,V)\in \Iep_{R_Q} A[[U,V]]$}.  
% fix: sollte hier wohl eher heissen: das von U^iV^j, (i,j) oberhalb N,
% erzeugte Ideal
\end{enumerate}
\end{definition}

%\noindent
Note that this definition is independent of the choice of the (generic) adapted
elements $u,v,U,V$. 

We write $\Defep_{R_Q}$ (resp.\ $\uDefep_{R_Q}$) for the category (resp.\ the
functor of adapted isomorphism classes) of equipolygonal deformations of
$R_Q$. The vector space $\Tep$ can then be identified with the tangent space to
$\uDefep_{R_Q}$.  

For \mbox{$Q=P$}, equiadapted and equipolygonal deformations of $R_Q$ are
nothing but equimultiple deformations of the equation (along the section
prescribed by $\eta$).  If $Q$ is arbitrary, 
{\em equiadapted\/} deformations preserve the points of intersection
of the Newton polygon $N_{Q,R}$ with the $u$- and $v$-axis, while equipolygonal
deformations preserve the Newton polygon (for generic adapted coordinates).\\

%fix: das Folgende glaube ich nicht (hier muss man erst mal eta geeignet
%fixieren -- ich denke, so dass u,v adapted sind (also keine Not U,V
%einzufuehren).
 
For each deformation \mbox{$\xi\in \Defep_{R_Q}(A)$}, there is a well-defined
{\em Kodaira-Spencer map\/} \mbox{$\Psi=\Psi_{\xi}:T_A\to \Tep$}, where
\mbox{$T_A=\Der_K(A,K)=\Hom(A,\Keps)$} is the Zariski--tangent space to
$A$. For \mbox{$G\in \Iep A[[U,V]]$} inducing $\xi$, the map $\Psi$ takes a
(local) homomorphism 
\mbox{$\delta:A\to \Keps$} to \mbox{$[h]\in \Tep=\Iep/J_{Q,R}$}, where 
\mbox{$g+\varepsilon h\in \Keps[[u,v]]$} defines \mbox{$\delta\xi\in 
  \Defep_{R_Q}(\Keps)$}.  

A deformation \mbox{$\xi\in \Defep_{R_Q}(A)$} is called {\em
  equipolygonal versal\/} (or, {\em ep-versal\/}) if the
corresponding Kodaira-Spencer map $\Psi_{\xi}$ is surjective.
If $\Psi_{\xi}$ is an isomorphism, we call $\xi$ {\em
  equipolygonal semiuniversal} (or, {\em ep--semiuniversal}).

\begin{proposition}\label{prop:7.3}
Let \mbox{$h_1,\dots,h_s\in \Iep_{R_Q}$} represent a basis (resp.\ a system
  of gene\-rators) of $\Tep_{Q,R}$.  Then, with \mbox{$T_1, \ldots, T_s$} new
  variables and \mbox{$\bT=(T_1, \ldots, T_s)$},

$$ G=g + \sum_{i=1}^s T_i h_i \in K[[\bT]][[u,v]]$$

defines an ep--semiuniversal (resp.\ ep-versal) deformation of
$Q/\langle g\rangle$.

% fix: hier muessen wir Diagramme betrachten, denn \eta muss hier trivialen
% Schnitt erlauben....
\end{proposition}

%\noindent
In particular, each ep-versal deformation of $R_Q$ has a smooth base space.

\begin{proof}
Since \mbox{$\{\frac{\partial}{\partial T_i}\ |\ i=1, \ldots, r\}$} is a basis of
\mbox{$T_{K[[\bT]]}$}  and since the Kodaira--Spencer map for $G$ maps
  \mbox{$\frac{\partial}{\partial T_i}$} to $h_i$, the statement is almost a tautology.
 % fix: to be inserted
\end{proof}

\begin{remark}\label{rem:7.4} We introduced the notion of equipolygonal
  (semiuni--)versal at this stage only in order to have a convenient
  notation. The fact that this notion is equivalent to the ususal definition
  of (semiuni--)versality for the functor \mbox{$\uDefep_{R_Q}$} is by no
  means trivial and follows from  the results of section \ref{sec:8}.
\end{remark}

We generalize the notions from above to multicurves (resp. curve diagrams).

\begin{definition}\label{def:curve diagram}
A {\em curve diagram\/} is a finite list $\sC$ of infinitely near points of $P$
on $R$ (repetitions are allowed) together with, for each \mbox{$Q\in \sC$},
\begin{enumerate}
\item [(1)]the set of exceptional branches $E_Q,D_Q$ (if these exist), and
\item [(2)]a non-exceptional curve \mbox{$Q\to Q/\langle g\rangle$} (that is, a
  curve without exceptional branches) 
\end{enumerate}
such that the following holds: if \mbox{$Q\to Q'$} is a formal blow-up among
\mbox{$Q,Q'\in \sC$}, then the non-exceptional curve  \mbox{$Q'\to Q'/\langle
  g'\rangle$} at $Q'$ is the strict
transform of $Q/\langle g\rangle$ under the formal blow-up. We denote such a
curve diagram by $(\sC,\sG)$, where \mbox{$\sG=(Q/\langle g\rangle)_{Q\in
    \sC}$}. 

If none of the points in $\sC$ are consecutive, we also refer to  a {\em curve
  diagram\/} as a {\em multicurve}.  
\end{definition}

\begin{example}
Let \mbox{$g\in Q$} decompose into $s$ {\em tangential components}, that is,
\mbox{$g=g_1\cdot \ldots\cdot g_s$}, where the $g_i$ are unitangential and have
pairwise different tangent directions. Then the {\em multicurve of tangential
components\/} of $Q/\langle g\rangle$ is given by the list
\mbox{$\sC=(Q,\dots,Q)$}, $Q$ repeated $s$ times, together with, for the $j$th
entry, the set of exceptional branches $E_Q,D_Q$ (if these exist), and
the curve \mbox{$Q\to Q/\langle g_j\rangle$}. 
\end{example}

\begin{definition}\label{def:6.9}
An {\em equipolygonal deformation of a curve diagram $(\sC,\sG)$} 
over $A$ is a list of objects in \mbox{$\Defsec_{H_Q\leftarrow Q} (A), Q\in
  \sC$}, and a list of equipolygonal deformations of \mbox{$Q/\langle
  g\rangle$}, \mbox{$Q\in \sC$}, \mbox{$Q/\langle g\rangle\in \sG$}, adapted to
the given deformations of $E_Q,D_Q$, 
such that the following holds: 

If \mbox{$Q\to Q'$} is a formal blow-up of points in $\sC$ and if the
equipolygonal deformation  
of \mbox{$Q/\langle g\rangle$} is defined by \mbox{$G\in Q_A$} and the section
\mbox{$\sigma:A\to Q_A$}, then
\begin{enumerate}
\item the equipolygonal deformation of \mbox{$Q'/\langle
  g'\rangle$} is given by the strict transform $G'$ of $G$ under the formal
blow-up of $I_\sigma$ and a section \mbox{$\sigma':A\to Q'$} which is
compatible to $\sigma$; 
\item the given deformation of $E_{Q'}$ is the 
exceptional divisor of the formal blow-up of $I_\sigma$;
\item if $Q'$ is satellite, then the given deformation of $D_{Q'}$ is the
  strict transform of the given deformation of the exceptional branch at $Q$
  whose strict transform at $Q'$ is $D_{Q'}$. 
\end{enumerate}

%\noindent
For a curve diagram $(\sC,\sG)$, we introduce \mbox{$\bTep=\bigoplus_{Q\in \sC}
  \Tep_{Q,Q/\langle g\rangle}$}. Further, for each equipolygonal deformation of
$(\sC,\sG)$ over 
$A$, we consider the Kodaira-Spencer map \mbox{$T_A\to \bTep$} given
componentwise as above. An equipolygonal deformation of a curve diagram is
called {\em 
  equipolygonal versal\/} (or, {\em ep-versal\/}) if the Kodaira-Spencer map is
surjective.
\end{definition}

\begin{proposition}\label{prop:new9.4}
Let \mbox{$g=g_1\cdot \ldots\cdot g_s$} be the decomposition of \mbox{$g\in Q$}
into tangential components, and assume that \mbox{$G\in Q_A$} defines an
equipolygonal deformation of $Q/\langle g\rangle$, then there exists a unique
factorization \mbox{$G=G_1\cdot \ldots\cdot G_s$} such that the following
hold:
\begin{enumerate}
\item $G_j$ defines an equipolygonal deformation of $g_j$,
  \mbox{$j=1,\dots,s$}. 
\item If $G$ defines an equipolygonal versal deformation $\xi$ of $Q/\langle
  g\rangle$, then $(G_1,\dots,G_s)$ defines an equipolygonal versal deformation
  of the multicurve of tangential components of $Q/\langle g\rangle$.
\end{enumerate}
\end{proposition}

\begin{proof}
  First, assume that \mbox{$\kp=\kp_{Q,R}$} consists of just the point $P$,
  that is, \mbox{$e=d=1$}. In this case, since the coordinates are generic, equipolygonal means 
 nothing but equimultiple. We may assume that \mbox{$u,v\in Q$} are adapted
elements and $g$ is a Weierstra\ss\ polynomial, 
 \mbox{$g\in K[[u]][v]$}, of degree $m$ equal to the multiplicity of $g$. 
 Then, up to a unit,  \mbox{$G\in A[[U]][V]$} is also a Weierstra\ss\
 polynomial of degree $m$. The formal
 transformation \mbox{$(U,V)\to (U,V')$} given by \mbox{$U\mapsto
   U$},  \mbox{$V\mapsto V'U$} leads to \mbox{$g=u^m g'$},
 \mbox{$G=U^m G'$}, with \mbox{$g'\in K[[u]][v']$}, resp. \mbox{$G'\in
   A[[U]][V']$} polynomials of degree $m$ in $v'$ resp. $V'$. On the other hand,
 the factors $g_i$ of $g$ can be assumed to be Weierstra\ss\
 polynomials, too, giving rise to transforms $g'_i$ such that
 \mbox{$g'=g'_1\cdot \ldots\cdot g'_s$}. Now, the residues modulo
 \mbox{$uK[[u,v']]$} of the $g'_i$ are relative prime polynomials in $K[v']$. So,
 Hensel's lemma provides us with a factorization \mbox{$G'=G'_1\cdot
   \ldots\cdot G'_s$}, where \mbox{$G'_i\in A[[U]][V']$} defines a
 deformation of \mbox{$Q'/\langle g'_i\rangle$}. Since the $G'_i$ are
 polynomials, we may apply the backward transformation, giving rise 
to the required 
 factorization \mbox{$G=G_1\cdot \ldots\cdot G_s$}, with $G_i$
 defining an equimultiple deformation of \mbox{$Q/\langle g_i\rangle$} at $P$.
The uniqueness of the factorization follows from the uniqueness of
Hensel's lifting. 

For a general $\kp$, first factorize \mbox{$G=G_1\cdot \ldots\cdot
  G_s$} as equimultiple deformation in terms of tangential components
  according to the previous step.
Now, for those components \mbox{$Q/\langle g_j\rangle$} which are not
tangential to one of the coordinate curves \mbox{$Q/\langle v\rangle$} or
\mbox{$Q/\langle u\rangle$} (with $u,v$ adapted to $H_Q$), the deformation
given by $G_j$ is an equipolygonal deformation as it is equimultiple. If some
component \mbox{$Q/\langle g_j\rangle$} is tangential to  \mbox{$Q/\langle
  v\rangle$} (respectively to \mbox{$Q/\langle u\rangle$}), then the constance
of the Newton polygon for $G$ (condition (2) of Definition \ref{def:6.5alt}) implies the constance of the Newton
polygon $N_i$ for $G_i$, as one can easily deduce from the
factorization of $G$. Thus, each $G_i$ defines an equipolygonal deformation of
$Q/\langle g_j\rangle$, which proves (1).

For (2), we have to prove the surjectivity of the Kodaira-Spencer
map of the multicurve, \mbox{$(\Psi_1,\dots,\Psi_s):T_A\to \bTep$}, 
where  \mbox{$\Psi_j:T_A\to \Tep_{Q,Q/\langle g_j\rangle}$} is the
Kodaira-Spencer map for the equipolygonal deformation $\xi_j$ of
\mbox{$Q/\langle g_j\rangle$} defined by $G_j$. To do so, take elements
\mbox{$a_j\in \Iep_{Q/\langle g_j\rangle}$} and consider
\mbox{$a=\sum_{j=1}^s a_jh_j\in \Iep_{Q/\langle g\rangle}$}, where
\mbox{$h_j=g/g_j$}, 
\mbox{$j=1,\dots,s$}. Since $G$ defines an ep-versal deformation, there exists
some \mbox{$\delta\in T_A$} such that \mbox{$\Psi_\xi(\delta)=[a]\in
  \Tep_{Q,Q/\langle g\rangle}$}. We claim that \mbox{$\Psi_j(\delta)=[a_j]$}.

To show this claim, let \mbox{$\rho\in J_{Q,Q/\langle
    g\rangle}$} be such that the induced deformation $\delta\xi$ be given by
\mbox{$g+\varepsilon (a+\rho)$}. Since \mbox{$g=g_jh_j$}, we may write $\rho$
as \mbox{$\rho=\rho'_jh_j+b_jg_j$} for some \mbox{$\rho'_j\in J_{Q,Q/\langle
    g_j\rangle}$} and \mbox{$b_j\in P$} (by definition of the adapted Jacobian
ideal). On the other hand, if the induced deformations $\delta\xi_j$ are given
by \mbox{$g_j+\varepsilon a'_j$}, where \mbox{$a_j'\in K[[u,v]]$}, then
the equality \mbox{$G=G_1\cdot \ldots G_s$} implies
\mbox{$ a+\rho = \sum_{j=1}^s a'_j h_j$}. 
Together with the above, we get 
\mbox{$ (a'_j-a_j-\rho'_j) h_j \in \langle g_j \rangle \subset  K[[u,v]] $}
(note that $h_i$, \mbox{$i\neq j$}, is divisible by $g_j$). Since $h_j$ has no
common divisor with $g_j$, this shows that \mbox{$a'_j-a_j\in J_{Q,Q/\langle
    g_j\rangle}$} as required. 
\end{proof}

%\noindent
Another important example of curve diagrams are blow-up
diagrams of non-exceptional curves: Let \mbox{$Q$} be an infinitely near point
of $P$ on $R$, and let
\mbox{$Q\to Q/\langle g\rangle$} be a non-exceptional curve. Then the 
{\em blow-up diagram $\sR_{Q,g}=(\sC, \sG)$ of \mbox{$Q\to Q/\langle g\rangle$}} is
defined  
as follows: the entries of the
list  \mbox{$\sC$} are $Q$ and each infinitely near point $Q'$
on $R$ in the first neighbourhood of $Q$ (no repetition). The 
list of non-exceptional curves \mbox{$\sG$} consists of \mbox{$Q\to
  Q/\langle g\rangle$} and its strict transforms \mbox{$Q'\to Q'/\langle
  g'\rangle$} under the formal blow-up \mbox{$Q\to Q'$}.

The following lemma provides us with necessary conditions for an equipolygonal
deformation of \mbox{$Q\to Q/\langle g\rangle$} to lift to an
equipolygonal deformation of the blow-up diagram:

\begin{lemma}\label{lem:6.11}
Let \mbox{$Q/\langle g\rangle$} be unitangential, let \mbox{$U,V\in Q_A$} be
adapted 
elements such that \mbox{$u\equiv U\!\!\mod \fm_A$} is transversal to $g$,
and let \mbox{$G\in A[[U,V]]$} define an equipolygonal deformation of
\mbox{$Q\to Q/\langle g\rangle$} along the trivial section. Let $Q'$
be an infinitely near point on $R$ in the first neighbourhood of $Q$
corresponding to the linear factor \mbox{$v+\alpha u$}, \mbox{$\alpha\in K$},
of the tangent cone of $g$. 
Finally, let $\eta$ be the deformation of $E_{Q'}$ given by the exceptional
divisor of the formal blow-up of \mbox{$I_{\sigma}=\langle U,V\rangle$} and
the section $\sigma'$. Then the following
are equivalent:
\begin{enumerate}
\item[(i)] The strict transform \mbox{$G'\in A[[U,V']]$},
  \mbox{$V'=\frac{V}{U}+\alpha$}, of $G$ defines a deformation which is
  equiadapted to $\eta$. 
\item[(ii)] \mbox{$G=c(V+aU+\alpha U)^m +\,(\text{terms of order $>m$ in
      $U,V$})$}, for some \mbox{$a\in \fm_A$}, \mbox{$c\in 
    A^\ast$}, and \mbox{$I_{\sigma'}=\langle 
    U,V'\!+\!\!\:a\rangle$}. 
\end{enumerate}
Moreover, if (ii) holds, $a$ is uniquely determined unless $A$
is non-reduced and the characteristic of $K$ is a divisor of the multiplicity
$m$ of $g$.  
\end{lemma}

\begin{proof}
 Since $\sigma$ is compatible with $\sigma'$, we have
 \mbox{$I_{\sigma'}=\langle U,V'\!+\!\!\:a\rangle$}. As equiadapted implies
 equimultiple, we may assume that the $m$-jet $L$ of $G$ is a homogeneous
 polynomial of degree $m$. Then the strict transform $G'$ satisfies
 \mbox{$G'(U,V')\equiv L(1,V')\!\!\mod 
   \langle U\rangle$}, hence the induced deformation along $\sigma'$ is 
 equiadapted to $\eta$ iff  \mbox{$L(1,V')=c(V'+a)^m$}. This proves the
 equivalence of (i) and (ii). The uniqueness follows by comparing coefficients.
\end{proof}

%\noindent
For a given \mbox{$P\to R=P/\langle f\rangle$} and any \mbox{$k\geq 0$}, one
has the curve diagram $\bR^{(k)}$ consisting of all infinitely near
points $P'$ of $P$ on $R$ which lie in a neighbourhood of order \mbox{$\nu\leq
  k$} and  such that either $R'$ is singular at $P'$ or $P'$ is not
consecutive to a point $P''$ such that $R''$ is smooth. Notice that $\bR^{(0)}$
consists of the data $(P,R)$ and that $\bR^{(1)}$ is nothing but the blow-up
diagram $\sR_{P,f}$ of \mbox{$P\to P/\langle f\rangle$}. Denote by $h$ the
maximum $k$ such that  \mbox{$\bR^{(k-1)}\neq\bR^{(k)}$}.
%contains all the essential points of $P$ on$R$. 
Set \mbox{$\bR=\bR^{(h)}$}.

The category (resp.\ functor of isomorphism classes) of ep-deformations of $\bR^{(k)}$ is denoted by
\mbox{$\Defep_{\bR^{(k)}}$} (resp.\ \mbox{$\uDefep_{\bR^{(k)}}$}). For
\mbox{$k=h$}, we simply write \mbox{$\Defep_{\bR}$} (resp.\
\mbox{$\uDefep_{\bR}$}) for the category \mbox{$\Defep_{\bR^{(h)}}$} (resp.\
the functor  \mbox{$\uDefep_{\bR^{(h)}}$}).

The following lemma shows that the deformations in  \mbox{$\Defep_{\bR}$} can
also be considered as deformations of the parametrization.

\begin{lemma}\label{lem:6.12}
There is a natural functor \mbox{$\Defep_{\bR} \to \Defsec_{\Rbar \leftarrow P}$}
which identifies \mbox{$\Defep_{\bR}$} with a full subcategory of
\mbox{$\Defsec_{\Rbar \leftarrow P}$}. 
\end{lemma}

\begin{proof}
Each \mbox{$\xi\in \Defep_{\bR}$} is equipped with a section of \mbox{$A\to
  Q_A$} for each $Q$ in $\bR$, in particular, with sections
\mbox{$\sigma_i:Q_{i,A}\to A$} for the maximal points $Q_1,\dots,Q_r$ in $\bR$,
that is, those essential points where $R_Q$ is smooth. After relabelling,
$Q_1,\dots,Q_r$ correspond one to one to the branches $R_1,\dots,R_r$ of $R$,
in the sense that $Q_i$ is on $R_i$ and not on $R_j$ for \mbox{$j\neq i$}. Now,
the deformation \mbox{$Q_{i,A}\to R_{Q_i,A}=\Rbar_{Q_i,A}$} given by $\xi$ becomes a deformation 
\mbox{$Q_{i,A}\to \Rbar_{i,A}$} of $\Rbar_i$, and $\sigma_i$ induces a
section \mbox{$\osigma_i:\Rbar_{i,A}\to A$} as it lies inside $R_{i,A}$. The
$\osigma_i$ all together give rise to \mbox{$\osigma:\Rbar_A\to A$} and the
statement of the lemma follows from this fact.
\end{proof}
\medskip

%\noindent
\begin{center}
{\em Equipolygonal deformations of the parametrization}
\end{center}

Next, we consider equipolygonal deformations in the parametric case. For this,
choose an infinitely near point $Q$ of $P$, a parametrization
\mbox{$\varphi_Q:Q\to \Rbar_Q$} of $R_Q$, and adapted elements $u,v$ to
$H_Q$. Up to permutation we may assume that \mbox{$i(E_Q,Q/\langle
  u\rangle)\geq e$}, \mbox{$i(D_Q,Q/\langle v\rangle)\geq d$} if $E_Q$, $D_Q$ exist. 

For \mbox{$i\in \Lambda_Q$}, denote by \mbox{$m_{e,i}=\ord \varphi_{Q,i}(u)$}
(resp.\ \mbox{$m_{d,i}=\ord \varphi_{Q,i}(v)$}) if $E_Q$ (resp.\ $D_Q$) exists,
and \mbox{$m_{e,i}=m_i$} (resp.\ \mbox{$m_{d,i}=m_i$}) otherwise (see
Definition \ref{def:2.1} for notations). The following
objects depend only on $Q$ and \mbox{$R_Q\to \Rbar_Q$}: 
\begin{enumerate}
\item[(i)] the submodule \mbox{$I^{\ep}_{\Rbar_Q\leftarrow R_Q}=\bigoplus_{i\in
      \Lambda_Q} \bigl(\fmbar^{m_{e,i}}_i\oplus \fmbar^{m_{d,i}}_i\bigr)$} of
  \mbox{$\Rbar_Q\oplus \Rbar_Q$},
\item[(ii)] the adapted Jacobian module given by $$J_{Q,\Rbar\leftarrow
      R}=\fmbar(\dot{u}, \dot{v})+(\langle u,v^d\rangle \oplus \langle
    u^e,v\rangle)\,,$$ 
\item[(iii)] the finite dimensional $K$-vector space 
  \mbox{$T^{\ep}_{Q\Rbar \leftarrow R}\!=\Iep_{\Rbar_Q\leftarrow
      R_Q}/J_{Q,\Rbar\leftarrow R}$}.
\end{enumerate}

\begin{definition}
  Let \mbox{$\eta=(\phi_A,\osigma,\sigma)\in
    \Defsec_{\Hbar_Q\leftarrow Q}(A)$} and \mbox{$\zeta\in
    \Defsec_{\Rbar_Q\leftarrow Q}(A)$}. Then we say that $\zeta$ is 
{\em adapted to $\eta$\/} if \mbox{$\zeta=(\varphi_A,\osigma',\sigma)$} fits
into a commutative diagram  
$$
\UseComputerModernTips
\xymatrix@C=15pt@R=4pt{
\Rbar_{Q,A} \ar@/_/[ddddr]_{\osigma'} && 
{\Hbar}_{Q,A}\ar@/^/[ddddl]^{\osigma}\\   
& Q_A \ar[ur]^-{\phi_A} \ar@{->}_-{\varphi_A}[ul]\ar@/^/[ddd]^(0.45){\!\sigma}\\
 \\
\\
& \,A\,, \ar[uuu]
}
$$
that is, $\zeta$ and $\eta$ are deformations with the same section $\sigma$.
\end{definition}

%\noindent
Choose generic adapted elements \mbox{$u,v\in Q$} and  \mbox{$U,V\in
  I_\sigma\subset Q_A$} adapted to $\eta$ such that \mbox{$u\equiv U\!\!\mod \fm_A$} and 
\mbox{$v\equiv V\!\!\mod \fm_A$}. Assume that $U$ (resp.\ $V$) corresponds to
$E_{Q,A}$ (resp.\ $D_{Q,A}$) if $E_Q$ (resp.\ $D_Q$) exists. Then the adapted
deformation $\zeta$ is called {\em equiadapted to\/} $E_{Q,A}$ (resp.\
$D_{Q,A}$) if for any \mbox{$i\in \Lambda_Q$} one has \mbox{$\ord
  \varphi_{Q,A,i}(U)=\ord \varphi_{Q,i}(u)$} (resp.\ \mbox{$\ord
  \varphi_{Q,A,i}(V)=\ord \varphi_{Q,i}(v)$}).

The adapted deformation $\zeta$ is said to be {\em equipolygonal\/} if it is
equimultiple in the sense of Remark \ref{rmk:2.5}\:(2) and if it is equiadapted
to $E_{Q,A}$ and to $D_{Q,A}$ whenever $E_Q$ and $D_Q$ exist. Notice that if
$\zeta$ is equipolygonal, then its image in $\Defsec_{R_Q}(A)$ is equipolygonal
in the sense of Definition \ref{def:6.5alt}. The converse is not true as Remark
\ref{rmk:2.5}\:(3) shows.

The above definitions of equiadapted and equipolygonal deformations do not
depend on the choice of the (generic) adapted elements $u,v,U,V$.

Denote by \mbox{$\Defep_{\Rbar_Q\leftarrow R_Q}$} (resp.\
\mbox{$\uDefep_{\Rbar_Q\leftarrow R_Q}$}) the category of equipolygonal
deformations (resp.\ of adapted
isomorphism classes of equipolygonal deformations) of the parametrization. The
vector space \mbox{$T^{\ep}_{\Rbar_Q\leftarrow R_Q}$} can be considered as the
tangent space to \mbox{$\uDefep_{\Rbar_Q\leftarrow R_Q}$}. For each
\mbox{$\zeta\in \Defep_{\Rbar_Q\leftarrow R_Q}(A)$}, there is a well-defined linear
(Kodaira-Spencer) map \mbox{$\psi_{\zeta}:T_A\to T^{\ep}_{Q,\Rbar\leftarrow
    R}$} similar as in the non--parametric case. The deformation $\zeta$ is said to be {\em ep-versal\/} if its
Kodaira-Spencer map $\psi_{\zeta}$ is surjective. If $\psi_{\zeta}$ is an
isomorphism, then we call $\zeta$ equipolygonal semiuniversal or {\em ep--semiuniversal}.

\begin{proposition}
Let \mbox{$(a_1,b_1),\dots,(a_s,b_s)\in  \Iep_{\Rbar_Q\leftarrow
    R_Q}$} represent a basis (resp.\ a system of generators) of
\mbox{$T^{\ep}_{Q, \Rbar\leftarrow R}$}. Then

$$ U=u+\sum_{i=1}^s T_ia_i\,,\quad    V=v+\sum_{i=1}^s T_ib_i $$

defines an ep--semiuniversal (resp.\ ep--versal) deformation of
\mbox{$Q\to \Rbar$}. 
\end{proposition}

The proof is trivial, cf. Proposition \ref{prop:7.3}. Moreover, Remark
\ref{rem:7.4} applies mutatis mutandis.

In particular, each ep-versal deformation of $\Rbar_Q\leftarrow R_Q$ has a
smooth basis.

\begin{definition}
A {\em parametric curve diagram\/} $(\sC,\overline{\sG})$ is a curve diagram as
in Definition \ref{def:curve diagram} for which the curve \mbox{$Q\to Q/\langle
  g\rangle$} at each point $Q$ is given by a specified parametrization
\mbox{$\varphi_Q:Q\to \Rbar_Q$}, that is, a diagram consisting of $\sC$ and the
list of parametrizations \mbox{$(\varphi_Q)_{Q\in \sC}$}.
\end{definition}

%\noindent
An {\em equipolygonal deformation\/} of a parametric curve diagram
$(\sC,\overline{\sG})$ is a list of equipolygonal deformations of \mbox{$Q_A\to
  \Rbar_{Q,A}$}, \mbox{$Q\in \sC$}, \mbox{$(\varphi_Q:Q\to \Rbar_Q)\in
  \overline{\sG}$}, adapted to the given deformations of $E_Q,D_Q$ such that
the obvious analogues to (1)\,--\,(3) in Definition \ref{def:6.9} are satisfied.

For a parametric curve diagram $(\sC,\overline{\sG})$ one defines
$$\bTep=\bTep_{(\sC,\overline{\sG})}= \bigoplus_{Q\in \sC} 
  \Tep_{Q,\Rbar\leftarrow R}\,.$$ 
An ep-deformation of $(\sC,\overline{\sG})$ is
said to be {\em ep-versal\/} if the obvious Kodaira-Spencer map \mbox{$T_A\to
  \bTep$} is surjective. 

The parametric analogues to Proposition \ref{prop:new9.4}, Lemma \ref{lem:6.11}
and Lemma \ref{lem:6.12} also hold and they are rather trivial as shown below.

First, the parametric multicurve of tangential components of a given
\mbox{$\varphi:Q\to \Rbar_Q$} is given by \mbox{$\sC=(Q,\dots,Q)$} ($s$
entries, where $s$ is the number of tangential components) and the lists $\overline{\sG}$ of the parametrizations $\varphi_j$ of
the curves given by the branches which share one of the tangents. 

\begin{proposition} Let \mbox{$\zeta\in \Defep_{\Rbar_Q\leftarrow
      R_Q}$}. Denote by $\zeta_{(\sC,\overline{\sG})}$ the deformation of the
  parametric multicurve of tangential components obtained by distributing the
  deformations of the parametrizations of the branches of $\Rbar_Q$ accordings
  to their tangents. Then one has: (1) $\zeta_{(\sC,\overline{\sG})}$ is an
  equipolygonal deformation of $(\sC,\overline{\sG})$, (2) If $\zeta$ is
  ep-versal, then  $\zeta_{(\sC,\overline{\sG})}$ is ep-versal, too. 
\end{proposition}

\begin{proof} (1) follows from the definitions and (2) from the fact that the
  map \mbox{$\Tep_{Q, \Rbar\leftarrow R}\to \bTep_{(\sC,\overline{\sG})}$} is
  obviously surjective.  
\end{proof}

%\noindent
Now, assume that $R_Q$ has only one tangential component and denote by
$\widetilde{\zeta}$ the transform at $Q'_A$ of \mbox{$\zeta\in
  \Defep_{\Rbar_Q\leftarrow R_Q} (A)$} under the blowing up of the section
$\sigma$ given by $\zeta$. One has either \mbox{$m_i=m_{e,i}\leq m_{d,i}$} for all
\mbox{$i\in \Lambda_Q$} or \mbox{$m_i=m_{d,i}\leq m_{e,i}$} for all
\mbox{$i\in \Lambda_Q$}. Without loss of generality, we assume the first
case. Then, since the leading term of $\varphi_{Q,A,i}(U)$ is a unit, one has
that the leading term $a_i$ of $\varphi_{Q,A,i}(V)/\varphi_{Q,A,i}(U)$ is a
well-defined element of $A$. Then the following lemma is again trivial:

\begin{lemma}
Let $R_Q$ be unitangential. Then the deformation $\widetilde{\zeta}$ is equiadapted to $E_{Q',A}$ along the
section $\sigma'$ if and only if one has \mbox{$a_i=a_{i'}$} for every couple
\mbox{$i,i'\in \Lambda_Q$}. In that case, if \mbox{$a=a_{i}$}, \mbox{$i\in
  \Lambda_Q$}, is the common value, then the section $\sigma'$ is given by
\mbox{$\langle u,\frac{v}{u}\!\!\:-\!\!\:a\rangle$}.  
\end{lemma}

%\noindent
Finally, for a given \mbox{$P\to \Rbar$} and \mbox{$k\geq 0$}, denote by
$\overline{\bR}^{(k)}$ the parametric curve diagram consisting of the points $Q$ of
$\bR^{(k)}$ together with the parametrizations $Q\to \Rbar_Q$. The category (resp.\ functor) of ep-deformations of
$\overline{\bR}^{(k)}$ will be denoted by $\Defep_{\overline{\bR}^{(k)}}$
(resp.\ $\uDefep_{\overline{\bR}^{(k)}}$). For \mbox{$k=h$}, we simply
write $\Defep_{\overline{\bR}}$ (resp.\ $\uDefep_{\overline{\bR}}$) for
$\Defep_{\overline{\bR}^{(h)}}$ (resp.\
$\uDefep_{\overline{\bR}^{(h)}}$). Notice that one has natural maps
\mbox{$\Defep_{\overline{\bR}^{(k)}} \to \Defep_{R^{(k)}}$}.

\begin{lemma} \label{lem:7.15} There is a natural functor \mbox{$\Defep_{\overline{\bR}}\to
    \Defsec_{\Rbar\leftarrow P}$} which identifies $\Defep_{\overline{\bR}}$
  with a full subcategory of $\Defsec_{\Rbar\leftarrow P}$.
  
\end{lemma}

\begin{proof} The composite of the functor in Lemma \ref{lem:6.12} with the natural
  functor \mbox{$\Defep_{\overline{\bR}}\to \Defep_{\bR}$} corresponding to
  \mbox{$k=h$} yields the required functor in the lemma.  
\end{proof}

%\noindent
Now, we come to the main result of this section which shows that
ep-deformations of $\bR$, ep-deformations of $\overline{\bR}$ and equisingular
deformations of the parametrization are essentially the same objects.

\begin{theorem}\label{thm:7.16}
For a given parametrization \mbox{$P \to \Rbar$}, the categories
\mbox{$\Defep_{\overline{\bR}}$}, \mbox{$\Defep_{\bR}$} and
\mbox{$\Defes_{\Rbar \leftarrow P}$} (resp.\ the functors
\mbox{$\uDefep_{\bR}$}, \mbox{$\uDefep_{\bR}$} and 
\mbox{$\uDefes_{\Rbar \leftarrow P}$}) are pairwise equivalent (resp. pairwise
isomorphic).
\end{theorem}

\begin{proof}
By lemma
\ref{lem:6.12} and \ref{lem:7.15}
  there are natural maps \mbox{$\Defep_{\overline{\bR}}\to\Defep_{\bR}\to
    \Defsec_{\Rbar\leftarrow P}$} which identify $\Defep_{\overline{\bR}}$
    and $\Defep_{\bR}$ with
    respective full subcategories of \mbox{$\Defsec_{\Rbar\leftarrow P}$}. We claim
    that both subcategories are equal and in fact also equal to the
    subcategory \mbox{$\Defes_{\Rbar\leftarrow P}$} of
      \mbox{$\Defsec_{\Rbar\leftarrow P}$}. The theorem follows from the claim.

To prove the claim, it is enough to show two statements:
\begin{enumerate}
\item[(a)] The image of any \mbox{$\xi\in \Defep_{\bR}(A)$} in
  \mbox{$\Defsec_{\Rbar\leftarrow P}(A)$} is an object of
  \mbox{$\Defes_{\Rbar\leftarrow P}$},
\item[(b)] For any \mbox{$\zeta\in \Defes_{\Rbar\leftarrow P} (A)$} there exists
  \mbox{$\xi\in \Defep_{\overline{\bR}}(A)$} having $\zeta$ as image in
  \mbox{$\Defsec_{R\leftarrow P}(A)$}.
\end{enumerate}

First, for \mbox{$\xi\in \Defep_{\bR}(A)$} consider its image
\mbox{$\zeta\in\Defsec_{\Rbar\leftarrow P}(A)$} and, for each $Q$, the deformation
\mbox{$\xi_Q\in\Defsec_{\Rbar_Q\leftarrow Q}(A)$} given by the data $(\xi, \zeta)$. To
check (a) one needs to show that $\xi_Q$ is equiadapted to $E_{Q, A}$, $D_{Q,
    A}$ ( if $E_Q, D_Q$ exist) and, moreover, an equimultiple
  deformation of the parametrization. We will show this by recurrence on the
  integer $h$. For $h=0$, it is trivial. Now, take $Q$ and $U,V$ adapted to
  $Q$. By recurrence, assume that $\xi_{Q'}$ is equiadapted to \mbox{$E_{Q',
      A}, D_{Q', A}$} (if they exist) and an equimultiple deformation of the
  parametrization for all $Q'$ in the first neigbourhood of $Q$. For fixed
  $Q'$, assume (without loss of generality) that  \mbox{$u=U \mod \frak{m}_A$} is
  transversal to $R_{Q,i}$ for all \mbox{$i\in \Lambda_{Q'}$}. If $Q'$ is
  satellite, then one component of $H_Q$ will be tangent to all $R_{Q, i}$ with
    \mbox{$i\in \Lambda_{Q'}$}, and $V$ is nothing but an equation for that
    component. If $\varphi_A$ denotes the deformation of the parametrization
    given by $\zeta$ then, since $U$ becomes an equation for $E_{Q', A}$, one has
    \mbox{$\ord \varphi_{Q', A, i}(U)=\ord\varphi_{Q', i}(u)$} for \mbox{$i\in
      \Lambda_{Q'}$}. Since \mbox{$V = V'U+aU$} with $a\in A$ and $U,V$
    generators of $I_{\sigma'}$, one has \mbox{$\ord\varphi_{Q,
	A,i}(V)\geq \ord \varphi_{Q,A,i}(U)=\ord\varphi_{A,i}(u)$}. The
    parametric adaptedness to the deformation of the nontangent components of
    branches with \mbox{$i\in \Lambda_{Q'}$}, as well as the parametric
    equimultiplicity for $\xi_Q$ with respect to those branches, follows from
    the above fact. On the other hand, if $Q'$ is satellite, then
    \mbox{$I_{\sigma'}=\langle U, V'\rangle$} with \mbox{$V'=\frac{V}{U}$}. By recurrence
    hypothesis one has \mbox{$\ord\varphi_{Q', A,i}(V)=\ord\varphi_{Q',i}(v')$}, hence
    \mbox{$\ord\varphi_{Q,A,i}(V)=\ord\varphi_{Q,i}(v)$} which shows that the
    parametric equiadaptness condition also holds for the deformation of the
    tangent component of $H_Q$ for $\xi_Q$ with respect to the branches with
    \mbox{$i\in \Lambda_{Q'}$}. This shows (a).

Second, take \mbox{$\zeta\in \Defes_{\Rbar\leftarrow P}(A) $}. Then $\zeta$ induces
deformations \mbox{$\zeta_Q\in \Defsec_{\Rbar_Q\leftarrow Q}(A)$} which are parametric
equimultiple. To check (b) one needs to prove that they are also parametric
equiadapted to \mbox{$E_{Q,A}, D_{Q, A}$}. From this one sees that $\zeta_Q$
gives rise to a list of ep--deformations of the parametrization which defines
a \mbox{$\xi\in \Defes_{\overline{\bR}}(A)$}   whose image in
  \mbox{$\Defsec_{\Rbar\leftarrow R}(A)$} is nothing but $\zeta$. Now, we prove
  this by recurrence on the order of the infinitesimal neighbourhood $\nu$ of
  $Q$. For $\nu=0$, this is trivial. Assume, it is true for $Q$ and take $Q'$ in the first
  neighbourhood of $Q$ and choose adapted $U,V$ at $Q$ such that the equations
  of the components of \mbox{$H_{Q,A}$} are either $U$ or $V$. Also assume
  \mbox{$u=U(\!\mod \frak{m}_A)$} is transversal to $R_{Q,i}$ for \mbox{$i\in\Lambda_Q$}
  so that one has \mbox{$\ord\varphi_{Q,A,i}(U)=\ord\varphi_{Q,i}(u)$}
    ($\varphi_A$ being the data given by $\zeta$). Since $U$ is an equation
    for \mbox{$E_{Q', A}$} it follows that $\zeta_{Q'}$ is equiadapted to
      \mbox{$E_{Q',A}$}. If $Q'$ is satellite, then \mbox{$V'=\frac{V}{U}$} is an
	equation for \mbox{$D_{Q', A}$} and $V$ is an equation for either
	  \mbox{$D_{Q', A}$} or \mbox{$D_{Q, A}$}, By recurrence one has
	      \mbox{$\ord\varphi_{Q,A,i}(V)$}\mbox{$=\ord\varphi_{Q,i}(v)$} hence one deduces \mbox{$\ord\varphi_{Q,A,i}(V')=\ord\varphi_{Q',i}(v')$},
		where \mbox{$v=V(\text{mod},\frak{m}_A),$} \mbox{$v'=V'(\text{mod }\frak{m}_A)$}. This proves the statement (b) and completes the proof of the theorem.
\end{proof}

%\newpage

\section{Proofs for Weakly Equisingular Strata} 
\label{sec:8}

%\noindent
To prove Theorem \ref{theo:6.3}, we need some preparations by means of some
constructions, properties and notations.
 
Let \mbox{$Q$} be an infinitely near point of $P$ on $R$, and let $u,v$ be
generic adapted coordinates (to $H_Q$). Let \mbox{$g\in Q=K[[u,v]]$} define a 
non-exceptional curve, and  
let \mbox{$\eta=\eta_Q\in \Defep_{Q/\langle g\rangle}(C)$} be defined by
\mbox{$G(U,V)\in C[[U,V]]$} and \mbox{$I_{\sigma}=\langle U,V\rangle$}, where
$U,V$ are adapted, \mbox{$u\equiv U\!\!\mod \fm_C$}, \mbox{$v\equiv
  V\!\!\mod \fm_C$}. Then \mbox{$G-g\in \Iep_{R_Q}C[[U,V]]$} and, by
Proposition \ref{prop:new9.4}, we can assume that $G$ decomposes as
\mbox{$G=G_1\cdot \ldots\cdot G_s$}, where $G_j$ defines an equipolygonal
deformation of the tangential component $g_j$ of $g$, \mbox{$j=1,\dots,s$}. Up
to relabeling the components, we may assume that \mbox{$g_1,\dots,g_{s'}$} are
not tangential to $H_Q$, and that \mbox{$g_{s'+1},\dots,g_{s}$} are tangential
to $H_Q$. Notice that \mbox{$0\leq s\!\!\:-\!\!\:s'\leq 2$}. As above, we 
assume additionally that $u$ is transversal to $g_j$, \mbox{$j=1,\dots,s'$}. 
 
We introduce new indeterminates \mbox{$\bW=(W_1,\dots,W_{s'})$} and consider
the following ideals \mbox{$I_j\subset C[[\bW]]$}, \mbox{$j=1,\dots,s$}: 
if \mbox{$j> s'$}, set \mbox{$I_j:=\langle 0\rangle$}.
If \mbox{$j\leq s'$}, the $m_j$-jet of $G_j$ ($m_j$ the multiplicity of
  $G_j$) reads
$$ L_j = c_j \bigl(V^{m_j}+c_{1,j}V^{m_j-1}U+ \ldots +
c_{m_j,j}U^{m_j}\bigr)\,,\quad c_j\in C^\ast, \ c_{ij}\in C\,.$$ 
Write \mbox{$m_j=q_jm'_j$}, such that $q_j$ is the largest power of the
characteristic $p$ dividing $m_j$ (if \mbox{$p=0$}, we set
\mbox{$q_j=1$}). Then we set $I_j$ to be the ideal generated by the following
elements: 
\begin{enumerate}
\item[(i)] $c_{i,j}$, for $i$ not a multiple of $q_j$,
\item[(ii)] $c_{\ell q_j,j} - \dbinom{m'_j}{j} c_{q_j,j}^{\ell}$, for
  \mbox{$\ell=2,\dots,m'_j$}, 
\item[(iii)] \mbox{$W_j^{q_j} + \gamma_{j}^{q_j} - \dfrac{1}{m'_j} c_{q_j,j}$},
    where \mbox{$\gamma_{j}\in K$}, \mbox{$\gamma_{j}^{q_j}\equiv
      \dfrac{1}{m'_j}c_{q_j,j}$ mod $\fm_C$}, 
\item[(iv)] the coefficients $d^{(j)}_{i,k}$ of $U^iV^{\prime k}$,
  \mbox{$(i,k)\not\in N_{Q',R}$}, in the transformed polynomial \mbox{$G'_j= 
U^{-m} G(U,UV'\!-\!\!\:UW_j\!-\!\!\:U\gamma_j)$}, \mbox{$m=\sum_{j=1}^s
m_j$}.  
\end{enumerate}
Set \mbox{$C'=C[[\bW]]/(I_1+\ldots+I_s)$}, and let \mbox{$\gamma:C\to
  C'$} be the natural morphism. Then the induced deformation
\mbox{$\gamma\eta$} is an ep-deformation which can be extended to an
ep-deformation \mbox{$\eta^{(1)}\!\in \Defep_{\sR_{Q,g}}(C')$} of the blow-up
diagram $\sR_{Q,g}$. In fact, this follows by construction: let $Q'$
be the point in the first neighbourhood of $Q$ corresponding to the tangent
direction of $g_j$. Then the vanishing of the polynomials in (i)\,--\,(iii)
guarantees that $L_j$ is a pure power, \mbox{$L_j=c_j
  (V+U\gamma_j+UW_j)^{m_j}$}, and the corresponding strict transform at $Q'$
reads \mbox{$c_j (V'\!+\!\!\:W_j)^{m_j}$}. 
Thus, we may choose as section $\sigma'$
through $Q'$ the section given by \mbox{$\langle U,V'\!+\!\!\:W_j\rangle$} if
\mbox{$j\leq s'$}, respectively the one given by \mbox{$\langle
  U,\frac{V}{U}\rangle$} if \mbox{$j>s'$} and $g_j$ is tangent to $v$,
 respectively  \mbox{$\langle
  \frac{U}{V},V\rangle$} if \mbox{$j>s'$} and $g_j$ is tangent to $u$.
The vanishing of the coefficients in (iv) gives that the deformation is,
indeed, an ep-deformation.

\begin{remark}\label{rmk:8.1}
The above construction has the following universal property: For any morphism
\mbox{$\chi:C\to A$} and any extension of the induced deformation $\chi(\eta)$
to an ep-deformation $\xi$ of $\sR_{Q,g}$ over $A$, there exists a unique map
\mbox{$\chi':C'\!\to A$} such that \mbox{$\chi'\circ \gamma=\chi$} and
\mbox{$\xi=\chi'(\eta^{(1)})$}. In fact, it follows from the property that the
elements in $A$ analogous to the generators (i)\,--\,(iv) of the ideals
\mbox{$I_j\subset C$} are equal to zero, if one takes for $W_j$  the value $a$
in Lemma \ref{lem:6.11}\:(ii) which corresponds to the section of $\xi$ at
$Q'_{j,A}$. For \mbox{$j>s'$}, $Q'_j$ is satellite and, therefore, the section
is given by the intersection of $E_{Q'_{j,A}}$ and $D_{Q'_{j,A}}$, so by
Definition \ref{def:6.9}\:(2) and (3) one has \mbox{$a=0$}. This shows that for
\mbox{$j>s'$} no variables $W_j$ are needed.

Denote by $\sR'_{Q,g}$ the multicurve obtained by deleting $Q$ and $R_Q$ from
the blow-up diagram $\sR_{Q,g}$. Denote by $\eta'$ the ep-deformation of
$\sR'_{Q,g}$ obtained by deleting the deformation associated to $Q$ (that is,
$\gamma \eta$) from the list of deformations of  $\eta^{(1)}$.
\end{remark}

The following proposition is the key step for proving Theorem \ref{theo:6.3}.
 
\begin{proposition}\label{prop:8.2}
With the above assumptions and notations, if \mbox{$\eta\in
  \Defep_{R_Q}(B)$} is an ep-versal deformation of $R_Q$ with smooth base $B$,
then the deformation \mbox{$\eta'\in \Defep_{\sR'_{Q,g}}(B')$} is ep-versal and
its base $B'$ is also smooth.
\end{proposition}

\begin{proof}
 First, assume that \mbox{$s=1$} and set \mbox{$Q'_1=Q'$}, \mbox{$g'_1=g'$},
 \mbox{$G'_1=G'$}. One has two possibilities: (i) \mbox{$s'=0$}, (ii) \mbox{$s'=1$}.

\medskip%\noindent
\mbox{{\it Case (i)\/}}: Take adapted $U,V$ with \mbox{$v\equiv V (\text{mod }
  \fm_B)$} tangent to $g$. Then the formal blow-up \mbox{$Q\to Q'$} (resp.\
\mbox{$Q_B\to Q'_B$}) is given by \mbox{$u=u$}, \mbox{$v=v'u$} (resp.\
\mbox{$U=U$}, \mbox{$V=V'U$}). The blowing up transformation $\Xi$ maps $u^iv^j$ to
$u^{i+j-m}v^{\prime j}$. Moreover, one has \mbox{$\Xi(\Iep_{R_Q})\subset
  \Iep_{R_{Q'}}$}, \mbox{$\Xi(J_{Q,R})\subset J_{Q',R}$} and, hence, $\Xi$
induces a linear map \mbox{$\overline{\Xi}:\Tep_{Q,R}\to \Tep_{Q',R}$}. Since
\mbox{$g'\in J_{Q',R}$} and since the monomial $v^{\prime m}$ occurs in the
support of $g'$, the vector space $\Tep_{Q',R}$ can be generated by monomials
of the image of $\Xi$ (for instance, by the basic monomials with respect to a
monomial ordering for which $v^{\prime m}$ is the leading term of $g'$). It
follows that $\overline{\Xi}$ is surjective.

On  the other hand, one has \mbox{$\psi_{\eta'}=\overline{\Xi}\circ
  \psi_{\eta}$}, $\psi_{\eta}, \psi_{\eta'}$ being the Kodaira-Spencer maps for
$\eta$ and $\eta'$, respectively. Since  $\overline{\Xi}$ is surjective, one
concludes that $\eta'$ is ep-versal and defined over the same base $B$ as
$\eta$. So, \mbox{$B'=B$} is smooth.

\medskip%\noindent
\mbox{{\it Case (ii)\/}:} Take adapted $U,V$ with \mbox{$u\equiv U(\text{mod
  }
  \fm_B)$} transversal to $g$. Since $Q'$ is free and $U$ is an equation for
$E_{Q',B}$ at $Q'_B$, in place of $U,V$ we will consider the couple of
generators $U,Z$ of $I_{\sigma}$ where \mbox{$Z=V+\lambda U$} and
\mbox{$\lambda\in K$} is such that \mbox{$z\equiv Z(\text{mod } \fm_B)$} is tangent
to $g$. If \mbox{$b=\lambda+a$}, one has
\mbox{$L=c(V-aU)^m$} iff \mbox{$L=c(Z-bU)^m$} and, in that case,
\mbox{$I_{\sigma'}=\langle U,V'\!-a\rangle=\langle U,Z'\!-b\rangle$} with
\mbox{$b\in \fm_B$}. Moreover, $L$ needs not to be a pure power. Condition
  (ii) from the beginning of this section 
for (the construction of) $B'$ guarantees that $L$ is a pure power of the above type for the
induced deformation on $B'$.

With the choice of $U,Z$, we set \mbox{$\kp_0=\{Q,Q'\}$}, \mbox{$e_0=1$},
\mbox{$d_0=2$}. This is an auxiliar notation (for the proof) which emphasizes
that $u,z$ can be seen as adapted parameters to $\kp_0$ as both $Q$ and $Q'$
are on the curves $R_Q$ and $Q/\langle z\rangle$. Other couples $u,z_0$ such
that $Q,Q'$ are the only points in common on $R_Q$ and $Q/\langle z\rangle$ can
be thought as generic for $\kp_0$. The Newton polygon $N_0$ of $g$ with respect
to generic couples for $\kp_0$ is constant. Thus, one has the associated
objects as follows: 

(1) the Newton polygon $N_0$, 

(2) the ideal $\Iep_{Q,0}$
generated by the monomials with exponents in \mbox{$N_0+\Z^2_{\geq 0}$}, 

(3)
\mbox{$J_{Q,0}\subset \Iep_{Q,0}$} the ideal generated by  $g,u\frac{\partial  
    g}{\partial u},u^2\frac{\partial
    g}{\partial z}, z\frac{\partial g}{\partial u},z\frac{\partial
    g}{\partial z}$, 

(4) the vector space \mbox{$\Tep_{Q,0}=\Iep_{Q,0}/J_{Q,0}$}.

Notice that one has \mbox{$J_{Q,0}\subset J_{Q,R}$} and
\mbox{$J_{Q,R}=J_{Q,0}+(u\frac{\partial g}{\partial z}) K$} (but not
necessarily \mbox{$u\frac{\partial g}{\partial z}\not\in 
  J_{Q,0}$}). We distinquish three subcases: 
\begin{itemize}
\item[(ii)-1:] \ $p\nmid m$ or p=0\,, 
\item[(ii)-2:] \ $p\mid m$ and \mbox{$J_{Q,R}\subset \Iep_{Q,0}$}\,, 
\item[(ii)-3:] \ $p\mid m$ and \mbox{$J_{Q,R}\not\subset \Iep_{Q,0}$}\,. 
\end{itemize}
Note that \mbox{$J_{Q,0}=J_{Q,R}$} is only possible in case (ii)-2. 

Fix an appropriate monomial ordering (for instance, the weighted degree reverse
lexicographic ordering with respect to the steepest segment of the Newton
polygon $N_{0}$). Let $\{g_1,\dots,g_{\ell}\}$ be a standard basis for the ideal
$J_{Q,0}$ with respect to the fixed ordering. If \mbox{$J_{Q,0}=J_{Q,R}$}
then, of course, $\{g_1,\dots,g_{\ell}\}$ is also a standard basis for the ideal
$J_{Q,R}$; if \mbox{$J_{Q,0}\subsetneq  J_{Q,R}$}, then
$\{g_1,\dots,g_{\ell},u\frac{\partial g}{\partial z}\}$ is a standard basis for
$J_{Q,R}$. Those monomials 
$n_1,\dots, n_d$ in $\Iep_{R_Q}$ which are not in the initial ideal of $J_{Q,R}$
induce a vector space basis of $T^{\ep}_{Q,R}$. Adding $u\frac{\partial
  g}{\partial z}$ (in case 
\mbox{$J_{Q,0}\subsetneq  J_{Q,R}$}), we get a vector space basis for
\mbox{$\Iep_{R_Q}/J_{Q,0}$}. 

The division theorem for standard bases now provides a unique way to write
any element \mbox{$h\in \Iep_{R_Q}$} as 
\[ 
h= \sum_{i=1}^d t_i \cdot n_i + b\cdot u\frac{\partial g}{\partial z} +
\sum_{j=1}^{\ell} c_j(u,z)\cdot 
g_j\,,\quad t_i,b\in K\,, 
\]
where the term \mbox{$b\cdot u\frac{\partial g}{\partial z}$} does not appear if
\mbox{$J_{Q,0}=J_{Q,R}$}, and where the $c_j(u,z)$  are power series with
support in some regions $R_j$ determined by the division procedure. Such
division also applies to $B[[u,z]]$. In fact, if \mbox{$H\in \Iep_{R_Q} B[[u,z]]$},
then $H$ can be uniquely expressed in the form
\[ 
H= \sum_{i=1}^d t_i \cdot n_i + b\cdot u\frac{\partial g}{\partial z} +
\sum_{j=1}^{\ell} C_j(u,z)\cdot g_j\,,\quad t_i,b\in \fm_B\,, 
\]
where \mbox{$C_j(u,z)\in B[[u,z]]$} has support in the region $R_j$.

Now, consider the equation \mbox{$G=g+H$} of the given ep-deformation
$\eta$. One has \mbox{$H\in \Iep_{R_Q} B[[u,z]]$}. By the above
discussion, we have
\[ 
G= g+ \sum_{i=1}^d t_i \cdot n_i + b\cdot u\frac{\partial g}{\partial z} +
\sum_{j=1}^{\ell} C_j(u,z)\cdot g_j\,,\quad t_i,b\in \fm_B\,, 
\]
with \mbox{$C_j(u,z)\in B[[x,y]]$} as above and $b$ only occurring if
\mbox{$J_{Q,0}\subsetneq  J_{Q,R}$}. Since $\eta$ is ep-versal, there are
derivations \mbox{$\delta_1,\dots,\delta_d\in T_B$} such that
\mbox{$\delta_j(t_i)=\delta_{ji}$}. 
This follows from the fact that, for each \mbox{$\delta\in T_B$}, the
Kodaira-Spencer map $\psi_{\eta}$ is surjective and given by
\mbox{$\psi_{\eta}(\delta)=\sum_{i=1}^d \delta(t_i)\cdot n_i$} modulo 
$J_{Q,R}$.  Since \mbox{$T_B=(\fm_B/\fm^2_B)^\ast$}, it follows that 
$t_1,\dots,t_d$ are part of a regular system of parameters for $B$. 

Notice also that, without loss of generality, we can assume that $g$ and $G$ can be
chosen as Weierstra\ss\ polynomials in $z$ of degree $n$.

\smallskip%\noindent
\mbox{\it Case (ii)-1\/}: Here, the leading term of $u\frac{\partial
  g}{\partial z}$ is given by $mbuz^{m-1}$. Thus, $uz^{m-1}$ is not a basic monomial, while
$u^2z^{m-2}, \dots, u^{m-1}z,u^m$ are so. Assume the latter are the basic
monomials $n_2,\dots,n_m$. Then $t_2,\dots,t_m$ are part of a regular system of
parameters, and the leading form of $G$ is given by 
\[
L= z^m +mbuz^{m-1} + t_2 u^2z^{m-2} + \ldots + t_mu^m\,.
\]
Conditions (i) and (ii) for $B'$ impose that \mbox{$L=(z-au)^m$} for some
\mbox{$a\in B$}, which is equivalent to \mbox{$b=-a$} and \mbox{$t_i -
  \binom{m}{i} b^i =0$} for \mbox{$i=2,\dots,m$}. Since \mbox{$b^i\in \fm_B^2$}
for \mbox{$i\geq 2$}, these conditions are analytically independent. Thus, the
quotient ring $\widetilde{B}$ of $B$ modulo the ideal generated by these
\mbox{$m-1$} smooth conditions is again regular. The induced ep-deformation on
it has an equation $\widetilde{G}$ whose leading term is
the pure power \mbox{$(z+bu)^m$}. Now, substituting \mbox{$\widetilde{z}:=
  z+bu$} (which preserves adaptation to the data $\kp_0$, 
hence, does not change the analytical description of the Kodaira-Spencer map),
we get
\begin{align*}
\widetilde{G} (u,\widetilde{z}-bu)\, &=\,  g(u,\widetilde{z}-bu) +
t_1n_1(u,\widetilde{z}-bu) + \sum_{i=m}^d t_in_i(u,\widetilde{z}-bu)
\\ 
&\qquad +  b u\frac{\partial g}{\partial z}(u,\widetilde{z}-bu) +
\sum_{j=1}^{\ell} C_j(u,\widetilde{z}-bu)  
\cdot g_j(u,\widetilde{z}-bu)\\
&=\, g(u,\widetilde{z}) + t_1n_1(u,\widetilde{z}) + \sum_{i=m}^d
t_in_i(u,\widetilde{z}) + \sum_{j=1}^{\ell} C_j(u,\widetilde{z})g_j(u,\widetilde{z})
+ r(u,\widetilde{z})\,, 
\end{align*}
where $C_j$ has coefficients in $\fm_{\widetilde{B}}$ and $r$ has
coefficients in $\fm_{\widetilde{B}}^2$. 
Let $r_i$ be the coefficient of $r$ for the
monomial $n_i$. Then, to insure that one gets
an ep-deformation after blowing up, one has conditions (iv) for $B'$ which
are of type 
\[
t_i+r_i = 0 \quad \text{ for all $i$\, with $n_i \not\in \Iep_{Q,0}$}\,.
\]
All those monomials are basic ones and \mbox{$r_i\in \fm_B^2$}, therefore
the latter conditions are also smooth ones. The ring $B'$ is nothing but the
quotient of $\widetilde{B}$ by the ideal generated by the above conditions. So,
$B'$ is a regular local ring and the transform $\eta'$ of the deformation
$\eta$ is equipolygonal. 

\smallskip%\noindent
\mbox{\it Cases (ii)-2,3\/}: Set \mbox{$m=:q\cdot m'$} with $q$ a power of $p$
and $m'$ prime to $p$. Now, the monomials $uz^{m-1}, \dots, u^{m-1}z,z^m$ are all
basic monomials, and we can assume that they coincide with
$n_1,\dots,n_m$ and, hence, $t_1,\dots,t_m$ are part of a regular system of
parameters for $B$. The leading form \mbox{$z^m+\sum_{i=1}^m t_i u^iz^{m-i}$}
should be a pure power of type \mbox{$(z-au)^m$}, so conditions (ii) and (iii)
for $B'$ can be written in the form
\mbox{$t_i=0$} for \mbox{$i\not\equiv 0$} modulo $q$, \mbox{$t_q=-m'w^q$}
for a new parameter $w$ in a new regular local ring, and
\mbox{$t_jq-(-1)^j\cdot \binom{m'}{j} w^j=0$} for
\mbox{$j=2,\dots,m'$}. Then the quotient ring $\widetilde{B}$ of $B[[w]]$
modulo the ideal generated by the equations given by the above conditions is a
regular local ring. The equation $\widetilde{G}$ of the induced deformation
$\widetilde{\eta}$ of $\eta$ on $\widetilde{B}$ has the pure power
\mbox{$(z-wu)^m$} as leading form. The change \mbox{$\widetilde{z}=z-wu$}
gives rise to 
\[
\widetilde{G} (u,\widetilde{z}+wu) = g(u,\widetilde{z}) + (w+b)u\frac{\partial
  g}{\partial z} + 
\!\sum_{i=m+1}^d \!t_in_i(u,\widetilde{z}) + \sum_{j=1}^{\ell}
C_j(u,\widetilde{z})g_j(u,\widetilde{z}) + r(u,\widetilde{z})\,, 
\]
with $C_j$ having coefficients in $\fm_{\widetilde{B}}$ and \mbox{$r$} having
coefficients in $\fm_{\widetilde{B}}^2$.

In Case (ii)-2, all the monomials \mbox{$n\in \Iep_{R_Q}\setminus \Iep_{Q,0}$} are
basic ones. Conditions (iv) in at the beginning of section \ref{sec:8} are of type 
\[
t_i+r_i = 0 \quad \text{ for all $i>m$  and $n_i \not\in \Iep_{Q,0}$}\,,
\]
where \mbox{$r_i\in \fm_{\widetilde{B}}^2$} is the coefficient of $n_i$ in
$r$. The ring $B'$ is nothing but the quotient of $\widetilde{B}$ modulo the
equations given by the above conditions, so it is a regular local ring.

In Case (ii)-3, all the  monomials \mbox{$n\in \Iep_{R_Q}\setminus \Iep_{Q,0}$}
are basic ones except the leading monomial $n_0$ of $u\frac{\partial
  g}{\partial z}$. Here, conditions (iv) are of type 

\begin{eqnarray*}
&&  t_i+r_i = 0 \quad \text{ for all $i>m$  and $n_i \not\in \Iep_{Q,0}$}\,,\\
&& w+b+r'=0
\end{eqnarray*}

with \mbox{$r_i,r'\in \fm_{\widetilde{B}}^2$}.
These conditions are, again, smooth, since $w$ is analytically independent
from the $t_i$ involved in the other conditions and since the linear part of
\mbox{$w+b+r'$} involves the parameter $w$ with coefficient $1$.  The ring $B'$
is again a regular local ring as it is nothing but the quotient of
$\widetilde{B}$ modulo the equations given by the above conditions.

In all cases, $\eta'$ is the transform of the deformation
$\gamma(\eta)$ induced by $\eta$ on $B'$. Now, we shall prove that
\mbox{$\eta'\in \Defep_{R_{Q'}}(B')$} is ep-versal. For this, we proceed as in
case (i).

\noindent
\rm For it, note that the transformation $\Xi$ associated to the formal blow up
\mbox{$\widetilde{z}'=\frac{\widetilde{z}}{u}$} maps
\mbox{$u^i\widetilde{z}^i$} to \mbox{$u^{i+j-m}\widetilde{z}'^j$} and
satisfies \mbox{$\Xi(\Iep_{Q,0})\not\subset \Iep_{R_{Q'}},
  \Xi(J_{Q,0})\subset J_{Q'R}$}, so that it induces a linear map
\mbox{$\otau:\Tep_{Q,0}\to \Tep_{Q,R}$} which is surjective by the same argument than
in (i). Moreover, as in (i), the Kodaira--Spencer map $ \psi_{\eta'}$ is the
composite of the Kodaira--Spencer map $\psi_0$ of $\gamma\eta$ considered as
ep--deformation with respect to auxiliar data in $\kp_0$ and the linear map
$\Xi$. Thus, the surjectivity of $\psi_{\eta'}$ will be proved if one checks
that $\psi_0$ is surjective. But, this follows from the fact that all the
basic monomials \mbox{$n_i\in \Iep_{Q,0}$} plus the term $\frac{\partial
  g}{\partial z}$ in case (ii)
when \mbox{$J_{Q,0}\not\subseteq J_{Q,R}$ } appear explicitly in the equation
of the deformation \mbox{$\gamma\eta$} with coefficients
\mbox{$t_i=-r_i\in m^2_{B'}$} or \mbox{$w+b=-r'\in m^2_{B'}$} (for
$u \frac{\partial g}{\partial z}$ ). These conditions are part of a regular
system of parameters and therefore, there exist derivatives \mbox{$\delta_i\in T_{B'}$} whose images in
\mbox{$\Tep_{Q,0}$} are equal to \mbox{$n_i\mod J_{Q,0}$}. This shows the
surjectivity of \mbox{$\psi_0:T_{B'}\to \Tep_{Q,0}$}, and, hence, the
ep--versality of $\eta'$ as required.

Now, assume $s>1$. Then, applying proposition \ref{prop:new9.4}, $\eta$ gives rise to an
ep--versal deformation $\eta_T$ of the multicurve of tangential components of
the given curve. Thus, conditions for $B'$ are applied simultaneously to all
tangential components, so that those which correspond to a single component,
are smooth which follows from the case $s=1$. Moreover, they are also smooth all
together as the ep--versality for $\eta_T$ provides stronger hypothesis than the
ep--versality of the individual components. In fact, the ep--versality of
$\eta_T$ guarantees that all coefficients of the basic monomials of all
tangential components form part of a regular system of parameters. This
follows from the fact that the surjectivity of the Kodaira--Spencer map
$\psi_{\eta_T}$ creates independent derivatives for them.

Actually, there are only two types of conditions: those of type
\mbox{$t_i+r_i=0$} with \mbox{$r_i\in m_{B^2}$} (for the coefficient $t_i$ of
some basic monomial); and thos of type \mbox{$w+b+r'=0$} with \mbox{$b\in
  m_B$} and \mbox{$w^q=t_i$} (for $q$ some power of the characteristic and
$t_i$ the coefficient of some basic monomial). In the latter case, the
coefficients of $w$ (in the linear part of the condition) is $1$, whereas the
coefficients of the other $w'$ of similar type which appear are necessarily
$0$. This guarantees that the whole set of conditions for $B'$ are
independent and, therefore, $B'$ is smooth.

The ep--versality of the multicurve deformation $\eta'$ follows from the
independence of the various $t_i$ involved in the construction of $B'$ and the
proof of ep--versality for $s=1$, applied to each component deformation. In
fact, from that proof one deduces that the image of the Kodaira--Spence map
$\psi_{\eta'}$ contains every \mbox{$\Tep_{Q'_{j'}R}$}, so $\psi_{\eta'}$ is
surjective.
\end{proof}

Now we come to the proof of Theorem \ref{theo:6.3}.
\begin{proof}
Take \mbox{$\eta=(C\to R_C, \tau)\in\Defsec_R(C)$}. First, we will construct
the object $\zeta/ \eta$ in the statement.

We may assume that $\eta$ is given by \mbox{$F(x,y)=\sum\limits_{i+j>0}
  c_{ij}x^iy^j\in \frak{m}_C C[[x,y]]$}. Then, consider the ideal $I^{(0)}$ generated
  by the element \mbox{$c_{ij}$} with \mbox{$i+j<m$}, set
  \mbox{$\gamma^{(0)}:C\to C^{(0)}=C/I^{(0)}$} the natural map and
  \mbox{$\eta^{(0)}:=\gamma^{(0)}(\eta)\in\Defep_R(C^{(0)})$}. Now, by \ref{def:7.1} applied to $\eta^{(0)}$ one gets morphismus
  \mbox{$\gamma^{(1)}:C^{(0)}\to C^{(1)}=C'$} and
  \mbox{$\eta^{(1)}\in\Defep_{\bR^{(1)}} (C')$} which extends
  \mbox{$\gamma^{(1)}(\eta^{(0)})$}. By iterating the construction
  simultaneously for all points in  the infinitesimal  neighbourhood of the
  same order one gets
  morphisms \mbox{$\eta^{(k)}\in\Defep_{\bR^{(k)}}(C^{(k)})$} such
  that $\eta^{(k)}$ extends to \mbox{$\gamma^{(k)}(\eta^{(k-1)})$} to points
  in the $k$--th neighbourhood. Then, we set \mbox{$C_\eta=C^{(h)},
  \chi=\gamma^{(h-1)}\circ\cdots\circ\gamma^{(0)}:C\to C_\eta $} and
  \mbox{$\zeta=\eta^{(h)}$}. Theorem \ref{thm:7.16} shows that one in fact has
  \mbox{$\zeta/\eta\in\Defwes_{R,\eta} (C_\eta)$}

Now the functoriality of the construction and (1) in the theorem follow from the construction and from Remark \ref{rmk:8.1} (one
can proceed by recurrence on $k$) taking into account Theorem \ref{thm:7.16}. Statements (2) and (3) follow from the construction
and the successive applications of \ref{rmk:8.1} and $I$ for the construction of
$\zeta/\eta$. Since the extension $C\to C_\eta$ is finite, each condition can
increase the codimension at most by one, allowing to show (4) in this way. If, in
particular, $\eta$ is versal in \mbox{$\Defsec_R$}, then proposition \ref{prop:8.2} shows
  that all the conditions are smooth and, therefore, $C_\eta$ is smooth of
  codimension equal to the number of conditions which contribute for the
  construction of $B'$.

Proposition \ref{prop:8.2} also allows to show that each free point $Q$ contributes with
\mbox{$\frac{1}{2}m_Q(m_Q+1)-1$} such conditions, whereas each satellite
contributes with \mbox{$\frac{1}{2}m_Q(m_Q+1)$}. The sum can be extended to any subset of
infinitely near points on $R$ which contains all satellite points and those
necessary to create them, in particular to \mbox{$\Ess(R)$}
(the minimal subset with the above properties).
\end{proof}
\medskip

Next let us prove Theorem \ref{thm:6.3} in a parallel way. To prove it, we need
some preparations as above.

Let $Q$ be an infinitely near point and fix \mbox{$\zeta=\zeta_Q\in
  \Defep_{\Rbar_Q\leftarrow R_Q}(C)$}. Take adapted $U, V$ and consider the
  tangential components \mbox{$\zeta_{Q,j}\in\Defep_{\Rbar_{Q_j}\leftarrow
  R_{Q_j}}(C), 1\leq j\leq s$}. After relabeling, assume that for
  \mbox{$1\leq j\leq s'$} the $j$--th component is not tangent to an exeptional
  branch and that for \mbox{$s'\leq j\leq s$} it is tangent.

Consider the following ideals \mbox{$I_j\subset C$}. If \mbox{$j>s',
  I_j=(0)$}. If $j\leq s'$, fix a concrete \mbox{$i_j\in\Lambda_{Q'_j}$}. For
  each $i\in\Lambda_{Q'_j}$ denote by $c_i$ the leading term of
  \mbox{$\varphi_{Q,C,i}(V)/\varphi_{Q,C,i}(U)$} and set
\[
\varphi_{Q,C,i}(V)-c_i\varphi_{Q,C,i}(U)=\sum\limits_{l>m_i} c_{il} t_i^l.
\]
Then $I_j$ is the ideal generated by the following elements
\begin{enumerate}
\item [(i)]  \mbox{$c_i-c_{i_j} \text{ for } i\in \Lambda_{Q'_j}, i\neq i_j$}
\item [(ii)] \mbox{$c_{il} \text{ for } i\in\Lambda_{Q'_j} \text{ and }
  m_i=m_{e,i}=m_{d,i}<l<m_i+m'_i,\ m'_i$} being the multiplicity of the strict
  transform of the $i$--th branch at $Q'_j$.
\end{enumerate}

Set \mbox{$C'=C/(I_1+\cdots +I_s)$}. Then, by construction, the natural
morphism \mbox{$\gamma:C\to C'$} induces a deformation of the parametrization
\mbox{$\zeta^{(1)}\in\Defep_{\overline{\sR}_{Q, g}}(C)$} which extends the
deformation $\gamma(\zeta)$ to the parametric blow up diagram \mbox{$\overline{\sR}_{Q, g}(C)$}.

\begin{remark} The above construction has the following universal property:
  For any map $\chi:C\to A$ and an extension of $\chi(\zeta)$ to an
  ep--deformation $\xi$ of $\overline{\sR}_{Q,g}$ there exists a unique map
  $\chi':C'\to A$ such that $\chi'\circ\gamma=\chi$ and
  $\xi=\chi'(\zeta^{(1)})$. In fact, this property follows from the fact that
  the analogous elements to (i) and (ii) in $A$ have to be zero
  because of the existence of $\xi$. Notice also that for $j>s$ no condition
  is required (in concordance with the fact that $I_j=(0)$).

Now, denote by $\zeta'$ the deformation of the multicurve $\overline{\sR}'_{Q_i}$
obtained from the data $(\overline{\sR}_{Q,g}, \zeta^{(1)})$ by deleting $Q$ from the
diagram and the assignation of $\gamma(\zeta)$ at $Q$ to the list of deformations
defining $\zeta^{(1)}$.
\end{remark}

\begin{proposition}\label{prop:8.4} With the assumptions and notations as above,  if \mbox{$\zeta\in \Defep_{\Rbar_Q\leftarrow Q}(B)$} is an ep--versal deformation of
  the parametrization with smooth base $B$, then the deformation
  \mbox{$\zeta'\in \Defep_{\overline{\sR'}_{Q,g}} (B')$} is ep--versal and its base $B'$
  is smooth too.
\end{proposition}

\begin{proof} First assume $s=1$, and consider the two possibilites (i)
  $s'=0$, (ii) $s'=1$.

Case (i).\ Assume, for instance, \mbox{$m_{e,i}<m_{d,i}$} for all $i\in
  \Lambda_Q$ (without loss of generality). One has $e=1<d$ at $Q$ and
  \mbox{$e'=1\leq d'=d-1$} at $Q'$. Then $\sigma'$ is nothing but the
  intersection of $E_{Q',A}$ and $D_{Q',A}$ and $\zeta'$ is equipolygonal with
  $B'=B$. The result follows from the fact that the Kodaira--Spencer map
  $\psi_{\zeta'}$ is nothing but the composite of two surjective maps, namely
  the Kodaira--Spencer map $\psi_\zeta$ and the linear map
\[
\bar{\Xi}:\Tep_{\Rbar_Q\leftarrow R_Q}\to \Tep_{\Rbar_{Q'}\leftarrow R_{Q'}}
\]
induced by \mbox{$\Xi(a_1, \ldots, a_{r'}, b_1, \ldots, b_{r'})=(a_1, \ldots,
a_{r'}, \frac{b_1}{a_1}, \cdots\frac{b_{r'}}{a_{r'}})$} where $r'=\# \Lambda_Q$.

Case (iii).\ One has $d=e=1,\ m_i=m_{i,e}=m_{i, d}$. Take adapted $U, V$
with \mbox{$u=U(\text{mod } \fm_B)$} transversal to $Q\to \Rbar_Q$. Since $Q'$ is free
and $U$ an equation of \mbox{$E_{Q',B}$} at $Q'_B$, one can replace $U,V$
by non adapted \mbox{$U, Z=V+\lambda U, \lambda \in K$}, such that
\mbox{$z=Z(\text{mod } \fm_B)$} is tangent to the branches of $Q\to \Rbar_Q$.

Thus, for \mbox{$i\in\Lambda_Q=\Lambda_{Q'}$} one has \mbox{$c_i=\lambda+
  b_i$} where $b_i\in \fm_B$ is the leading term of \mbox{$\varphi_{Q, B,
  i}(Z)/\varphi_{Q,B,i}(U)$} and
\[
\varphi_{Q,B,i}(Z)-b_i\varphi_{Q,B,i}(U)=\varphi_{Q,B,i}(V)-c_i\varphi_{Q,B,i}(U)=\sum\limits_{i>m_i}c_{il}t_i^l.
\]
One also  has \mbox{$\ord\varphi_{Q,i}(Z)\geq m_i+m'_i$} with equality if
  $m'_i<m_i$.

Consider the vector space \mbox{$\overline{T}^\ep_{Q,0}=\frac{\overline{I}^\ep_{Q, 0}}{\overline{J}^\ep_{Q,0}}$}, where:
\begin{enumerate}
\item [(1)] \mbox{$\overline{I}^\ep_{Q,0}=\bigoplus\limits_{i\in\Lambda_Q}\overline{\fm}_i^{m_i}\oplus \overline{\fm}_i^{m_i+m'_i}\subset
    \Iep_{\Rbar_Q\leftarrow R_Q}$}
\item [(2)] \mbox{$\overline{J}_{Q,0}=\overline{\fm}(\dot{u}, \dot{z})+\fm\oplus \fm_0,
  \text{ with }\fm_0=(u^2,z)\subset \fm=\fm_Q\subset R_Q$}.
\end{enumerate}

Notice that \mbox{$J_{Q,\Rbar\leftarrow R}/\bar{J}_{Q,0}$} is generated by $(0,u)$ as
vector space. Hence, the obvious linear map \mbox{$\Phi:\overline{T}^\ep_{Q,0}\to
  \Tep_{Q,\Rbar\leftarrow R}$} is injective as one has
\mbox{$J_{Q,\Rbar\leftarrow R}\cap
    \Iep_0=J_0$}. 

Now, one has $B'=B/I$ where $I$ is generated by the elements of $\fm_B$ given by
(i) \mbox{$b_i-b_{i_1}, i, i_1\in\Lambda_Q=\Lambda_Q, i\neq i_1, i_1$} fixed,
(ii) \mbox{$c_{il}, i\in\Lambda_Q$} and \mbox{$m_i<l<m_i+m'_i$}. Let $F$ be the
set of indices $(i,l)$ with $i\in\Lambda_Q$ and \mbox{$m_i<l<m_i+m'_i$} or
$l=m_i, i\neq i_1$. For each $(i,l)\in F$, let $z_{i,l}$ be the element in
$\Iep_{\Rbar_Q\leftarrow R_Q}$ whose $j$--th component modulo $J_{Q,R\to\Rbar}$ is $(0,
t_i^{m_i})$. Then, the set consisting of the elements $z_{i,l}$ with $(i,l)\in
  F$ gives a set of linearly independent classes in $\Tep_{Q,\Rbar \to R}$. Since
  the Kodaira--Spencer map \mbox{$T_A\to \Tep_{Q,\Rbar \leftarrow R}$} is surjective,
  there exist derivations \mbox{$\delta_{i,l}\in T_A$} whose images in
  \mbox{$\Tep_{Q,R\to\Rbar}$} are the classes of the corresponding $z_{i,l}$.

Set \mbox{$\varphi_{Q,B,i}(Z)-b_{i_1}\varphi_{Q,B,i}(U)=\sum\limits_{l\geq
    m_i}s_{il}t_i^l$}. If \mbox{$\varphi_{Q,B,i}(U)=\sum\limits_{l\geq
    m_{i_1}}u_{il}t_{i_1}^l$} one has
    \mbox{$s_{i_1m_i}=0,s_{i,m_i}=(b_i-b_{i_1})u_{im_i} \text{ for }i\neq
    i_1, \text{ and } s_{il}=c_{il}+(b_i-b_{i_1})u_{m_{i_1}}$} for
    $l>m_i$. Since $u_{im_i}\not\in \fm_B$ one
    has that the ideal $I$ is also generated by the elements
    $s_{il}\text{ with } (i, l)\in F$.

Now, by construction one has $\delta_{il}(s_{i'l'})=0$ if \mbox{$(i, l)\neq (i',
  l') \text{ and } \delta_{il}(s_{il})=1$}, so the elements $s_{il}\in \fm_B$ are
  linearly independent modulo $\fm^2_B$ and, hence, they are part of a regular
  system of parameters of $B$. This shows that $B'$ is smooth.

It remains to check that $\zeta'$ is ep--versal. For it, notice that one has 
\[
T_{B'}=\{\delta\in T_B\ |\ \delta(s_{il})=0 \text{ for all } (i,l)\in F\}\ .
\]
Hence, one has a induced linear map \mbox{$\Psi_0:T_{B'}\to
  \overline{T}^{ep}_{Q,0}$} which is nothing but the Kodaira--Spencer map for
  $\gamma(\zeta)$ considered as ep--deformation with respect to the data
  $\kp_0$. Since $\Phi \circ \Psi_0$ is equal to the restriction of $\Psi_\zeta$ to
  $T_{B'}$, $\Phi$ is injective and $\Psi_\zeta$ surjective, it follows that
  $\Psi_0$ is surjective. On the other hand, considering $\gamma(\zeta)$
  ep--deformation with respect to $\kp_0$, the situation is exactly the same
  as in case (i), so that, in particular, the analogous map \mbox{$\oXi_0$}
  to \mbox{$\oXi$} is surjective. It follows that $\Psi_{\zeta'}=\oXi_0\circ \Psi_0$ is surjective.

For $s>1$, taking into account \ref{thm:7.16}, the situation is completely analogous
to that of proposition \ref{prop:8.2}. The smoothness of $B'$ and ep--versality of
the multicurve deformation $\zeta'$ follows in the same way. We do not repeat
the arguments here.

%\noindent
Now, the proof of Theorem \ref{thm:6.3} follows in a similar way as for
Theorem \ref{theo:6.3} .
\end{proof}

\begin{proof} (of Theorem 6.3). Take \mbox{$\zeta=(\varphi_C, \sigma,
    \osigma)\in \Defsec_{\Rbar\leftarrow P}(C)$}. First, we will construct the
    ideal \mbox{$I_\zeta$} in the statement.

Assume that $\zeta$ is given by \mbox{$X_i(t)=\sum\limits_{j>0} a_{ij}t^j,
  \ Y_i(t)=\sum\limits_{j>0} b_{ij}t^j,\ i\in\Lambda$} with \mbox{$a_{ij},
  b_{ij}\in \fm_C$}. Then, consider the ideal \mbox{$I^{(0)}$} generated by the
  elements $a_{ij}, b_{ij}$ with \mbox{$0<j<m_i$} and \mbox{$i\in\Lambda$}. Set
  \mbox{$C^{(0)}:$}\mbox{$=C/I^{(0)},\ \gamma^{(0)}$} the natural map and
  \mbox{$\zeta^{(0)}:$}\mbox{$=\gamma^{(0)}(\zeta)\in\Defep_{R\to\Rbar}(C^{(0)})$}. Now,
  the successive application of the construction  for $C'$ (starting from
  \mbox{$\zeta{(0)})$} gives rise to a sequence of ideals \mbox{$I^{(0)}\subset
  I^{(1)}$}\mbox{$\subset\cdots\subset\cdots$}  and deformations
  \mbox{$\zeta^{(k)}\in\Defep_{\overline{\bR}^{(k)}}(C^{(k)}=C/I^{(k)})$}
  such that, for each \mbox{$k$, $\zeta^{(k)}$} extends to the $k$--th order neighbourhood
  points the deformation
  \mbox{$\gamma^{(k)}(\zeta^{(k-1)}),\ \gamma^{(k)}=C^{(k-1)}\to C^{(k)}$} being the
  obvious map. Then, set \mbox{$I_\zeta=I^{(h)}, C_\eta=C/I_\zeta=C^{(h)},
  \pi:C\to C_\zeta$} the natural map, and $\zeta=\pi(\zeta)$. Notice that
  \mbox{$\zeta\in \Defes_{\Rbar\leftarrow P}(C_\zeta)$} as it is the image of
  \mbox{$\zeta^{(h)}\in \Defep_{\overline{\bR}}(C_\zeta)$} by the isomorphism
  in Theorem \ref{thm:7.16}.

Statement (1) follows from the construction, Theorem \ref{thm:7.16} and Remark
8.3. Statements (2) and (3) follow from the construction and successive
applications of Proposition \ref{prop:8.4}. The integer
\mbox{$\con^\es_{\Rbar\leftarrow P}$} is the number of conditions used in the
successive applications to form $I^{(0)}$ for
the construction of $I_\zeta$, so this number bounds the codimension of the
equisingular stratum, showing (2). If $\zeta$ is versal for
\mbox{$\Defsec_{\Rbar\leftarrow P}$}, then the successive applications of
\ref{prop:8.4} show that all those conditions are smooth and transversal, hence, $C_\zeta$
is smooth and of codimension \mbox{$\con^\es_{\Rbar\leftarrow P}$}.

The construction of $I^{(0)}$ and the successive applications of
\ref{prop:8.4} starting from $\zeta$ versal for
\mbox{$\uDefsec_{\Rbar\leftarrow P}$} allow to compute
\mbox{$\con^\es_{\Rbar\leftarrow P}$} giving rise to the formula
\mbox{$\con^\es_{\Rbar\leftarrow P}=\sum\limits_{Q\in \Ess(R)}m_Q-ef_R-(r-1)$}.
\end{proof}

\section{Geometry of Equisingular Strata}

In this section, we study the geometry and show relations among the different equisingularity
strata and objects related to them. In particular, we prove that the dimension
of the weak equisingularity stratum is related to the terms in the
equisingularity exact sequences of Section \ref{sec:exact seq}.

Let \mbox{$\eta_{su}=\eta=(B\to R_B, \tau)\in \Defsec_R(B)$} be a semiuniversal
deformation in $\Defsec_R$. Then $B$ is a regular local ring. For the weakly
equisingular deformation $\zeta/\eta$ of $R$ based in $\eta$, one has that
$B_{\eta}$ is a regular local ring, the map \mbox{$B/I_{\eta}\to B_{\eta}$}
induced by \mbox{$B\to B_{\eta}$} is finite and the weak equisingularity
stratum \mbox{$S^{\wes,\sec}_R=\Spec(B/I^\wes_{\eta})$} is irreducible with \mbox{$\dim
  S^{\wes,\sec}_R=\dim \Spec(B_{\eta})=\dim B_{\eta}$}.

\begin{theorem}\label{thm:8.1}
With the above assumptions and notations from Section \ref{sec:exact seq}, 
$$\dim S^{\wes,\ses}_R=\dim B_{\eta}=\dim T^{1,\es}_R+\dim T^{1,\es}_{\Rbar/R} =\dim
T^{1,\es}_{\Rbar\leftarrow R}+\dim M^{\sec}_R\,.$$  
\end{theorem}

\begin{proof}
From theorems \ref{theo:6.3}  and \ref{thm:6.3} one has 
\[
\dim B_\eta=\tau^{\sec}_R-\frac{1}{2}\sum\limits_{Q\in \Ess(R)}
  m_Q(m_Q+1)+ef_R,
\]  

\[\dim T^{1, \es}_{\Rbar\leftarrow R}=\tau^{1,\sec}_{\Rbar\leftarrow P}-\sum\limits_{Q\in
    \Ess(R)}m_Q+ef_R+(r-1).
\]
On the other hand, from Lemma 5.4 one has
\[
\dim M^{\sec}_R=\tau^{\sec}_R-\tau^{\sec}_{\Rbar\leftarrow
    R}-\delta-(r-1).
\]
Finally, one has the well-known formula
\mbox{$\delta=\frac{1}{2}\sum\limits_{Q\in \Ess(R)}m_Q(m_Q-1)$}. It follows that \mbox{$\dim B_\eta=\dim T^{1, \es}_{\Rbar\leftarrow
    R}+\dim M^{\sec}_R=\dim T^{1, \es}_R+\dim T^{1, \es}_{\Rbar/R}$}, the last equality being a consequence of the exact sequence in Proposition 5.5.
\end{proof}

%\noindent
Since the base spaces of the semiuniversal deformations for \mbox{$\uDefes_R$}
and  \mbox{$\uDefes_{\Rbar\leftarrow R}$} are smooth, the above theorem shows dimensional relations between the different
semiuniversal equisingular deformations and equisingular strata. Next, we will
see how also the geometric nature of these objects can be understood in
terms of the object $\zeta/\eta$.

For this, we first give a technical result. Consider a given diagram of type
$$
\UseComputerModernTips
  \xymatrix@C=7pt@R=4pt{
{\Rbar}\ar@/_2pc/[dddd] && \ar[ll]
{\Rbar}_{A}\ar@/^2.1pc/[dddd]  \\ 
&\; \scriptstyle{\Box}\\
R \ar[uu]  && \ar[ll]
R_{A} \ar[uu]\ar@/^/[dd] && \ar[ll]|\hole R_C
\ar@/^/[dd]^\tau\\ 
&\; \scriptstyle{\Box}&&\; \scriptstyle{\Box} && \\
K \ar[uu]&& \ar[ll] A \ar[uu] && \ar[ll]^-{\chi}
C, \ar[uu] \\ 
}
%\raisebox{-5ex}{$(d)$}
\eqno(d)
$$
that is, a weakly $\wes$ of $R$ based in $(C\to R_C, \tau)$, where $\chi$ is
finite and injective. Denote by $K(\chi)$ the kernel of
the tangent map \mbox{$T(\chi):T_A\to T_C$}. Further denote by
\mbox{$T^{1,\wes,sec}_R\subset T^{1,\sec}_R$} the Zariski tangent space to the
weak equisingularity stratum
$S^{\wes,\sec}_R$. One obviously has \mbox{$\Ties_R \subset
  T^{1,\wes,\sec}_R\subset \T^{1,\sec}_R$}. Note
that, if one fixes \mbox{$(C\to R_C,\tau)$}, then the weakly equisingular
deformation based on it is nothing but a diagram as above which satisfies a
universal property among such diagrams.

\begin{proposition}\label{prop:9.2}
For a given diagram as above, one has an induced commutative diagram of vector
spaces
$$
\UseComputerModernTips
  \xymatrix@C=7pt@R=4pt@M=5pt{
0 \ar[r] & K_0&\!\!\!\!\!\!\!\!:=K(\chi)\cap
K(\beta)\ar[dddd]_-{\varepsilon}\ar@{_{(}->}[ddrr]_-{i}
\ar@/^1pc/@{^{(}->}[ddrrrr]^-{\overline{i}} &&  \\  
\\
 && &&  K(\chi)\ar[dd]_-{\delta} \ar[rr] &&
 T_A\ar[dd]_-{\beta}\ar[ddrr]^-{\gamma} \ar[rr]^-{T(\chi)} && T_C \ar[rrd]^-{\alpha} \\ 
&& && && && && T^{1,\sec}_R\\
0 \ar[rr] && M^{\sec}_R \ar[rr] &&  T^{1,\es}_{\Rbar/R} \ar[rr] &&
T^{1,\es}_{\Rbar\leftarrow R} \ar[rr] && T^{1,\es}_R \ar[rr]\ar@{^{(}->}[urr]^-{I} && 0\,,
}
$$
where all linear maps $\alpha, \beta, \gamma,\delta,\varepsilon$ only depend on
the given diagram; the maps $\alpha,\gamma$ depend only on $\chi(\eta)$ and
$\beta,\delta,\varepsilon$ on the fixed es-deformation of the parametrization
for $\chi(\eta)$. The bottom row is nothing but the exact sequence in
Proposition \ref{5.1} and the other map $I,i,\overline{i}$,
are inclusions. 
\end{proposition}

\begin{proof}
(1) The map $\alpha$ is given by the fact that \mbox{$\eta\in
  \Defsec_R (C)$}.

(2) The map $\beta$ is given by the fact that $\chi(\eta)$ plus the
  chosen deformation of the parametrization is an element of
  \mbox{$\Defes_{\Rbar\leftarrow R}(A)$}.

(3) The map $\gamma$ is the composite of
  \mbox{$T^{1,\es}_{\Rbar\leftarrow R}\to T^{1, \es}_R$} with $\beta$. One has
  \mbox{$I\circ\gamma=\alpha\circ (T(\chi))$}.

(4) Take a vector \mbox{$v\in K(\chi)\subset T_A$}. Then $\delta(v)$
  is defined as follows. Since \mbox{$v: A\to k[\varepsilon]$} is a
  $k$--algebra morphism such that \mbox{$v\circ\chi=0$}, one has that the induced
  deformation on $k[\varepsilon]$ by the base change $v$ is the trivial one,
  therefore, defines an element of \mbox{$T^{1,\es}_{\Rbar/R}$} (since at the
  level of $A$ one has a given deformation of the parametrization too). This
  element of \mbox{$T^{1,\es}_{\Rbar/R}$} is just $\delta(v)$.

(5) If $v\in K_0$, then $\delta(v)$ is an element of
  \mbox{$T^{1,\es}_{\Rbar/R}$} whose image in \mbox{$T^{1,\es}_{R\to\Rbar}$}
  is zero (as \mbox{$v\in \Ker (\beta)$}), so one has \mbox{$\delta(v)\in
  \Msec_R$}. The map $\varepsilon$ is nothing but the restriction of $\delta$
  to $K_0$ taking as target $\Msec_R$ instead of \mbox{$T^{1,\es}_{\Rbar/R}$}.
\end{proof}

\begin{remark}\label{rem:9.3}
Because of the universal property of the weak stratum one has that the image under
$\alpha$ of the subspace $T_{C/\Ker(\chi)}$ of $T_C$ is in $T^{1,\wes,\sec}_R$.  
\end{remark}

%\noindent
Now, consider the particular case of the diagram 

\[
\UseComputerModernTips
  \xymatrix@C=7pt@R=4pt{
{\Rbar}\ar@/_2pc/[dddd] && \ar[ll]
{\Rbar}_{B_{\eta}}\ar@/^2.1pc/[dddd]  \\ 
&\; \scriptstyle{\Box}\\
R \ar[uu]  && \ar[ll]
R_{B_{\eta}} \ar[uu]\ar@/^/[dd] && \ar[ll]|\hole R_B
\ar@/^/[dd]^{\sigma}\\ 
&\; \scriptstyle{\Box}&&\; \scriptstyle{\Box} \\
K \ar[uu]&& \ar[ll] B_{\eta} \ar[uu] && \ar[ll]^-{\varphi}
B \ar[uu]  
}
\eqno(d_{su})
%\raisebox{-5ex}{$(d_{su})$}
\]
given by the weakly equisingular deformation $\zeta/\eta$ based on a
semiuniversal deformation $\eta=(B\to R_B, \sigma)$ in $\Defsec_R$. Since $\eta$ is
semiuniversal, it follows from the universal property of $\zeta/\eta$
that for any diagram $(d)$ there is a diagram map from $(d_{su})$ to $(d)$.
In other words, any weak $\es$--deformation with respect to any given
deformation can be induced (not in a unique way) from $\zeta/\eta$.   Moreover, one has the following particular situation
$$
\UseComputerModernTips
  \xymatrix@C=7pt@R=4pt@M=5pt{
0 \ar[rr] && K_0 \ar[dddd]_-{\varepsilon}^-{\cong}\ar@{_{(}->}[ddrr]_-{i}
\ar@/^1pc/@{^{(}->}[ddrrrr]^-{\overline{i}} &&  \\  
\\
 && &0 \ar[r]&  K(\chi)\ar[dd]_-{\delta}^-{\cong} \ar[rr] &&
 T_{B_{\eta}}\ar@{->>}[dd]_-{\beta}\ar@{->>}[ddrr]^-{\gamma} \ar[rr]^-{T(\chi)} && T_B
 \ar[rrd]^-{\alpha}_-{\cong} \\  
&& && && && && T^{1,\sec}_R\\
0 \ar[rr] && M^{\sec}_R \ar[rr] &&  T^{1,\es}_{\Rbar/R} \ar[rr] &&
T^{1,\es}_{\Rbar\leftarrow R} \ar[rr] && T^{1,\es}_R \ar[rr]\ar@{^{(}->}[urr]^-{I} && 0\,.
}
$$
Hence, one has the following lemma:

\begin{lemma} \label{lem:9.4} With the above assumptions and notations, one has the following
  properties for the object $\zeta/\eta$:
  \begin{enumerate}
  \item[(1)] $\alpha$ is an isomorphism. 
  \item[(2)] $\gamma$ is surjective. 
  \item[(3)] $\delta$ is an isomorphism. 
  \item[(4)] \mbox{$\ker(\beta)\subset \ker(\varphi)$}, so $K_0=\ker (\beta)$.
  \item[(5)] $\varepsilon$ is an isomorphism. 
  \item[(6)] $\beta$ is surjective. 
  \end{enumerate}
\end{lemma}

\begin{proof}
Statement (1) is obvious from the semiuniversality of $\eta$ with respect
to the functor \mbox{$\uDefsec_R$}. Statement (2) (resp. (6)) follows from the
fact that every vector in \mbox{$T^{1\es}_R$}
(resp. \mbox{$T^{1\es}_{\Rbar\leftarrow R}$}) gives rise to a concrete diagram
($d$) with \mbox{$A=k[\varepsilon]$}, so the existing diagram map from
$(d_{su}$) to ($d$) shows that $\gamma$ (resp. $\beta$) is surjective. The
same argument, applied to vectors in \mbox{$T^{1\es}_{\Rbar/R}$} shows that also
$\delta$ is surjective. Since, from Theorem \ref{thm:8.1}, \mbox{$\dim K(\chi)=\dim
  T^{1\es}_{\Rbar/R}$}, (3) follows. Moreover \mbox{$\dim K_0 = \dim \Msec_R$}, so
  (5) follows from (3).
\end{proof}

%\noindent
We come now to the theorem which gives the geometric counterpart of Theorem \ref{thm:8.1}:

\begin{theorem}\label{thm:9.5}
 With the above assumptions and notations, the following holds for the object $\zeta/\eta$: 
 \begin{enumerate}
  \item[(A)] There exists a natural diagram of vector spaces with three exact
    sequences 
$$
\UseComputerModernTips
  \xymatrix@C=7pt@R=10pt@M=4pt{
&& && & T_{B_{\eta}}\ar@/^1pc/@{=>}[rrrd]\ar@{-->}[rd]\\
0 \ar[rr] && M^{\sec}_R \ar[rr]\ar@/^1pc/@{-->}[rrru] &&  T^{1,\es}_{\Rbar/R}
\ar[rr]\ar@{=>}[ru] && 
T^{1,\es}_{\Rbar\leftarrow R} \ar[rr]\ar@{-->}[rd] && T^{1,\es}_R
\ar[rr]\ar@{=>}[rd] && 0\,.\\ 
& 0 \ar@{-->}[ru] && 0 \ar@{=>}[ru]&&&& 0&& 0
}
$$
  \item[(B)] The image of the subspace \mbox{$T^{1,\es}_{\Rbar/R}$} (resp.\
    $M^{\sec}_R$) in $T_{B_{\eta}}$ is a well-defined subspace which represents
    the tangent space to the trivial (resp.\ parametrically trivial)
    subfamilies of $\zeta$.
  \item[(C)] There exist smooth subschemes of $\Spec(B_{\eta})$ such that the
    induced family of $\zeta$ restricted to them gives a semiuniversal
    deformation for $\uDefes_R$. Such subschemes are exactly those smooth ones which are complementary to the image of
    \mbox{$T^{1,\es}_{\Rbar/R}$}, so they may have different tangent
    spaces. Their images in $S^{\wes,\sec}_R$ are also smooth, all share the
    same tangent space, and the induced deformation of $\eta$ on them is
    semiuniversal for $\uDefes_R$. Moreover, \mbox{$S^{\wes,\sec}_R$} is
    nothing but the Zariski closure of the union of those smooth subschemes.
  \item[(D)] There exist smooth subschemes of $\Spec(B_{\eta})$ such that the
    induced family of $\zeta$ restricted to them provides a semiuniversal deformation for
    $\Defes_{\Rbar\leftarrow R}$. These subschemes are exactly all those smooth ones
    which are transversal to
    the image of $M^{\sec}_R$, so they may have different tangent spaces.
  \end{enumerate}
\end{theorem}

\begin{proof}
(A) follows from Lemma 9.4, using \mbox{$\delta^{-1}$}
  (resp. \mbox{$\varepsilon^{-1}$}) to define the linear map
  \mbox{$T^{1,\es}_{\Rbar/R}\to T_{B_\eta}$} (resp. \mbox{$\Msec_R\to
  T_{B_\eta}$}) in the theorem. (B) is obvious. Part of (C) and (D) follow from the
  fact that a semiuniversal deformation $\eta$ (resp. $\zeta$) for the functor
  \mbox{$\uDefes_{R}$} (resp. \mbox{$\uDefes_{\Rbar\leftarrow R}$}) gives
  rise to a diagram\\

%$\hspace*{5.45cm}\eta:$
%\vspace*{-1ex}
\[
\UseComputerModernTips
  \xymatrix{\Rbar & \Rbar_{B^{\es}}\ar[l] \ar@/_2pc/[dd]|(0.5)\hole& \\
R\ar[u] & R_{B^{\es}_R}\ar[l]\ar[u]\ar@{=}[r]\ar@/^1pc/[d] & R_{B^{\es}_R}\ar@/^1pc/[d]\\
K \ar[u] & B^{\es}_R\ar[u]\ar[l]\ar@{=}[r] & B^{\es}_R\ar[u]
}
\eqno{(d^{\es}_R)}
%\qquad\qquad\qquad\raisebox{-5ex}{$D^{\es}_R)$}
\]

resp.\\
\medskip

%$\hspace*{5.2cm}\zeta:$
%\vspace*{-0.5ex}
\[
\UseComputerModernTips
  \xymatrix{
%& \zeta:  \\
\Rbar & \Rbar_{B^{\es}_{\Rbar\leftarrow R}}\ar@/_2pc/[dd]|(0.5)\hole \ar[l] &\\R\ar[u] & R_{B^{\es}_{\Rbar\leftarrow R}}\ar[l]\ar[u]\ar@{=}[r]\ar@/^1pc/[d]
& \Rbar_{B^{\es}_{\Rbar\leftarrow R}}\ar@/^1pc/[d]\\
K \ar[u] & B^{\es}_{\Rbar\leftarrow R}\ar[u]\ar[l]\ar@{=}[r]&
  B^{\es}_{\Rbar\leftarrow R} \ar[u]
}
\eqno{(d^{\es}_{\Rbar\leftarrow R})}
%\qquad\qquad\qquad\qquad\qquad\qquad\qquad\qquad\raisebox{-5ex}{$(D^{\es}_{\Rba%r\leftarrow R})$}
\]
so \mbox{$(d^{\es}_R)\text{ (resp. }(d^{\es}_{\Rbar\leftarrow R}))$} are images
of $(d_{su})$ by  a (not unique) diagram map. This allows to consider $B^{\es}_R$
(resp. $B^{\es}_{\Rbar\leftarrow R}$) as isomorphic to concrete smooth subspaces of
$B_\eta$ satisfying the properties stated in (C) (resp. (D)). For the case of
$\uDefes_R$ the space $B^{\es}_R$ can also be seen as isomorphic to a smooth
subspace of $B^{\sec}_R=B$. The fact that the diagram map is not unique
explains that one
may have several copies of $B^{\es}_R$ (resp.$B^{\es}_{\Rbar\leftarrow R}$)
inside $B_\eta$. All copies may be tangent or transversal one to each others, but they
need to be transversal to the image of \mbox{$T^{1,\es}_{\Rbar/R}$}
(resp. $M^{\sec}_R$) inside $T_{B_\eta}$. However, since \mbox{$T^{1,\es}_R$}
is a well defined vector subspace of \mbox{$T^{1, \sec}_R$}, the image of the
copies of $B^{\es}_R$ inside \mbox{$B^{\sec}_R$} share all  \mbox{$T^{1,\es}_R$} as their tangent space. Notice that, in Remark
4.3.(6) such a copy was considered a candidate for an es--stratum. All these
strata need to be inside the weak equisingular stratum $S^{\wes,\sec}_R$.

The remaining parts of (C) and (D) follow from general arguments applied to this
particular situation. In fact, the deformation $\xi$ on $B_\eta$ which is based
in $\eta$ is versal for both \mbox{$\uDefes_R$} and
\mbox{$\uDefes_{\Rbar\leftarrow R}$} (if one forgets the involved deformation
of $\Rbar$ in the second case). Thus, one has a submersion
\mbox{$\Spec(B_\eta)\to \Spec(B^\es_R)$}
(resp. \mbox{$\Spec(B_\eta)\to\Spec(B^\es_{\Rbar\leftarrow R})$}) such that 

(i) $\xi$ is isomorphic to the pull back of the semiuniversal equisingular
deformation (resp. equisingular deformation of the normalization) on
\mbox{$\Spec(B^\es_R)$} (resp. \mbox{$\Spec(B^\es_{\Rbar\leftarrow R})$}),

(ii) the kernel of the tangent map is the image of
\mbox{$\Ties_{\Rbar\leftarrow R}$} (resp. \mbox{$M^\sec_R$}). 

Then if
$\widetilde{S}$ is a smooth subscheme of \mbox{$\Spec(B_\eta)$} which is
complementary to the image of \mbox{$\Ties_{\Rbar/R}$}
(resp. \mbox{$\Msec_R$}), the deformation induced by $\xi$ on it is
isomorphic, via the induced isomorphism \mbox{$\widetilde{S}\to
  \Spec(B^\es_R)$} (resp. \mbox{$\widetilde{S}\to
  \Spec(B^\es_{\Rbar\leftarrow R})$}), to the pull back of the semiuniversal
equisingular deformation. So the induced deformation of $\xi$ on
$\widetilde{S}$ is itself seminuniversal equisingular. Moreover, if
$\widetilde{S}$ is complementary to the image of \mbox{$\Ties_{\Rbar/R}$} then
the image $S''$ of $\widetilde{S}$ in $S=\Spec(B)$ is also a smooth subscheme
of $S$ isomorphic to $\widetilde{S}$, as the tangent map to
\mbox{$\widetilde{S}\to S$} is injective. It follows that the deformation
  induced by $\eta$ on $S''$ is again semiuniversal for
  \mbox{$\uDefes_R$}. Finally, by construction, it is clear that
  \mbox{$S^{\wes,\sec}_R$} is the Zariski closure of the union of the
  subspaces $S''$ as above. 

This completes the proof of the theorem.
\end{proof}

The following result characterizes when the weak equisingularity stratum and
the strong equisingularity stratum coincide. In fact, it characterizes
numerically when a unique strong equisingularity stratum exists.

\begin{corollary}\label{cor:9.6}
The following conditions for a plane curve singularity are equivalent:
\begin{enumerate}
\item [(i)]    \mbox{$\Ties_{\Rbar/R}=(0)$}
\item [(ii)]   \mbox{$S^{\wes,\sec}_R \text{ and } \Spec(B^\es_R)$}  have the same dimension
\item [(iii)]  \mbox{$ S^{\wes,\sec}_R \text{ and } \Spec(B^\es_R)$}  are isomorphic
\item [(iv)]   \mbox{$\text{ The deformation induced by } \eta \text{ on } S^{\wes.\sec}_R$}
  is strongly equisingular
\item [(v)]    The morphism \mbox{$\Spec(B_\eta)\to S^{\wes,\sec}$} is an isomorphism.
\item [(vi)] There is only one smooth subscheme of \mbox{$S=\Spec(B)$} such
  that the deformation induced by $\eta$ on it is semiuniversal for
  \mbox{$\uDefes_R$}. 
\item [(vii)] There is a unique largest subscheme of $S$ such that the deformation
  induced by $\eta$ on it is stronlgy equisingular.
\end{enumerate}
\end{corollary}

\begin{proof}(i)$\Rightarrow$ (ii) follows from \ref{thm:8.1}, (ii)
  $\Rightarrow$ (iii) and (iii)$\Rightarrow$ (iv) follow from the fact that
  \mbox{$S^\wes_R$} contains subschemes isomorphic to
  \mbox{$\Spec(B^\es_R)$}. Statements (iv) and (v) are equivalent, since in
  both cases the inclusion \mbox{$S^{\wes,\sec}_R\hookrightarrow \Spec(B)$} provides
  a universal object for the category \mbox{$\Defes_{R,\eta}$}. Finally, from
  (iv) follows that \mbox{$\dim  S^\wes_R\leq \dim \Ties_{K, R}$}, so
  \mbox{$\dim_K T^1_{\Rbar/R}=(0)$}. Equivalence of (i) to (v) with (vi) and
  (vii) follows from (C) in Theorem \ref{thm:9.5}.
\end{proof}

\begin{example}\label{exa:9.7}
From above results one can review the cases in example \ref{exa:3.4}. In
particular, in the four cases the subscheme $S'$ can be shown to be nothing but
the weak equisingularity stratum $S^{\wes,\sec}_R$.

\textbf{Case 1:} The embedded resolution consists of the blow up of one point
of multiplicity $2p$ and $2p$ points of multiplicity $1$, with \mbox{$ef_R=2$}, so one
has \mbox{$\dim S^\wes_R = \dim\Ties_R-(2p^2+3p-2)=1+(p-1)(p-2)$}. By
\ref{thm:8.1} and \ref{exa:5.6} one has \mbox{$\dim B^\es_R=(p-1)(p-2)$} and
\mbox{$\dim B^\es_{\Rbar\leftarrow R}=1+(p-1)(p-2)$}. In particular $B_\eta$ is isomorphic
to \mbox{$B^\es_{\Rbar\leftarrow R}$} and each $S_h$ is isomorphic to
\mbox{$\Spec(B^\es_R)$}. The computation of $B_\eta$ gives rise to
\[
\begin{array}{ll}
p\neq 2: & u_{i,j}=0 \text{ for }(i,j)\in D\smallsetminus \{(2p, 0), (p,p)\}\\
         & 2u_{2p,0}-u_{p,p}^2=0\\
p=2:     & u_{i,j}=0 \text{ for } (i,j)\in D\smallsetminus \{(4,0)\}\\
         & u_{4,0}=w_{4,0}^4\\
\end{array}
\]

By eliminating $w_{pp}$ and $w_{40}$ one gets the equations for
$S^{\wes,\sec}_R$ which are exactly those for $S'$ in \ref{exa:3.4} (1). By
specializing $w_{pp}$ (resp. $w_{40}$) to $h$, one gets the
equations for $S_h$. A regular system of parameters for $B_\eta$ consists of
$u_{ij}$ with \mbox{$(i,j)\in D_1\cup D_2$} and \mbox{$w_{pp}$}
(resp. $\omega_{40}$). By specializing, now in \mbox{$\Spec(B_\eta)$}, one gets
$1$--codimensional smooth subschemes $\widetilde{S}_h$ of
$\widetilde{S}^\wes_R$ given by
\[
w_{pp}-h=0\qquad \text{ (resp. } q_{40}-h=0).
\]
Each $\widetilde{S}_h$ applies to $S_h$ in $\Spec(B)$. The deformation induced
by $\eta$ on $\widetilde{S}_h$ is again semiuniversal for \mbox{$\uDefes_R$}. The
image of \mbox{$\Ties_{\Rbar/R}$} in the tangent space to
\mbox{$\Spec(B_\eta)$} is the $1$--dimensional vector subspace $T$ generated
by \mbox{$\frac{\partial}{\partial\omega_{pp}}$}
  (resp. \mbox{$\frac{\partial}{\partial\omega_{40}}$}). Each  $\widetilde{S}_h$
    is transversal to $T$. In fact, any hyperplane transversal of $T$ is
    realized as tangent space to some $\widetilde{S}_h$ for some convenient
    choice of $h$ of type \mbox{$h=\sum\limits_{(i,j)\in D} a_{ij}u_{ij},
      a_{ij}\in K$}.

\textbf{Case 2}: The embedded resolution consists of the blow up of one point
of multiplicity 4, three points of multiplicity 2 and two points of
multiplicity 1, with \mbox{$ef_R=4$}. One has \mbox{$\dim S^\wes_R=\dim \Tisec_R-17=7,
  \dim B^\es_R=\dim B^\es_{\Rbar\leftarrow R}=5$}. In particular, every \mbox{$S_{h,
  h'}$}is isomorphic to \mbox{$(B^\es_R)$}. The computation of \mbox{$B_\eta$} gives
rise to

\[
\begin{array}{lcl}
u_{ij} & = & 0,\ (i,j)\notin\{(4,0), (4,2), (5,1), (6,0), (3,3), (4,3), (5,2),
(5,3)\\
u_{40} & = & w_{4,0}^4\\
u_{51} & = & w_{4,0}^2 u_{3,3}\\
u_{60} & = & w_{6,0}^2+u_{5,1} w_{4,0}+u_{4,2}^2+u_{3,3}w_{4,0}^3\\
(1+w^2_{6,0}) u_{4,2} & = & w_{4,2}^2+w_{6,0}\\
       &   & +(1+w_{6,0}^2) u_{3,3}u_{4,0}^3+(1+w_{6,0})(u_{4,3}w_{4,0}^3+u_{5,2}w_{4,0}^2)
\end{array}
\]
By eliminating $w_{4,0}, w_{6,0}$ and $w_{4,2}$ one gets as equations for
$S^{\wes,\sec}_R$ exactly those of $S$ in \ref{exa:3.4} (2). By specializing
\mbox{$w_{4,0}, w_{6,0},w_{4,2}$} respectively to $h, h'', h'$ one gets the
equations defining $S_{h, h'}$. A regular system of parameters for $B_\eta$ is
given by \mbox{$u_{3,3}, u_{4,2}, u_{4,3}, u_{5,2}, u_{5,3}, w_{4,0},
  w_{4,2}$}. By specializing $w_{4,0}$ and $w_{4,2}$ one gets smooth
subschemes $\widetilde{S}_{h, h'}$ of $S^{\wes,\sec}_R$ given by 
\[
w_{4,0}-h=0\ ,\ w_{4,2}-h'=0,
\]
which apply to $S_{h, h'}$, are isomorphic to \mbox{$\Spec(B^\es_R)$}, and,
moreover, the deformation induced by $\eta$ on them is semiuniversal for
$\uDefes_R$. All the schemes $S_{h, h'}$, are transversal to the image of
\mbox{$\Ties_{\Rbar/R}$} in $T_{B_\eta}$ which is nothing but the vector
subspace generated by \mbox{$\frac{\partial}{\partial w_{4,0}}$} and
  \mbox{$\frac{\partial}{\partial w_{4,2}}$}. The weak equisingularity stratum
    $S^{\wes,\sec}_R$ is singular in this case.

\textbf{Case 3:} The embedded resolution consists of the blow up of one point
of multiplicity \mbox{$p,\ l-1$} points of multiplicity 2, and one point of
multiplicity 1, with \mbox{$ef_R=2$}. So, one has \mbox{$\dim
  S^{\wes,\sec}_R=$} \mbox{$\dim_K\Tisec_R-(2l^2+2l-3)=(l-2)^2,$} \mbox{$\dim(B^\es_R)=\dim
  (B^\es_{\Rbar\leftarrow R})=(l-2)^2$}. The computation of $B_\eta$ gives
rise to 
\[
u_{ij}=0, (i,j)\in D, i+j\leq p \text{ or } i+j=p+1 \text{ and } j\leq l-1,
\]
so one has \mbox{$\Spec(B_\eta) = S^{\wes.\sec}_R=S'$}. On the other hand, the
deformation induced by $\eta$ on $S$ is semiuniversal for $\uDefes_R$, so
$S'=S^{\wes,\sec}_R$ is the strong equisingularity stratum for $R$, which exists in
this case.

\textbf{Case 4:} The~embedded resolution consists of the blow up of one point
of multiplicity $p+1$, one point of multiplicity $2$ and \mbox{$p-1$} points of
multiplicty $1$, with \mbox{$ef_R=2$}. Then one has 
$\dim
  S^{\wes,\sec}_R=\dim_K\Tisec_R-\frac{1}{2}(p^2+5p+2)=\frac{1}{2}(p-1)(p-2),\ \dim (B^\es_R)=\dim (B^\es_{\Rbar\leftarrow R})=\frac{1}{2}(p-1)(p-2)$. The computation of $B_\eta$ gives rise to
\[
u_{i,j}=0, (i,j)\in D, i+j\leq p+1,
\]
so \mbox{$\Spec (B_\eta)=S^{\wes,\sec}_R=S$}, the deformation induced by $\eta$
on $S$ is semiuniversal for $\uDefes_R$, and again \mbox{$S=S^{\wes,\sec}_R$} is the
strong equisingularity stratum, which exists in this case.

\end{example}

\newpage

\end{document}